\documentclass[final]{amsart}

\usepackage{pdfsync}
\usepackage{amsmath,amssymb}
\usepackage{enumerate}
\usepackage{xspace}
\usepackage{boxedminipage}
\usepackage{algorithmic}
\usepackage{listings}
\usepackage{tristan}
\usepackage{color}
\usepackage{subfigure}

\numberwithin{equation}{section}
\usepackage{showkeys}
\setboolean{showtodo}{false}
\setboolean{shownotes}{false}
\setboolean{showchanges}{false}
\setboolean{usemathrsfs}{false}
\author{Jan Giesselmann} 
\address{Jan Giesselmann\newline
Weierstrass Institute\newline
Mohrenstrasse 39\newline
D-10117 Berlin\newline
  Germany} 
\curraddr{}
\email{\linkedemail{jan.giesselmann@wias-berlin.de}} 

\author{Tristan Pryer} 
\address{Tristan Pryer\newline
  Department of Mathematics and Statistics\newline
  Whiteknights\newline
  PO Box 220\newline
  Reading RG6 6AX\newline
  UK
} 
\curraddr{}
\email{\linkedemail{T.Pryer@reading.ac.uk}} 

\thanks{T.P. was supported by the EPSRC grant EP/H024018/1. J.G. was
  supported by the German Research Foundation (DFG) project ``Modeling
  and sharp interface limits of local and non-local generalized
  Navier--Stokes--Korteweg Systems'' and by the EU FP7-REGPOT project
  ``Archimedes Center for Modeling, Analysis and Computation''.}

\title[Energy consistent DG methods for a quasi-incompressible
system]{Energy consistent discontinuous Galerkin methods for a
  quasi-incompressible diffuse two phase flow model}

\date{\today}
\pdfformat{true}
\begin{document}
\renewcommand{\vec}[1]{\geovec{#1}}

\begin{abstract}
  We design consistent discontinuous Galerkin finite element schemes
  for the approximation of a quasi-incompressible two phase flow model
  of Allen--Cahn/Cahn--Hilliard/Navier--Stokes--Korteweg type which
  allows for phase transitions. We show that the scheme is mass
  conservative and monotonically energy dissipative. In this case the
  dissipation is isolated to discrete equivalents of those
  effects already causing dissipation on the continuous level, that
  is, there is no artificial numerical dissipation added
    into the scheme. In this sense the methods are consistent with
  the energy dissipation of the continuous PDE system.
  \\\\
  \smallskip
  \noindent
  \textbf{\keywordsname :} Quasi-incompressibility, Allen--Cahn,
  Cahn--Hilliard, Navier--Stokes--Korteweg, phase transition, energy
  consistent/mimetic, discontinuous Galerkin finite element method.
\end{abstract}


\maketitle

\section{Introduction}

In this work we propose a discontinuous Galerkin (dG) finite element
method for a quasi-incompressible phase transition model of
Allen--Cahn/Cahn--Hilliard/Navier--Stokes--Korteweg type. These
discretisations are of arbitrarily high order in space and provide
\emph{energy consistent} approximations to the model studied. This
means the method is automatically endowed with a particular stability
property by construction.

{Diffuse interface models enjoy the advantage that there is
  only one set of partial differential equation governing the
  behaviour of the mixture over the entire domain. Additionally, no
  particular conditions need be imposed at the interface.}
Historically, the first {diffuse interface} model for a mixture of two
incompressible Newtonian fluids goes back to the so-called \emph{model
  H} proposed in \cite{HH77} where the model is based on the liquids
having the same density.  In \cite{GPV96,LT98} that model was modified
in a thermodynamically consistent way, to allow for liquids with
different densities.  This situation is known as
\emph{quasi-incompressibility}. While the constituents are
incompressible the density of the mixture may vary due to different
concentrations of the constituents.  In this work we will focus on a
model derived in \cite{ADGK} which bears many similarities to
\cite{LT98} while it differs in the choice of the energy functional
and allows for chemical reactions.

The models mentioned above include a \emph{phase field} which
determines which constituent is present at a certain point, for
example, the values $\pm 1$ correspond to the pure constituents.  All
fields (including the phase field) vary smoothly across the interface
between constituents, although steep gradients will usually
occur, hence the name \emph{diffuse interface} model. 

The models derived in \cite{LT98,ADGK} enjoy the advantages of being
thermodynamically consistent, \ie they are compatible with an
entropy function, which may also serve as a Lyapunov function provided
the proper boundary conditions hold, and are frame indifferent.  In
particular, these models are invariant under Galileian transformations
and the only effect of transformations to non-inertial
 coordinate systems is the introduction
of inertial forces, e.g., centrifugal force.  On the other hand they
have the drawback that they include a complicated constraint for the
barycentric (i.e., mass averaged) velocity field, which is no longer solenoidal.
Physically this is to be expected in the presence of exchange of mass
between both constituents. Given two constituents, A and B, if a
certain amount of mass of constituent A becomes constituent B the
different densities and the conservation of mass require a change of
occupied volume.

The divergence constraint makes the extension of (single phase)
incompressible Navier-Stokes solvers infeasible. In addition, the way
the Lagrange multiplier accounting for the incompressibility
constraints enters the equations in \cite{LT98,ADGK} makes the
derivation as well as the numerical analysis of potential schemes
challenging. Regardless, in case of \cite{LT98}, it is possible to
show the model is well-posed, see \cite{Abe09,Abe12}. Although an
extension of these results to \cite{ADGK} does not seem to be
straightforward and to the best of the knowledge of the
authors the well-posedness of \eqref{eq:quasi-crazy-pde-m1} has not
been investigated yet.

The difficulties caused by the divergence constraint have led to the
development of models which are built in such a way that the
considered (not necessarily barycentric) velocity field is solenoidal,
see \cite[e.g.]{AGG,Boy99}, which helps the authors of \cite{Gru,GK}
in the construction and analysis of a scheme.  In particular,
\begin{quote}
  \emph{a simplified version of this model [given in \cite{LT98}] has
    been successfully used for numerical studies \dots In contrast,
    there are -- to the best of the authors' knowledge -- no discrete
    schemes available which are based on the full model \dots This may
    be due to fundamental new difficulties compared with model H \dots
    For instance, the velocity field $\vec v$ is no longer
    divergence-free and therefore no solution concept is available
    which avoids \dots determin[ing] the pressure $p$ \cite{AGG}. }
\end{quote}
In addition,
\begin{quote}
  \emph{Lowengrub and Truskinovsky proposed \dots for the first time a
    diffuse-interface model consistent with thermodynamics. The gross
    velocity field is obtained by mass averaging of individual
    velocities. As a consequence, it is not divergence free, and the
    pressure $p$ enters the model as an essential unknown. However, no
    energy estimates are available to control $p$. Moreover, the
    pressure enters the chemical potential and is hence strongly
    coupled to the phase-field equation. This intricate coupling may
    be one reason why so far it has not been possible to formulate
    numerical schemes for [the] model [given in {\cite{LT98}}] \cite{GK}. }
\end{quote}
{During the review process of this work, a numerical scheme for the model of
Lowengrub--Truskinovsky \cite{LT98} was detailed in \cite{LowengrubetChina}.}

Let us give a short sketch of the derivation of the model in
\cite{ADGK}.  The authors start from the basic balances for mass,
momentum and energy of the mixture.  As an isothermal situation is
considered the latter is only used to determine the heat flux.  The
basic balances contain many quantities (e.g. reaction rates, diffusion
fluxes, stresses) which need to be modelled by constitutive relations.
These are derived by choosing an energy density, introducing a
Lagrange multiplier to account for the incompressibility of the
constituents and exploiting the requirement of thermodynamical
consistency. {Sharp interface limits of the model derived in
  \cite{ADGK} can be found in \cite{ADGK,pfmi}. In particular, the
  authors show that there is mass transfer across the phase boundary,
  hence volume of the phases is not conserved.}

For the derivation of a viable numerical scheme we use a similar
approach to that taken in
\cite{GiesselmannMakridakisPryer:2012}. Here, the authors designed an
approximation of the Navier--Stokes--Korteweg (NSK)/Euler--Korteweg
(EK) system to circumvent some of the numerical artefacts which occur
when applying ``standard'' numerical discretisations to the
problem. The numerical scheme derived was energy consistent in the
sense that for the NSK model it was monotonically energy dissipative
and for the EK model it was energy conservative. The underlying idea
behind the discretisation was to choose a mixed formulation such that
the energy argument at the continuous level could be mimicked at the
discrete level. The quasi-incompressible system we address in this
work has a similar monotone energy functional as the NSK system (see
Theorem \ref{the:energy-cont} and \cite[Lemma
2.3]{GiesselmannMakridakisPryer:2012}).  As such, it becomes possible
to design the numerical scheme to satisfy a discrete equivalent of
this, resulting in a monotonically energy dissipative numerical
scheme, \emph{without} the need for additional artificial dissipation.

Many numerical schemes have been used for the simulation of
{quasi-}incompressible multiphase flows described by \emph{sharp
  interface models}.  In this approach a lot of care is needed to
avoid so called \emph{parasitic currents} in a vicinity of the
interface. They are related to the discretisation of the surface
tension forces, \cite[e.g.]{BKZ92, SZ99, VC00,BGN}. {There
  is also a considerable amount of numerical schemes based on diffuse
  interface models for mixtures of two incompressible fluids with
  differing densities \cite[e.g.]{ALV10,Ding:2007,Dong:2012,Sussman,ZhangTang,LiuShen,ShenYang}
}

We like to point out that our
algorithm does not suffer from parasitic currents,
cf. \S\ref{sec:2d-sims}.

The paper is set out as follows: In \S\ref{sec:notation} we introduce
the quasi-incompressible model and some properties, ultimately leading
to the introduction of the mixed formulation, which is the basis of
designing appropriate numerical schemes. In
\S\ref{sec:spatially-discrete} we detail the construction of a
spatially discrete scheme, moving on to the temporally discrete case
in \S\ref{sec:temporally-discrete}. We combine the results in
\S\ref{sec:fully-discrete} to provide a fully discrete scheme. In
\S\ref{sec:numerics} we conduct various numerical experiments testing
convergence in a simple case as well as the energy consistency in one
and two spatial dimensions and a test on a rotating coordinate system.

\section{Notation and problem setup}
\label{sec:notation}

In this section we formulate the model problem, fix notation and give
some basic assumptions. Let $\W\subset \reals^{d}$, with $d=1,2,3$ be
a bounded domain with Lipschitz boundary. We then begin by
introducing the Sobolev spaces \cite{Ciarlet:1978,Evans:1998}
\begin{equation}
  \sobh{k}(\W) 
  := 
  \ensemble{\phi\in\leb{2}(\W)}
  {\D^{\vec\alpha}\phi\in\leb{2}(\W), \text{ for } \norm{\geovec\alpha}\leq k},
\end{equation}
which are equipped with norms and semi-norms
\begin{gather}
  \Norm{u}_{k}^2
  := 
  \Norm{u}_{\sobh{k}(\W)}^2 
  = 
  \sum_{\norm{\vec \alpha}\leq k}\Norm{\D^{\vec \alpha} u}_{\leb{2}(\W)}^2 
  \\
  \AND \norm{u}_{k}^2 
  :=
  \norm{u}_{\sobh{k}(\W)}^2 
  =
  \sum_{\norm{\vec \alpha} = k}\Norm{\D^{\vec \alpha} u}_{\leb{2}(\W)}^2
\end{gather}
respectively, where $\vec\alpha = \{ \alpha_1,...,\alpha_d\}$ is a
multi-index, $\norm{\vec\alpha} = \sum_{i=1}^d\alpha_i$ and
derivatives $\D^{\vec\alpha}$ are understood in a weak sense. In addition, let
\begin{equation}
  \hoz 
  := 
  \ensemble{\phi\in\sobh1(\W)}
  {\phi\vert_{\partial\W} = 0} 
  \AND
  \hon(\W)
  := 
  \ensemble{\geovec \phi\in\qb{\sobh1(\W)}^d}
  {\Transpose{\qp{\geovec \phi \vert_{\partial\W}}}\geovec n = 0} 
\end{equation}
where $\geovec n$ denotes the outward pointing normal to $\partial
\W$.

We use the convention that for a multivariate function, $u$, the
quantity $\nabla u$ is a column vector consisting of first order
partial derivatives with respect to the spatial coordinates. The
divergence operator, $\div{}$, acts on a vector valued multivariate
function and $\Delta u := \div\qp{\nabla u}$ is the Laplacian
operator. We also note that when the Laplacian acts on a vector valued
multivariate function, it is meant componentwise.  Moreover, for a
vector field $\geovec v$, we denote its Jacobian by $\D \geovec v$. We
also make use of the following notation for time dependant Sobolev
(Bochner) spaces:
\begin{equation}
  \leb{2}(0,T; \sobh{k}(\W))
  :=
  \ensemble{u : [0,T] \to \sobh{k}(\W)}
           {\int_0^T \Norm{u(t)}_{k}^2 \d t < \infty}.
\end{equation}

\subsection{Problem setup}

We consider a mixture of two Newtonian fluids, which might be two
phases of one substance, or two different substances. As both
situations are described by the same model, we will use the terms
phase and constituent interchangeably. In the domain $\W$ we denote
$\phi$ to be the volumetric \emph{phase fraction}, \ie it measures the
fraction of volume occupied by one of the phases. It is scaled in such
a way that $\phi=\pm 1$ corresponds to pure phases. We let $\rho_1 >
0$ and $\rho_2 >0$ be constants that represent the densities of the
incompressible constituents in the fluid. Thus the \emph{total
  density} of the mixture is
\begin{equation}\label{eq:dens}
  \rho(\phi)
  =
  \frac{1}{2} 
  \qb{
    \rho_1
    \qp{
      1+\phi
      }
    +
    \rho_2
    \qp{
      1-\phi
    }
  }.
\end{equation}
We also introduce the constants
\begin{equation}
  \label{eq:cpm}
  c_\pm := \frac {1}{\rho_1} \pm \frac{1}{\rho_2}.
\end{equation}
We let $\gamma > 0$ denote the \emph{capillarity constant} and
$W(\phi)$ be a \emph{double well potential} of $\phi$ then
\begin{equation}
  \begin{split}
    \mu(\phi) &:= W'(\phi) - \gamma \Delta \phi \AND
    \\
    p(\phi) &:= \phi W'(\phi) - W(\phi)
  \end{split}
\end{equation}
represent the \emph{chemical potential} and \emph{pressure}
respectively.  {Note that the thickness of the
  interfacial layer is proportional to $\sqrt{\gamma}$. This can be
  seen by $\Gamma$-limit techniques, cf. \cite{Ste88,ORS90}.} We
denote $\vec v$ to be the \emph{velocity} of the fluid and $\lambda$
is the Lagrange multiplier associated to the incompressibility of the
consitutent{s}.

\subsection{Quasi-incompressible phase transition model}
We then seek $\phi, \vec v \AND \lambda$ such that
\begin{equation}
  \label{eq:quasi-crazy-pde-m1}
  \begin{split}
  \pdt \phi + \div{\phi\vec v} 
  &=
  c_+ \qp{m_j \Delta - m_r}\qp{c_+ \mu(\phi) + c_- \lambda} 
  \\
  \rho(\phi)
  \qp{
    \pdt {\vec v}
    +
    \qp{\Transpose{\vec v} \nabla}\vec v
  }
  +
  \nabla \qp{p(\phi)+\lambda}
  &= 
   \div( \vec \sigma_{NS})
  +
  \gamma \phi \nabla\Delta \phi 
  \\
  \div{\vec v} 
  &=
  c_-
  \qp{m_j \Delta - m_r}\qp{c_+ \mu(\phi) + c_- \lambda} 
  \end{split}
\end{equation}
where
\begin{equation}
  \vec \sigma_{NS} := \eta_1 \div\qp{\vec v} \geomat I_d + \eta_2 \qp{\D \vec v + \Transpose{\qp{\D \vec v}} -\frac{2}{d}\div\qp{\vec v} \geomat I_d },
\end{equation}
is the \emph{Navier--Stokes tensor}, $\geomat I_d$ is the $d\times d$
identity matrix and $\eta_1, \eta_2 \geq 0$
denote bulk and shear viscosity coefficients and $m_j,m_r>0$ are
mobilities. For the derivation of the system
\eqref{eq:quasi-crazy-pde-m1} we refer the reader to \cite{ADGK}.

Note, for clarity of exposition we will not use the full
Navier--Stokes tensor, but the simplified model:
\begin{gather}
  \label{eq:quasi-crazy-pde-1}
  \pdt \phi + \div\qp{\phi\vec v} 
  =
  c_+ \qp{m_j \Delta - m_r}\qp{c_+ \mu(\phi) + c_- \lambda} 
  \\
  \label{eq:quasi-crazy-pde-2}
  \rho(\phi)
  \qp{
    \pdt {\vec v}
    +
    \qp{\Transpose{\vec v} \nabla}\vec v
  }
  +
  \nabla \qp{p(\phi)+\lambda}
  = 
  \eta \Delta \vec v
  +
  \gamma \phi \nabla\Delta \phi
  \\
  \label{eq:quasi-crazy-pde-3}
  \div\qp{\vec v} 
  =
  c_-
  \qp{m_j \Delta - m_r}\qp{c_+ \mu(\phi) + c_- \lambda} ,
\end{gather}
with $\eta >0$.
An energy consistent discretisation of the full model follows our
arguments given a standard (signed) discretisation of the
Navier--Stokes tensor and numerical experiments to this end are
  given in \S\ref{sec:rotate}.

\begin{Rem}[local conservation of mass]
  \label{rem:cons-of-mass-loc}
  It is important to observe that combining
  \eqref{eq:quasi-crazy-pde-1} and \eqref{eq:quasi-crazy-pde-3} gives
  \begin{equation}
    \label{eq:massj}
    \frac{c_-}{c_+} \qp{\pdt \phi + \div\qp{\phi\vec v}} -  \div{\vec v}  =0.
  \end{equation}
  Due to \eqref{eq:dens} and \eqref{eq:cpm} this is equivalent to
  \begin{equation}
    \label{eq:mass} 
    \pdt \rho(\phi)  +  \div(\rho(\phi)\vec v)=0,
  \end{equation}
  \ie the (local) conservation of mass is encoded in
  (\ref{eq:quasi-crazy-pde-1})--(\ref{eq:quasi-crazy-pde-3}).
\end{Rem}

\begin{Rem}[boundary conditions]
  \label{rem:bcs}
  We associate with (\ref{eq:quasi-crazy-pde-1})--(\ref{eq:quasi-crazy-pde-3}) the following boundary conditions:
  \begin{gather}
    \nabla \phi \cdot \vec n = 0
    \\
    \label{eq:bcs-vel}
    \vec v = \vec 0
    \\
    \qp{\nabla\qp{c_+ \mu(\phi) + c_- \lambda}} \cdot \vec n = 0.
  \end{gather}
  This choice yields \emph{global conservation of mass},  \emph{global momentum balance} and a \emph{entropy dissipation equality}
as we will see subsequently.
\end{Rem}
\begin{Pro}[Conservation of mass,balance of momentum]
  Let $\qp{\phi,\vec v,\lambda}$ be a strong solution to the system
  (\ref{eq:quasi-crazy-pde-1})--(\ref{eq:quasi-crazy-pde-3}) satisfying the boundary conditions in Remark \ref{rem:bcs} then
  \begin{equation}
    \label{eq:cons-of-mass-glob}
    \ddt\qp{ \int_\Omega \rho(\phi)} =0,
  \end{equation}
  and
  \begin{equation}
    \label{eq:cons-of-mom}
    \ddt \qp{\int_\Omega \rho(\phi)\vec v}
    =
    - \int_{\partial \Omega} 
    \qp{p(\phi) + \lambda - \phi \Delta \phi}\vec n
    - \qp{D\vec v}
    \cdot \vec n .
  \end{equation}
\end{Pro}
\begin{Proof}    
  The proof of (\ref{eq:cons-of-mass-glob}) can be seen using Remark
  \ref{rem:cons-of-mass-loc} and the boundary conditions
  (\ref{eq:bcs-vel}). To see (\ref{eq:cons-of-mom}) it is enough to
  use (\ref{eq:quasi-crazy-pde-2}), the identity
  \begin{equation}
    \phi \nabla \Delta \phi 
    =
    \div\qp{\qp{\phi \Delta \phi 
        +
        \frac{1}{2} \abs{\nabla \phi}^2 }\geomat{I}_d 
      -
      \nabla \phi \otimes \nabla \phi},
  \end{equation}
  and the boundary conditions.
  \end{Proof}
  
  {For completeness we formulate the energy dissipation equality in
    Theorem \ref{the:energy-cont}. Its validity is a direct
    consequence of the modeling paradigm employed in \cite{ADGK} and a
    proof can be found in \cite{pfmi}. We have organized the proof in
    such a way that it may serve as a guideline for the construction
    of a numerical discretisation which satisfies a discrete energy
    dissipation equality.}

\begin{The}[energy dissipation equality]
  \label{the:energy-cont}
  Let $\qp{\phi,\vec v,\lambda}$ be a strong solution to the system
  (\ref{eq:quasi-crazy-pde-1})--(\ref{eq:quasi-crazy-pde-3}) satisfying the boundary conditions in Remark \ref{rem:bcs}, then
  \begin{equation}
    \begin{split}
      \ddt 
      \qp{
        \int_\W
        W(\phi)
        +
        \frac{\rho(\phi)}{2} \norm{\vec v}^2
        +
        \frac{\gamma}{2} \norm{\nabla \phi}^2 
        }
      =
      -
      \int_\W
      m_j
      &
      \norm{\nabla \qp{c_+ \mu(\phi) + c_- \lambda}}^2
      \\
      &+
      m_r 
      \qp{
        c_+ \mu(\phi) + c_- \lambda
      }^2
      +
      \eta\norm{\D \vec v}^2.
    \end{split}
  \end{equation}
\end{The}
\begin{Proof}
  Let
  \begin{equation}
    \begin{split}
      a 
      &=
      c_+ 
      \mu(\phi)
      +
      c_- \lambda
      \AND
      \\
      b
      &= 
      \lambda 
      +
      \frac{\rho_1 + \rho_2}
      {4}
      \norm{\vec v}^2.
    \end{split}
  \end{equation}
  We proceed by testing (\ref{eq:quasi-crazy-pde-1}) with
  $\frac{a}{c_+}$ and (\ref{eq:quasi-crazy-pde-2}) with $\vec v$ and
  taking the sum, yielding
  \begin{equation}
    \begin{split}
      0 
      &=
      \int_\W
      \frac{a\pdt \phi}
           {c_+}
      +
      \frac{a\div\qp{\phi\vec v}}
           {c_+}
      -
      m_j a \Delta a 
      +
      m_r a^2
      +
      \rho(\phi)
      \bigg(
        \pdt{\vec v} \cdot \vec v
        +
        \qp{\qp{\vec v \cdot \nabla}\vec v }
        \cdot \vec v
        \\
        &
        \qquad\qquad 
        -
        \frac{1}{2}
        \nabla\qp{\norm{\vec v}^2}\cdot \vec v
        \bigg)
      +
      \nabla b \cdot \vec v
      +
      \frac{\phi}{c_+}
      \nabla
      \qp{a - c_- b} \cdot \vec v
      -
      \eta \vec v \cdot \Delta \vec v.
   \end{split}
  \end{equation}
  Integrating by parts and noting that
  \begin{equation}
    \qp{\qp{\vec v \cdot \nabla}\vec v }
    \cdot \vec v
    -
    \frac{1}{2}
    \nabla\qp{\norm{\vec v}^2} \cdot \vec v
    =
    0 
  \end{equation}
  gives
  \begin{equation}
    \begin{split}
      0
      &=
      \int_\W
      \frac{a\pdt \phi }
           {c_+}
      +
      \frac{a\div\qp{\phi\vec v}}
           {c_+}
      +
      m_j \norm{\nabla a}^2
      +
      m_r a^2
      +
      \rho(\phi)
      \pdt{\vec v} \cdot \vec v
      +
      \nabla b \cdot \vec v 
      \\
      &\qquad\qquad +
      \frac{\phi}{c_+}
      \nabla \qp{a - c_- b} \cdot \vec v
      +
      \eta \norm{\D \vec v}^2     
      - 
      \int_{\partial\W}
      m_j a\nabla a\cdot \vec n
      +
      \eta \qp{\D \vec v \cdot \vec n}\cdot \vec v.
    \end{split}
  \end{equation}
  Due to the boundary conditions given in Remark \ref{rem:bcs} the
  boundary terms are zero. In addition we note that
\begin{equation}
  \int_\W
  \frac{a\div\qp{\phi\vec v}}{c_+}
  +
  \frac{\phi\nabla a \cdot \vec v}{c_+}
  =
  \int_\W
  \frac{\div\qp{a\phi\vec v}}{c_+}
  =
  \int_{\partial\W} 
  \frac{a\phi\vec v \cdot \vec n}{c_+}
  =
  0
\end{equation}
again due to the boundary conditions, leaving
\begin{equation}
  \begin{split}
    0
    &=
    \int_\W
    \frac{a\pdt \phi }
         {c_+}
    +
    m_j \norm{\nabla a}^2
    +
    m_r a^2
    +
    \rho(\phi)
    \pdt{\vec v} \cdot \vec v
    +
    \nabla b \cdot \vec v 
    -
    \frac{c_- \phi}{c_+}
    \nabla b \cdot \vec v
    +
    \eta \norm{\D \vec v}^2     
    .
  \end{split}
\end{equation}
Using the definition of $a$ in the first term and integrating by parts
the two terms involving $b$, we see
\begin{equation}
  \begin{split}
    0
    &=
    \int_\W
    W'(\phi)\pdt \phi
    -
    \gamma 
    \pdt \phi
    \Delta \phi
    + 
    \frac{c_-\lambda\pdt \phi}{c_+}
    +
    m_j \norm{\nabla a}^2
    +
    m_r a^2
    +
    \rho(\phi)
    \pdt{\vec v} \cdot \vec v
    \\
    &\qquad\qquad 
    -
    b \div\qp{\vec v}
    +
    \frac{c_- b}{c_+}
    \div\qp{ \phi\vec v }
    +
    \eta \norm{\D \vec v}^2     
    +
    \int_{\partial\W}
    b\vec v \cdot \vec n
    - 
    \frac{c_- b\phi}{c_+}
    \vec v \cdot \vec n
    .
  \end{split}  
\end{equation}
The boundary terms vanish, again, due to Remark \ref{rem:bcs}. Using
the local conservation of mass \eqref{eq:massj}
\begin{equation}
  \begin{split}
    0
    &=
    \int_\W
    W'(\phi)\pdt \phi
    -
    \gamma 
    \pdt \phi
    \Delta \phi
    + 
    \frac{c_-\lambda\pdt \phi}{c_+}
    +
    m_j \norm{\nabla a}^2
    +
    m_r a^2
    +
    \rho(\phi)
    \pdt{\vec v} \cdot \vec v
    \\
    &\qquad\qquad 
    -
    \frac{c_-}{c_+}
    b \pdt \phi
    +
    \eta \norm{\D \vec v}^2     
    .
  \end{split}  
\end{equation}
Using the definition of $b$ and integrating the second term by parts,
it holds that
\begin{equation}
  \label{eq:cont-energy-penultimate-bound}
  \begin{split}
    0
    &=
    \int_\W
    W'(\phi)\pdt \phi
    +
    \gamma 
    \nabla \qp{\pdt \phi}
    \nabla \phi
    + 
    m_j \norm{\nabla a}^2
    +
    m_r a^2
    +
    \rho(\phi)
    \pdt{\vec v} \cdot \vec v
    \\
    &\qquad\qquad
    -
    \frac{c_-\qp{\rho_1+\rho_2}}{4c_+}
    \norm{\vec v}^2 \pdt\phi
    +
    \eta \norm{\D \vec v}^2     
    -
    \int_{\partial\W}
    \gamma \pdt \phi \nabla \phi \cdot \vec n
    .
  \end{split}  
\end{equation}
Due to the definition of $c_\pm$ (\ref{eq:cpm})
\begin{equation}
  \frac{c_-\qp{\rho_1 + \rho_2}}{4c_+} 
  =
  \frac{\rho_2 - \rho_1}{4}
  =
  - \frac{d\rho(\phi)}{d\phi},
\end{equation}
and hence
\begin{equation}
  \begin{split}
    \int_\W
    \ddt 
    \qp{
      W(\phi) 
      +
      \frac{\rho(\phi)}{2}
      \norm{\vec v}^2
      +
      \frac{\gamma}{2} \norm{\nabla \phi}^2
    }
    &=
    \int_\W
    W'(\phi)\pdt \phi
    +
    \rho(\phi) \pdt {\vec v}\cdot \vec v
    \\
    &\qquad+
    \frac{c_-\qp{\rho_1 + \rho_2}}{4c_+} 
    \norm{\vec v}^2 \pdt \phi
    +
    \gamma \nabla\qp{\pdt{\phi}} \cdot \nabla \phi.
  \end{split}
\end{equation}
Using the boundary conditions in Remark \ref{rem:bcs} one final time
to eliminate the boundary contributions from
(\ref{eq:cont-energy-penultimate-bound}) shows
\begin{equation}
  0 
  =
  \int_\W
  \ddt 
  \qp{
    W(\phi) 
    +
    \frac{1}{2}
    \rho(\phi)
    \norm{\vec v}^2
    +
    \frac{\gamma}{2} \norm{\nabla \phi}^2
  }
  +
  m_j \norm{\nabla a}^2
  +
  m_r a^2
  +
  \eta \norm{\D \vec v}^2.  
\end{equation}
The result then follows using the definition of $a$, concluding the proof.
\end{Proof}

\subsection{Continuous mixed formulation}

The proof of Theorem \ref{the:energy-cont} motivates the introduction
of the auxiliary variables $a,b, \geovec{q}$, transforming
(\ref{eq:quasi-crazy-pde-1})--(\ref{eq:quasi-crazy-pde-3}) into the
following mixed system:
\begin{equation}
  \begin{split}
    0
    &= 
    \pdt \phi 
    +
    \div\qp{\phi\geovec v}
    -
    c_+ m_j \Delta a
    +
    c_+ m_r a
    \\
    0 
    &= 
    \rho(\phi)
    \qp{
      \pdt{} \vec{v} 
      +
      \qp{{\vec{v}}\cdot \nabla}\vec{v}
      -
      \frac{1}{2}
      \nabla\qp{\norm{\vec{v}}^2}
    }
    - 
    \eta \Delta \geovec{v}
    +
    \nabla b
    +
    \frac{\phi}{c_+}\nabla( a - {c_-} b)
    \\
    0 
    &= 
    \div\qp{\geovec v}
    -
    \frac{c_-}{c_+}
    \qp{
      \pdt{} \phi
      +
      \div\qp{\phi\geovec v}
    }
    \\
    0
    &=
    a
    - 
    c_+W'(\phi)
    +
    c_+\gamma\div\qp{\vec{q}}
    -
    c_-\lambda
    \\
    0&= 
    b
    -
    \lambda
    -
    \frac{\rho_1+\rho_2}{4}
    \norm{\geovec{v}}^2
    \\
    0 
    &=
    \geovec{q} 
    -
    \nabla \phi,
  \end{split}
\end{equation}
coupled with boundary conditions
\begin{gather}
  \label{eq:mixed-bcs}
  \vec q \cdot \vec n = 0, \quad
  \vec v = \vec 0, \quad
  \nabla a \cdot \vec n = 0.
\end{gather}

\section{Spatially discrete approximation}
\label{sec:spatially-discrete}

In this section we design spatially discrete approximations of the
system (\ref{eq:quasi-crazy-pde-1})--(\ref{eq:quasi-crazy-pde-3}) of
arbitrary order using discontinuous Galerkin finite elements.

Let $\T{}$ be a conforming, shape regular triangulation of $\W$,
namely, $\T{}$ is a finite family of sets such that
\begin{enumerate}
\item $K\in\T{}$ implies $K$ is an open simplex (segment for $d=1$,
  triangle for $d=2$, tetrahedron for $d=3$),
\item for any $K,J\in\T{}$ we have that $\closure K\meet\closure J$ is
  a full subsimplex (i.e., it is either $\emptyset$, a vertex, an
  edge, a face, or the whole of $\closure K$ and $\closure J$) of both
  $\closure K$ and $\closure J$ and
\item $\union{K\in\T{}}\closure K=\closure\W$.
\end{enumerate}
We use the convention where $\funk h\W\reals$ denotes the
\emph{meshsize function} of $\T{}$, \ie 
\begin{equation}
  h(\geovec{x}):=\max_{\closure K\ni \geovec x}h_K,  
\end{equation}
where $h_K$ is the diameter of an element $K$. We let $\E{}$ be the
skeleton (set of common interfaces) of the triangulation $\T{}$ and
say $e\in\E$ if $e$ is on the interior of $\W$ and $e\in\partial\W$ if
$e$ lies on the boundary $\partial\W$.

\begin{Defn}[broken Sobolev spaces, trace spaces]
  \label{defn:broken-sobolev-space}
  We introduce the
  broken Sobolev space
  \begin{equation}
    \sobh{k}(\T{})
    :=
    \ensemble{\phi}
             {\phi|_K\in\sobh{k}(K), \text{ for each } K \in \T{}},
  \end{equation}
  similarly for $\sobh1_0(\T{})$ and $\hon(\T{})$.

  We also make use of functions defined in these broken spaces
  restricted to the skeleton of the triagulation. This requires an
  appropriate trace space
  \begin{equation}
    \Tr{\E}
    :=
    \prod_{K\in\T{}} \leb{2}(\partial K) 
    \subset
    \prod_{K\in\T{}} \sobh{\frac{1}{2}}(K).
  \end{equation}
\end{Defn}

Let $\poly p(\T{})$ denote the space of piecewise polynomials of
degree $p$ over the triangulation $\T{}$ we then introduce the
\emph{finite element spaces}
\begin{gather}
  \label{eqn:def:finite-element-space}
  \fes := \dg{p} = \poly p(\T{})
  \\
  \feszero := \fes \cap \hoz(\T{})
  \\
  \fesn := \fes^d \cap \hon(\T{})
\end{gather}
to be the usual spaces of (discontinuous) piecewise polynomial
functions. For simplicity we will assume that $\fes$ is constant in time.

\begin{Defn}[jumps and averages]
  \label{defn:averages-and-jumps}
  We may define average and jump operators over $\Tr{\E}$ for
  arbitrary scalar, $v\in\Tr{\E}$, and vector valued functions, $\geovec
  v\in\Tr{\E}^d$.
  \begin{equation}
    \label{eqn:average}
    \dfunkmapsto[.]
	        {\avg{\cdot}}
	        v
	        {\Tr{\E\cup \partial\W}}
	        {\frac{1}{2}\qp{v|_{K_1} + v|_{K_2}}}
	        {\leb{2}(\E\cup \partial\W)}
  \end{equation}
  \begin{equation}
    \label{eqn:average-vec}
    \dfunkmapsto[.]
	        {\avg{\cdot}}
	        {\geovec v}
	        {\qp{\Tr{\E\cup \partial\W}}^d}
	        {\frac{1}{2}\qp{\geovec{v}|_{K_1} + \geovec{v}|_{K_2}}}
	        {\qp{\leb{2}(\E\cup \partial\W)}^d}
  \end{equation}
  \begin{equation}
    \label{eqn:jump}
    \dfunkmapsto[.]
	        {\jump{\cdot}}
	        {v}
	        {\Tr{\E\cup \partial\W}}
	        {{{v}|_{K_1} \geovec n_{K_1} + {v}|_{K_2}} \geovec n_{K_2}}
	        {\qp{\leb{2}(\E\cup \partial\W)}^d}
  \end{equation}
  \begin{equation}
    \label{eqn:jump-vec}
    \dfunkmapsto[.]
	        {\jump{\cdot}}
	        {\geovec v}
	        {\qp{\Tr{\E\cup \partial\W}}^d}
	        {{\Transpose{\qp{\geovec{v}|_{K_1}}}\geovec n_{K_1} + \Transpose{\qp{\geovec{v}|_{K_2}}}\geovec n_{K_2}}}
	        {\leb{2}(\E\cup \partial\W)}
  \end{equation}
  \begin{equation}
    \label{eqn:tensor-jump}
    \dfunkmapsto[,]
	        {\tjump{\cdot}}
	        {\geovec v}
	        {\qp{\Tr{\E\cup \partial\W}}^d}
	        {{\geovec{v}|_{K_1} }\otimes \geovec n_{K_1} + \geovec{v}|_{K_2} \otimes\geovec n_{K_2}}
	        {\qp{\leb{2}(\E\cup \partial\W)}^{d\times d}}
  \end{equation}
  where $\geovec n_{K_i}$ denotes the outward pointing normal to $K_i$.
  Note that on the boundary of the domain $\partial\W$ the jump and
  average operators are defined as
  \begin{gather}
    \jump{v}
    \Big\vert_{\partial\W}
    := 
    v\geovec n
    \qquad 
    \jump{\geovec v}
    \Big\vert_{\partial\W} 
    :=
    \Transpose{\geovec v}\geovec n 
    \qquad
    \tjump{\geovec v}
    \Big\vert_{\partial\W} 
    :=
    \geovec v\otimes\geovec n 
    \\
    \avg{v}
    \Big\vert_{\partial\W} 
    := v
    \qquad 
    \avg{\geovec v}
    \Big\vert_{\partial\W}
    :=
    \geovec v.
  \end{gather}
\end{Defn}

\subsection{Discrete mixed formulation}

We propose the following semidiscrete (spatially discrete) formulation
of the system: To find $\big(\phi_h$, $\vec v_h$, $\lambda_h$, $a_h$, $b_h$,
$\vec q_h \big)\in C^1([0,T),\fes)\times C^1([0,T),\feszero^d)\times
C^0([0,T),\fes )\times C^0([0,T),\fes ) \times C^0([0,T),\fes )
\times C^0([0,T),\fesn )$ such that 
\begin{equation}
  \label{eq:discrete-mixed-form}
  \begin{split}
    0
    &=
    \int_\W \qp{\pdt{}\phi_h + \div\qp{\phi_h\geovec v_h} + c_ + m_r a_h} \Chi 
    - c_+ m_j \cA_1(a_h,\Chi)
    -
    \int_\E \jump{\phi_h \geovec v_h}\avg{\Chi}
    \\
    0 &=
    \int_\W 
    \rho(\phi_h) \pdt{}\geovec v_h \cdot \geovec \Xi
    +
    \rho(\phi_h) \qp{ \qp{\geovec v_h \cdot \nabla }\geovec v_h} \cdot \geovec \Xi
    \\
    &\qquad\qquad-
    \frac{1}{2} \rho(\phi_h) \nabla\qp{\norm{\geovec v_h}^2} \cdot \geovec \Xi
    +
    \nabla b_h \cdot \geovec \Xi
    +
    \frac{\phi_h}{c_+}\nabla(a_h - c_- b_h) \cdot \geovec \Xi
    -
    \eta\cA_2\qp{{\geovec v_h},{\geovec \Xi}}
    \\
    & \qquad +
    \int_\E
    \qp{{-}\avg{ \geovec \Xi} \otimes \avg{ \rho(\phi_h) \geovec v_h}}: \jump{\geovec v_h}_\otimes 
    {+}
    \frac{1}{2} \jump{\norm{\geovec v_h}^2}\cdot \avg{\rho(\phi_h) \geovec \Xi}
    \\ 
    &\qquad\qquad
    -
    \jump{b_h} \cdot \avg{\geovec \Xi}
    -
    \frac{1}{c_+}\jump{ a_h - c_- b_h} \cdot \avg{\phi_h \geovec \Xi}
    \\
    0 &= 
    \int_\W 
    \div\qp{\geovec v_h} \Zeta
    -
    \frac{c_-}{c_+}\pdt{}\phi_h \Zeta
    -
    \frac{c_-}{c_+} \div\qp{\phi_h\geovec v_h} \Zeta
    +
    \int_\E \jump{\frac{c_-}{c_+}\phi_h\geovec v_h - \geovec v_h}\avg{\Zeta} 
    \\
    0 &=
    \int_\W
    \qp{a_h - c_+ W'(\phi_h) - c_- \lambda_h} \Psi + c_+ \gamma \div\qp{\vec q_h} \Psi
    - c_+\gamma\int_\E \jump{\vec q_h}\avg{\Psi}
    \\
    0 &=
    \int_\W
    \qp{b_h - \lambda_h -\frac{\rho_1 + \rho_2}{4}\norm{\geovec v_h}^2}\Upsilon
    \\
    0 &=
    \int_\W 
    \vec q_h \cdot \vec \Tau
    -
    \nabla\phi_h \cdot \vec \Tau
    +
    \int_\E 
    \jump{\phi_h} \cdot \avg{\vec \Tau}
    \\
    &\qquad\qquad \Foreach \qp{\Chi, \geovec \Xi, \Zeta, \Psi, \Upsilon, \geovec \Tau} 
    {\in
      \fes\times\feszero^d\times\fes\times\fes\times\fes\times\fesn.}
  \end{split}
\end{equation}
Where 
\begin{equation}
  \begin{split}
    \cA_1\qp{a_h, \Chi}
    &=
    -\int_\W 
    \nabla a_h \cdot \nabla \Chi
    + 
    \int_{\E} \avg{\nabla \Chi} \cdot \jump{a_h}
    \\
    & \qquad +
    \int_\E \jump{\Chi}\cdot \avg{\nabla a_h}
    -
    \frac{\sigma}{h} \jump{a_h} \cdot \jump{\Chi}
    \\
    \cA_2\qp{\geovec v_h, \geovec \Xi}
    &=
    -\int_\W 
    \frob{\D \geovec v_h}{\D \geovec \Xi}
    +
    \int_{\E\cup\partial\W} 
    \frob{\avg{\D \geovec \Xi}}{\tjump{\geovec v_h}}
    \\
    &\qquad +
    \int_{\E\cup\partial\W} 
    \frob{\avg{\D \geovec v_h}}{\tjump{\geovec \Xi}}
    -
    \frac{\sigma}{h}\frob{\tjump{\geovec v_h}}{\tjump{\geovec \Xi}}
  \end{split}
\end{equation}
represent symmetric interior penalty discretisations of the scalar and
vector valued Laplacians respectively, which are signed (coercive)
when the \emph{penalty parameter} $\sigma$ is chosen sufficiently
large.

\begin{Rem}[discrete boundary conditions]
  The boundary conditions (\ref{eq:mixed-bcs}) are encoded in the
  finite element spaces for the Dirichlet type conditions on $\vec
  v_h$ and $\vec q_h$. For $a_h$ the Neumann condition is encoded in
  the bilinear form $\cA_1$.
\end{Rem}

\begin{Rem}[alternative bilinear forms]
  We may choose $\cA_{1,2}$ to be any discretisation of scalar and
  vector valued Laplacian, the only requirement is that they are
  coercive.
\end{Rem}

Throughout the calculations in this section we will regularly refer to
the following proposition.
\begin{Pro}[elementwise integration]
  \label{Pro:trace-jump-avg}
{
  Let 
  \begin{equation}
  \sobh{\div}(\T{}) 
  := 
  \ensemble{\geovec p \in(\leb{2}(\T{}))^d}{\div\qp{\geovec p|_K} \in \leb{2}(K) \text{ for each } K \in \T{}}.
  \end{equation}
}
  Suppose $\geovec p 
  \in \sobh{\div}(\T{})$
  and $\varphi\in\sobh1(\T{})$ then 
  \begin{equation}
    \begin{split}
      \sum_{K\in\T{}}
      \int_K \div\qp{\geovec p} \varphi \d \geovec x
      =
      \sum_{K\in\T{}}
      \qp{
        -
        \int_K
        \geovec p \cdot \nabla \varphi \d \geovec x
        +
        \int_{\partial K}
        \varphi \geovec p \cdot \geovec n_K \d s
        }.
    \end{split}
  \end{equation}
  In particular we have $\geovec p \in \Tr{\E}^d$ and $\varphi \in
  \Tr{\E}$, and the following identity holds
  \begin{equation}
    \sum_{K\in\T{}}
    \int_{\partial K}
    \varphi \Transpose{{\geovec p}}
    \geovec n_K \d s
    =
    \int_\E 
    \jump{\geovec p} 
    \avg{\varphi}
    \d s
    +
    \int_{\E\cup\partial\W}
    \jump{\varphi}\cdot
    \avg{\geovec p}
    \d s
    =
    \int_{\E\cup\partial\W}
    \jump{\geovec p \varphi}
    \d s.
  \end{equation}
\end{Pro}

\begin{Pro}[discrete conservation of mass]  
  \label{pro:cons-of-mass-space}
  The semi discrete scheme (\ref{eq:discrete-mixed-form}) is mass conserving, that is,
  \begin{equation}
    \ddt \qp{\int_\Omega \rho(\phi_h)} =0.
  \end{equation}
\end{Pro}
\begin{Proof}
  Let $1$ be the scalar function which is one everywhere on $\Omega$.
  Then using  $\Zeta= 1$ in \eqref{eq:discrete-mixed-form}$_3$ we see
  \begin{equation}
    0 = 
    \int_\W 
    \div\qp{\geovec v_h} 
    -
    \frac{c_-}{c_+}\pdt{}\phi_h 
    -
    \frac{c_-}{c_+} \div\qp{\phi_h\geovec v_h} 
    +
    \int_\E \jump{\frac{c_-}{c_+}\phi_h\geovec v_h - \geovec v_h}.
  \end{equation}
  We have, using integration by parts, that
  \begin{equation}
    \frac{c_-}{c_+}\ddt \qp{\int_\Omega \phi_h} =0.
  \end{equation}
  This infers the desired result.
\end{Proof}
\begin{Rem}[conservation of momentum]
Note that we have employed a non-conservative discretisation of the momentum equation. Therefore a discrete version of the global momentum balance does not hold in general.
It does not seem feasible to have conservation of momentum and the discrete energy dissipation equality below at the same time.
The situation is similar to the one in \cite{GiesselmannMakridakisPryer:2012} where this problem is elaborated upon in more detail.
\end{Rem}

\begin{The}[discrete energy dissipation equality]
  \label{the:spat-disc-energy}
  Let $\qp{\phi_h, \vec v_h, \lambda_h, a_h, b_h, \vec q_h}$ solve the
  semidiscrete problem (\ref{eq:discrete-mixed-form}) then we have
  that
  \begin{equation}
    \begin{split}
      \ddt 
      \bigg(\int_\W W(\phi_h) 
        +
        \frac{1}{2}\rho(\phi_h) \norm{\vec v_h}^2
        & +
        \frac{1}{2}\gamma \norm{\vec q_h}^2
      \bigg)
      \\
      &=
      \int_\W 
      - 
      m_r \norm{a_h}^2 
      +
      m_j
      \cA_1 \qp{a_h, a_h}
      +
      \eta \cA_2 \qp{\vec v_h, \vec v_h}.
    \end{split}
  \end{equation}
\end{The}
\begin{Proof}
  The proof mimics that of the continuous argument in Theorem
  \ref{the:energy-cont}. To that end we proceed to take the sum of
  (\ref{eq:discrete-mixed-form})$_1$ and
  (\ref{eq:discrete-mixed-form})$_2$ with $\Chi = a_h/c_+$ and $\Xi =
  \vec v_h$, yielding
  \begin{equation}
    \label{eq:energy-1}
    \begin{split}
      0 
      &=
      \int_\W \qp{\pdt{}\phi_h + \div\qp{\phi_h\geovec v_h} + c_ + m_r a_h} \frac{a_h}{c_+} 
      +
      \rho(\phi_h) \pdt{}\geovec v_h \cdot \geovec v_h
      +
      \rho(\phi_h) \qp{ \qp{\geovec v_h \cdot \nabla }\geovec v_h} \cdot \geovec v_h
      \\
      &\qquad
      +
      \int_\W 
      -
      \frac{1}{2} \rho(\phi_h) \nabla\qp{\norm{\geovec v_h}^2} \cdot \geovec v_h
      +
      \nabla b_h \cdot \geovec v_h
      +
      \frac{\phi_h}{c_+}\nabla(a_h - c_- b_h) \cdot \geovec v_h
      \\
      &\qquad
      - 
      c_+ m_j \cA_1(a_h,\frac{a_h}{c_+}) 
      -
      \eta\cA_2\qp{{\geovec v_h},{\geovec v_h}}
      \\
      & \qquad 
      +
      \int_\E -\jump{\phi_h \geovec v_h}\avg{\frac{a_h}{c_+}}
      {-}
      \qp{\avg{ \geovec v_h} \otimes \avg{ \rho(\phi_h) \geovec v_h}}: \jump{\geovec v_h}_\otimes 
      \\ 
      &\qquad
      +
      \int_\E
      \frac{1}{2} \jump{\norm{\geovec v_h}^2}\cdot \avg{\rho(\phi_h) \geovec v_h}
      -
      \jump{b_h} \cdot \avg{\geovec v_h}
      -
      \frac{1}{c_+}\jump{ a_h - c_- b_h} \cdot \avg{\phi_h \geovec v_h}.
    \end{split}
  \end{equation}
  Note that
  \begin{gather}
    \label{eq:energy-2}
      \int_\W \rho(\phi_h) \qp{ \qp{\geovec v_h \cdot \nabla }\geovec
        v_h} \cdot \geovec v_h
      - \frac{1}{2} \rho(\phi_h)
      \nabla\qp{\norm{\geovec v_h}^2} \cdot \geovec v_h = 0 \AND
      \\
      \int_\E
      \qp{\avg{ \geovec v_h} \otimes \avg{ \rho(\phi_h) \geovec v_h}}:
      \jump{\geovec v_h}_\otimes - \frac{1}{2} \jump{\norm{\geovec
          v_h}^2}\cdot \avg{\rho(\phi_h) \geovec v_h} 
      = 0
    \end{gather}
    In addition, we have that
  \begin{equation}
    \label{eq:energy-3}
    \begin{split}
      &\int_\W
      \frac{a_h}{c_+} \div\qp{\phi_h\geovec v_h}
      +
      \frac{\phi_h}{c_+}\nabla a_h \cdot \geovec v_h
      \\
      &\qquad -
      \int_\E \jump{\phi_h \geovec v_h}\avg{\frac{a_h}{c_+}}
      +
      \frac{1}{c_+}\jump{ a_h} \cdot \avg{\phi_h \geovec v_h}
      =
      \frac{1}{c_+}
      \int_\W
      \div\qp{\phi_h a_h \vec v_h} 
      -
      \int_\E
      \jump{\phi_h a_h \vec v_h}
      \\
      &\qquad\qquad\qquad\qquad\qquad\qquad\qquad\qquad\qquad\quad=
      \frac{1}{c_+} \int_{\partial\W} \phi_h a_h \vec v_h \cdot \vec n = 0.
    \end{split}
  \end{equation}
  Taking the observations from (\ref{eq:energy-2}) and
  (\ref{eq:energy-3}) and substituting them into (\ref{eq:energy-1}), we see
  \begin{equation}
    \label{eq:energy-4}
    \begin{split}
      0 &= 
      \int_\W \pdt{}\phi_h  \frac{a_h}{c_+}  + m_r a_h^2
      +
      \rho(\phi_h) \pdt{}\geovec v_h \cdot \geovec v_h
      +
      \nabla b_h \cdot \geovec v_h
      -
      \frac{c_- \phi_h}{c_+}\nabla b_h \cdot \geovec v_h
      \\
      &\qquad 
      - 
      m_j \cA_1(a_h,a_h) 
      -
      \eta\cA_2\qp{{\geovec v_h},{\geovec v_h}}
      \\
      &\qquad
      -
      \int_\E
      \jump{b_h} \cdot \avg{\geovec v_h}
      -
      \frac{c_-}{c_+}\jump{b_h} \cdot \avg{\phi_h \geovec v_h}.
    \end{split}
  \end{equation}
  Now we make use of (\ref{eq:discrete-mixed-form})$_4$ with $\Psi =
  \frac{\pdt{} \phi_h}{c_+}$ on the first term in (\ref{eq:energy-4}) and find
  that
  \begin{equation}
    \label{eq:energy-5}
    \begin{split}
      0
      &= 
      \int_\W \pdt{}\phi_h 
      \qp{ W'(\phi_h) 
        +
        \frac{c_-}{c_+} \lambda_h 
        -
        \gamma \div{\vec q_h}
      }
      +
      m_r a_h^2
      \\
      &\qquad 
      +
      \int_\W
      \rho(\phi_h) \pdt{}\geovec v_h \cdot \geovec v_h
      +
      \nabla b_h \cdot \geovec v_h
      -
      \frac{c_- \phi_h}{c_+}\nabla b_h \cdot \geovec v_h
      \\
      &\qquad 
      - 
      m_j \cA_1(a_h,a_h) 
      -
      \eta\cA_2\qp{{\geovec v_h},{\geovec v_h}}
      \\
      &\qquad
      -
      \int_\E
      \jump{b_h} \cdot \avg{\geovec v_h}
      -
      \frac{c_-}{c_+}\jump{b_h} \cdot \avg{\phi_h \geovec v_h}
      -
      \gamma\jump{\vec q_h} \avg{\pdt \phi_h}.
    \end{split}
  \end{equation}
  Using (\ref{eq:discrete-mixed-form})$_3$ with $\Zeta = b_h$ and integration by parts we have
  that
\begin{equation}
    \label{eq:energy-6}
    \begin{split}
      0
      &= 
      \int_\W \pdt{}\phi_h 
      \qp{ W'(\phi_h) 
        +
        \frac{c_-}{c_+} \lambda_h 
        -
        \gamma \div{\vec q_h}
        {-}
        \frac{c_-}{c_+} b_h
      }
      +
      m_r a_h^2
      +
      \int_\W
      \rho(\phi_h) \pdt{}\geovec v_h \cdot \geovec v_h
      \\
      &\qquad 
      - 
      m_j \cA_1(a_h,a_h) 
      -
      \eta\cA_2\qp{{\geovec v_h},{\geovec v_h}}
      +
      \int_\E
      \gamma\jump{\vec q_h} \avg{\pdt \phi_h}.
    \end{split}
  \end{equation}
  Now using (\ref{eq:discrete-mixed-form})$_5$ with $\Upsilon = \pdt
  \phi_h$ on the second term in (\ref{eq:energy-6}) and integrating
  the third term by parts we see
  \begin{equation}
    \label{eq:energy-7}
    \begin{split}
      0
      &= 
      \int_\W \pdt{}\phi_h 
      \qp{ W'(\phi_h) 
        {-}
        \frac{c_-\qp{\rho_1+\rho_2}}{4 c_+} \norm{\vec v_h}^2
      }
      +
      \gamma {\vec q_h} \cdot \nabla \pdt \phi_h
      +
      m_r a_h^2
      \\
      &\qquad 
      +
      \int_\W
      \rho(\phi_h) \pdt{}\geovec v_h \cdot \geovec v_h
      - 
      m_j \cA_1(a_h,a_h) 
      -
      \eta\cA_2\qp{{\geovec v_h},{\geovec v_h}}
      -
      \int_\E
      \gamma\avg{\vec q_h} \cdot \jump{\pdt \phi_h}.
    \end{split}
  \end{equation}
  Taking the time derivative of (\ref{eq:discrete-mixed-form})$_6$,
  inserting $\vec \Tau = \vec q_h$ and using this on the fourth term
  in (\ref{eq:energy-7}) we find
  \begin{equation}
    \label{eq:energy-8}
    \begin{split}
      0
      &= 
      \int_\W \pdt{}\phi_h 
      \qp{ W'(\phi_h) 
        {-}
        \frac{\rho_2-\rho_1}{4} \norm{\vec v_h}^2
      }
      +
      \gamma {\vec q_h} \cdot \pdt{} \vec q_h
      +
      m_r a_h^2
      +
      \rho(\phi_h) \pdt{}\geovec v_h \cdot \geovec v_h
      \\
      &\qquad 
      - 
      m_j \cA_1(a_h,a_h) 
      -
      \eta\cA_2\qp{{\geovec v_h},{\geovec v_h}}
      ,
    \end{split}
  \end{equation}
  which infers the desired result, concluding the proof.
\end{Proof}

\begin{Rem}[uniqueness of fluxes]
  The choice of fluxes in the spatially discrete formulation is not
  unique. Indeed, using the more general framework given in
  \cite{GiesselmannMakridakisPryer:2012} we may give conditions for
  families of fluxes which admit energy consistent schemes.
\end{Rem}

\section{Temporally discrete approximation}
\label{sec:temporally-discrete}

In this section we present a methodology for designing temporally
discrete energy consistent discretisations of the system
(\ref{eq:quasi-crazy-pde-1})--(\ref{eq:quasi-crazy-pde-3}). We do this
by appropriately modifying a Crank--Nicolson type temporal
discretisation. The resultant scheme is of $2$nd order. Higher order
energy consistent discretiations can be designed based on
appropriately modifying symplectic Gauss--Legendre type Runge--Kutta
schemes.

Let $[0,T]$ be the time interval in which we approximate the
quasi-incompressible system. We subdivide the time interval $[0,T]$
into a partition of $N$ consecutive adjacent subintervals whose
endpoints are denoted $t_0=0<t_1<\ldots<t_{N}=T$.  The $n$-th timestep
is defined as ${k_n := t_{n+1} - t_{n}}$.  We will consistently use
the shorthand $F^n(\cdot):=F(\cdot,t_n)$ for a generic time function
$F$. We also denote $\nplush{F} := \frac{1}{2}\qp{F^n + F^{n+1}}$.

The semidiscrete (temporally discrete) formulation of the system
(\ref{eq:quasi-crazy-pde-1})--(\ref{eq:quasi-crazy-pde-3}) is: Given
initial conditions $\rho^0$, $\vec v^0$, $\lambda^0$, $a^0$, $b^0$ and
$\vec q^0$, for each
{$n\in\naturals_{0}$} find
$\rho^{n+1}$, $\vec v^{n+1}$, $\lambda^{n+1}$, $a^{n+1}$, $b^{n+1}$
and $\vec q^{n+1}$ such that
\begin{equation}
  \label{eq:temporally-discrete}
  \begin{split}
    0
    &= 
    \frac{\phi^{n+1} - \phi^n}{k_n} 
    +
    \div\qp{\nplush{\phi} \nplush{\geovec v}}
    -
    c_+ m_j \Delta \nplush{a}
    +
    c_+ m_r \nplush{a}
    \\
   \vec 0 
    &= 
    \rho(\nplush{\phi})
    \qp{
      \frac{\vec{v}^{n+1} - \vec{v}^n}{k_n}
      +
      \qp{\nplush{\vec{v}}\cdot \nabla}\nplush{\vec{v}}
      -
      \frac{1}{2}
      \nabla\qp{\norm{\nplush{\vec{v}}}^2}
    }
    \\
    &\qquad - 
    \eta \Delta \nplush{\geovec{v}}
    +
    \nabla \nplush{b}
    +
    \frac{\nplush{\phi}}{c_+}\nabla( \nplush{a} - {c_-} \nplush{b})
    \\
    0 
    &= 
    \div\qp{\nplush{\geovec v}}
    -
    \frac{c_-}{c_+}
    \qp{
      \frac{\phi^{n+1} - \phi^n}{k_n}
      +
      \div\qp{\nplush{\phi}\nplush{\geovec v}}
    }
    \\
    0
    &=
    \nplush{a}
    - 
    c_+\frac{W(\phi^{n+1}) - W(\phi^n)}{\phi^{n+1}-\phi^n}
    +
    c_+\gamma\div\qp{\nplush{\vec{q}}}
    -
    c_-\nplush{\lambda}
    \\
    0&= 
    \nplush{b}
    -
    \nplush{\lambda}
    -
    \frac{\rho_1+\rho_2}{8}
    \qp{\norm{\geovec{v}^{n+1}}^2 + \norm{\geovec{v}^{n}}^2}
    \\
    \vec 0 
    &=
    {\nplush{\vec q}} 
    -
    \nabla \nplush{\phi},
  \end{split}
\end{equation}
satisfying the boundary conditions
\begin{equation}
  \label{eq:temp-bcs}
  \vec q^n \cdot \vec n = 0 , \qquad \vec v^n =\vec  0, \qquad \nabla a^n \cdot \vec n = 0,
\end{equation}
for each $n\in [0, N]$.
\begin{Pro}[temporally discrete mass conservation]
  \label{pro:mass-cons-temp}
  The temporally discrete scheme (\ref{eq:temporally-discrete})
  satisfies
  \begin{equation}
    \int_\W \rho(\phi^{n+1}) = \int_\W \rho(\phi^n) \Foreach n\in [0,N-1
  \end{equation}
\end{Pro}
\begin{Proof}
  For $\rho_1=\rho_2$ the assertion is trivial. Thus, we may assume
  $c_- \not=0$ for the rest of this proof. Integrating
  (\ref{eq:temporally-discrete})$_3$ over the domain we have that
  \begin{equation}
    0 
    = 
    \int_\W
    \div\qp{\nplush{\geovec v}}
    -
    \frac{c_-}{c_+}
    \qp{
      \frac{\phi^{n+1} - \phi^n}{k_n}
      +
      \div\qp{\nplush{\phi}\nplush{\geovec v}}
    }.
  \end{equation}
  In view of Stokes Theorem and making use of the boundary conditions
  (\ref{eq:temp-bcs}) we see that
  \begin{equation}
    0
    =
    \int_\W 
    \frac{c_-}{c_+}
    \frac{\phi^{n+1} - \phi^n}{k_n}.
    \end{equation}
    This infers that
    \begin{equation}
      \int_\W \phi^{n+1} = \int_\W \phi^n,
    \end{equation}
    which, in view of the linearity of $\rho(\phi^n)$, yields the desired result.
\end{Proof}

\begin{The}[temporally discrete energy dissipation equality]
  \label{the:temp-discrete-energy}
  Let $\{\rho^{n}$, $\vec v^{n}$, $\lambda^{n}$, $a^{n}$, $b^{n}$,
  $\vec q^{n}\}_{n\in [0, N]}$ be the sequence generated by the
  semidiscrete scheme (\ref{eq:temporally-discrete}) then we have that
  for any $n \in [0,N]$
  \begin{equation}
    \label{eq:inductive-hyp}
    \begin{split}
      \int_\W W(\phi^n) + \frac{1}{2}\rho(\phi^n) \norm{\vec v^n}^2 +
      \frac{\gamma}{2} \norm{\vec q^n}
      &= \int_\W W(\phi^0) +
      \frac{1}{2}\rho(\phi^0) \norm{\vec v^0}^2 + \frac{\gamma}{2}
      \norm{\vec q^0} 
      \\&\qquad - \sum_{i=0}^{n-1} \Bigg( k_i \int_\W m_j
        \norm{\nabla \iplush{a}}^2 + m_r \norm{\iplush{a}}^2 \\
        &\qquad \qquad \qquad \qquad + \eta
        \norm{\D \iplush{\vec v}}^2 \Bigg).
    \end{split}
  \end{equation}
\end{The}
\begin{Proof}
  We will prove this using induction. Our inductive hypothesis is
  given by (\ref{eq:inductive-hyp}). It is clear that
  (\ref{eq:inductive-hyp}) holds in the case $n=0$. We then assume
  that (\ref{eq:inductive-hyp}) holds for all $k \leq n$ and make our
  inductive step.

  Using the semidiscrete scheme (\ref{eq:temporally-discrete}),
  testing the first equation (\ref{eq:temporally-discrete})$_1$ with
  $\nplush{a}$ and the second (\ref{eq:temporally-discrete})$_2$ with
  $\nplush{\vec v}$ and taking the sum we have
  \begin{equation}
    \begin{split}
      0
      &= 
      \int_\W
      \frac{\nplush{a}}{c_+}\qp{\frac{\phi^{n+1} - \phi^n}{k_n} 
      +
      \div\qp{\nplush{\phi} \nplush{\geovec v}}
      -
      c_+ m_j \Delta \nplush{a}
      +
      c_+ m_r \nplush{a}}
      \\
      &\qquad + \nplush{\vec v}\cdot \bigg(\rho(\nplush{\phi})
      \qp{
        \frac{\vec{v}^{n+1} - \vec{v}^n}{k_n}
        +
        \qp{\nplush{\vec{v}}\cdot \nabla}\nplush{\vec{v}}
        -
        \frac{1}{2}
        \nabla\qp{\norm{\nplush{\vec{v}}}^2}
      }
      \\
      &\qquad - 
      \eta \Delta \nplush{\geovec{v}}
      +
      \nabla \nplush{b}
      +
      \frac{\nplush{\phi}}{c_+}\nabla( \nplush{a} - {c_-} \nplush{b})\bigg).
    \end{split}
  \end{equation}
  In view of the same arguments given in the proof of Theorem
  \ref{the:energy-cont} we see, upon integrating by parts, that
  \begin{equation}
    \label{eq:temp-energy-0}
    \begin{split}
      0 
      &=
      \int_\W
      \qp{\phi^{n+1}-\phi^n} \frac{\nplush a}{c_+}
      +
      k_n \qp{m_j\norm{\nabla \nplush a}^2
        +
        m_r \norm{\nplush a}^2
        +
        \eta \norm{\D \nplush{\vec v}}^2}
      \\
      &\qquad\qquad
      +
      k_n \nabla \nplush b \cdot \nplush{\vec v}
      -
      k_n \nplush \phi \frac{c_-}{c_+} \nabla \nplush b \cdot \nplush{\vec v}
      \\
      &\qquad\qquad 
      +
      \rho(\nplush \phi) \qp{\vec v^{n+1} - \vec v^n} \cdot \nplush{\vec v}
      \\
      &\qquad
      -
      k_n
      \int_{\partial \W} 
      m_j \nabla \nplush a \cdot \vec n \nplush a
      +
      \eta \qp{\D \nplush{\vec v} \vec n} \cdot \nplush{\vec v}
      \\
      &\qquad \qquad
      + \frac{1}{c_+} \nplush\phi \nplush a \nplush{\vec v} \cdot \vec n.
    \end{split}
  \end{equation}
  Note that the boundary terms vanish due to (\ref{eq:temp-bcs}). Now
  testing (\ref{eq:temporally-discrete})$_3$ with $\nplush b$ we see
  \begin{equation}
    \label{eq:temp-energy-1}
    \begin{split}
      0 
      &=
      \int_\W
      k_n \div\qp{\nplush{\vec v}} \nplush b
      -
      \frac{c_-}{c_+}
      \qp{\phi^{n+1} - \phi^n}
      \nplush{b}
      -
      \frac{k_n c_-}{c_+}
      \div\qp{\nplush\phi \nplush{\vec v}} \nplush b
      \\
      &=
      \int_\W 
      -
      k_n \nplush{\vec v} \cdot \nabla \nplush b
      -
      \frac{c_-}{c_+}
      \qp{\phi^{n+1} - \phi^n}
      \nplush{b}
      +
      \frac{k_n c_-}{c_+}
      \nplush\phi \nplush{\vec v} \cdot \nabla\nplush b
      \\
      &\qquad
      +
      \int_{\partial \W}
      k_n \nplush{\vec v} \cdot \vec n \nplush b
      -
      \frac{k_n c_-}{c_+}
      \nplush\phi \nplush{\vec v} \cdot \vec n \nplush b.
    \end{split}
  \end{equation}
  Notice again that the boundary terms vanish due to
  (\ref{eq:temp-bcs}). Testing (\ref{eq:temporally-discrete})$_5$ with
  $\qp{\phi^{n+1} - \phi^n}$ we have that
  \begin{equation}
    \label{eq:temp-energy-2}
    0= 
    \int_\W
    \qp{\phi^{n+1} - \phi^n}
    \qp{\nplush{b}
    -
    \nplush{\lambda}
    -
    \frac{\rho_1+\rho_2}{8}
    \qp{\norm{\geovec{v}^{n+1}}^2 + \norm{\geovec{v}^{n}}^2}}.
  \end{equation}
  Substituting (\ref{eq:temp-energy-1}) and (\ref{eq:temp-energy-2})
  into (\ref{eq:temp-energy-0}), we have
  \begin{equation}
    \begin{split}
      0
      &= 
      \int_\W
      W(\phi^{n+1}) - W(\phi^n)
      -
      \gamma \qp{\phi^{n+1} - \phi^n} \div\qp{\nplush {\vec q}}
      \\
      &\qquad +
      \rho(\nplush{\phi}) \qp{\vec v^{n+1} - \vec v^n} \cdot \nplush{\vec v} 
      -
      \frac{c_-\qp{\rho_1 + \rho_2}}{8 c_+} \qp{\norm{\vec v^{n+1}}^2 + \norm{\vec v^n}^2} \qp{\phi^{n+1} - \phi^n}
      \\
      &\qquad 
      +
      k_n \qp{m_j\norm{\nabla \nplush a}^2
        +
        m_r \norm{\nplush a}^2
        +
        \eta \norm{\D \nplush{\vec v}}^2}
      \\
      &=
      \int_\W
      W(\phi^{n+1}) - W(\phi^n)
      -
      \gamma \qp{\phi^{n+1} - \phi^n} \div\qp{\nplush {\vec q}}
      \\
      &\qquad +
      \frac{1}{2}\rho(\nplush{\phi})\qp{\norm{\vec v^{n+1}}^2 - \norm{\vec v^n}^2} - \frac{\rho_2 - \rho_1}{8} \qp{\phi^{n+1} - \phi^n} \qp{\norm{\vec v^{n+1}}^2 + \norm{\vec v^n}^2} 
      \\
      &\qquad 
      +
      k_n \qp{m_j\norm{\nabla \nplush a}^2
        +
        m_r \norm{\nplush a}^2
        +
        \eta \norm{\D \nplush{\vec v}}^2}.
    \end{split}    
  \end{equation}
  Using the identities
  \begin{gather}
    \rho(\nplush \phi) = \frac{1}{2} \qp{\rho(\phi^{n+1}) + \rho(\phi^n)}
    \\
    -\frac{\rho_2-\rho_1}{8} \qp{\phi^{n+1}-\phi^n}
    =
    \frac{1}{4} \qp{\rho(\phi^{n+1}) - \rho(\phi^n)},
  \end{gather}
  we have
  \begin{equation}
    \begin{split}
      0
      &=
      \int_\W
      W(\phi^{n+1}) - W(\phi^n)
      -
      \gamma \qp{\phi^{n+1} - \phi^n} \div\qp{\nplush {\vec q}}
      \\
      &\qquad +
      \frac{1}{2}\qp{\rho(\phi^{n+1})\norm{\vec v^{n+1}}^2 - \rho(\phi^n) \norm{\vec v^n}^2}    
      \\
      &\qquad 
      +
      k_n \qp{m_j\norm{\nabla \nplush a}^2
        +
        m_r \norm{\nplush a}^2
        +
        \eta \norm{\D \nplush{\vec v}}^2}.
    \end{split}
  \end{equation}
  Now using the fact that
  \begin{equation}
    \begin{split}
      \int_\W - \gamma \qp{\phi^{n+1} - \phi^n} \div\qp{\nplush {\vec
          q}} &= 
      \int_\W \gamma \nabla \qp{\phi^{n+1} - \phi^n}
      \nplush {\vec q}
      \\
      &\qquad
      - 
      \int_{\partial\W} \gamma \qp{\phi^{n+1} - \phi^n}
      \nplush {\vec q} \cdot \geovec n
      \\
      &=
      \int_\W \frac{\gamma}{2} \qp{\vec q^{n+1} - \vec q^n} \cdot
      \qp{\vec q^{n+1} + \vec q^n}
      \\
      &=
      \int_\W \frac{\gamma}{2} \qp{\norm{\vec q^{n+1}}^2 - \norm{\vec q^n}^2},
    \end{split}
  \end{equation}
  by (\ref{eq:temporally-discrete})$_6$, we see
  \begin{equation}
    \begin{split}
      \int_\W
      W(\phi^{n+1})
      &+
      \frac{\gamma}{2} \norm{\vec q^{n+1}}^2
      +
      \frac{1}{2}\rho(\phi^{n+1})\norm{\vec v^{n+1}}^2
      \\
      &=
      \int_\W
      W(\phi^n)
      +
      \frac{\gamma}{2} \norm{\vec q^{n}}^2
      +
      \frac{1}{2}\rho(\phi^n)\norm{\vec v^n}^2
      \\
      &\qquad\qquad
      +
      \int_\W
      k_n \qp{m_j\norm{\nabla \nplush a}^2
        +
        m_r \norm{\nplush a}^2
        +
        \eta \norm{\D \nplush{\vec v}}^2},
    \end{split}
  \end{equation}
  which, using the inductive hypothesis (\ref{eq:inductive-hyp}), concludes the proof.
\end{Proof}

\section{A fully discrete approximation}
\label{sec:fully-discrete}

In this section we present a fully discrete approximation of
(\ref{eq:quasi-crazy-pde-1})--(\ref{eq:quasi-crazy-pde-3}) which is
energy consistent.

Collecting the results of \S\ref{sec:spatially-discrete} and
\S\ref{sec:temporally-discrete} we propose the following scheme:

\begin{equation}
  \label{eq:fully-discrete-mixed-form}
  \begin{split}
    0
    &=
    \int_\W \qp{\frac{\phi_h^{n+1}-\phi_h^n}{k_n} + \div\qp{\nplush{\phi_h}\nplush{\geovec v_h}} + c_ + m_r \nplush{a_h}} \Chi 
    \\
    &\qquad
    - c_+ m_j \cA_1(\nplush{a_h},\Chi)
    {-}\int_\E \jump{\nplush{\phi_h} \nplush{\geovec v_h}}\avg{\Chi}
    \\
    0 &=
    \int_\W 
    \rho(\nplush{\phi_h}) \frac{\geovec v_h^{n+1} - \vec v_h^n}{k_n} \cdot \geovec \Xi
    +
    \rho(\nplush{\phi_h}) \qp{ \qp{\nplush{\geovec v_h} \cdot \nabla }\nplush{\geovec v_h}} \cdot \geovec \Xi
    \\
    &\qquad\qquad-
    \frac{1}{2} \rho(\nplush{\phi_h}) \nabla\qp{\norm{\nplush{\geovec v_h}}^2} \cdot \geovec \Xi
    -
    \eta\cA_2\qp{{\nplush{\geovec v_h}},{\geovec \Xi}}
    +
    \nabla \nplush{b_h} \cdot \geovec \Xi
    \\
    &\qquad\qquad +
    \frac{\nplush{\phi_h}}{c_+}\nabla(\nplush{a_h} - c_- \nplush{b_h}) \cdot \geovec \Xi
    \\
    & \qquad +
    \int_\E
    \qp{{-}\avg{ \geovec \Xi} \otimes \avg{ \rho(\nplush{\phi_h}) \nplush{\geovec v_h}}}: \tjump{\nplush{\geovec v_h}} 
    \\&\qquad \qquad {+}
    \frac{1}{2} \jump{\norm{\nplush{\geovec v_h}}^2}\cdot \avg{\rho(\nplush{\phi_h}) \geovec \Xi}
    \\ 
    &\qquad\qquad
    -
    \jump{\nplush{b_h}} \cdot \avg{\geovec \Xi}
    -
    \frac{1}{c_+}\jump{\nplush{ a_h} - c_-\nplush{ b_h}} \cdot \avg{\nplush{\phi_h} \geovec \Xi}
    \\
    0 &= 
    \int_\W 
    \div\qp{\nplush{\geovec v_h}} \Zeta
    -
    \frac{c_-}{c_+}\frac{\phi_h^{n+1}-\phi_h^n}{k_n} \Zeta
    -
    \frac{c_-}{c_+} \div\qp{\nplush{\phi_h}\nplush{\geovec v_h}} \Zeta
    \\&\qquad +
    \int_\E \jump{\frac{c_-}{c_+}\nplush{\phi_h}\nplush{\geovec v_h} - \nplush{\geovec v_h}}\avg{\Zeta} 
    \\
    0 &=
    \int_\W
    \qp{\nplush{a_h} - c_+ \frac{W(\phi_h^{n+1})-W(\phi_h^n)}{\phi_h^{n+1}-\phi_h^n} - c_- \nplush{\lambda_h}} \Psi + c_+ \gamma \div\qp{\nplush{\vec q_h}} \Psi
    \\
    &\qquad - c_+\gamma\int_\E \jump{\nplush{\vec q_h}}\avg{\Psi}
    \\
    0 &=
    \int_\W
    \qp{\nplush{b_h} - \nplush{\lambda_h} -\frac{\rho_1 + \rho_2}{8}\qp{\norm{\geovec v_h^{n+1}}^2 + \norm{\vec v_h^n}^2}}\Upsilon
    \\
    0 &=
    \int_\W 
    \nplush{\vec q_h} \cdot \vec \Tau
    -
    \nabla \nplush{\phi_h} \cdot \vec \Tau
    +
    \int_\E 
    \jump{\nplush{\phi_h}} \cdot \avg{\vec \Tau}
    \\
    &\qquad\qquad \Foreach \qp{\Chi, \geovec \Xi, \Zeta, \Psi, \Upsilon,\geovec \Tau} 
    \in
    {\fes\times\feszero^d\times\fes\times\fes\times\fes\times\fesn.}
  \end{split}
\end{equation}
\begin{Pro}
The fully discrete scheme (\ref{eq:fully-discrete-mixed-form}) is mass conservative, \ie
\begin{equation}
  \int_\Omega \rho(\phi_h^{n+1}) = \int_\Omega \rho(\phi_h^n).
\end{equation}
\end{Pro}
\begin{Proof}
  The proof is given by combining Propositions
  \ref{pro:cons-of-mass-space} and \ref{pro:mass-cons-temp} which
  yield the spatial and temporal semidiscrete mass conservation
  results respectively.
\end{Proof}

\begin{The}[fully discrete energy consistent approximation]
  \label{the:fully-discrete-energy}
  The sequence of solutions generated by the fully discrete
  approximation (\ref{eq:fully-discrete-mixed-form}) satisfies the
  following energy identity:
  \begin{equation}
    \begin{split}
      \int_\W 
      W(\phi_h^{n+1}) 
      &+
      \frac{1}{2} \rho(\phi_h^{n+1})\norm{\vec v_h^{n+1}}^2
      +
      \frac{\gamma}{2} \norm{\vec q_h^{n+1}}^2
      \\
      &=
      \int_\W 
      W(\phi_h^{n}) 
      +
      \frac{1}{2} \rho(\phi_h^{n})\norm{\vec v_h^{n}}^2
      +
      \frac{\gamma}{2} \norm{\vec q_h^{n}}^2
      \\
      &\qquad -
      k_n 
      \Bigg(\int_\W m_r \qp{\nplush{a_h}}^2
        -
        m_j \cA_1\qp{\nplush{a_h}, \nplush{a_h}} 
        \\
        &\qquad\qquad -
        \eta \cA_2\qp{\nplush{\vec v_h}, \nplush{\vec v_h}}
      \Bigg).
    \end{split}
  \end{equation}
\end{The}
\begin{Proof}
  The proof follows those of Theorem \ref{the:spat-disc-energy} and
  Theorem \ref{the:temp-discrete-energy}.
\end{Proof}

\begin{Rem}[Adaptive interface tracking]
  { Resolution of the diffuse interface is of paramount
    importance for both stability and long time accuracy of the
    numerical method. The restrictions placed upon $\T{}$ in
    \S\ref{sec:spatially-discrete} do not proclude the use of
    adaptivity to refine the mesh in proximity of the
    interface. Indeed, it is possible to design heuristic adaptive
    schemes based on local adaptive refinement/coarsening routines as
    dictated by gradient aposteriori indicators for $\phi$, for example. }
\end{Rem}

\section{Numerical experiments}
\label{sec:numerics}

In this section we conduct a series of numerical experiments aimed at
testing the robustness of the method.

\subsection{Implementation issues}

The numerical experiments were conducted using the \dolfin interface
for \fenics \cite{LoggWells:2010}. The graphics were generated using
\gnuplot and \paraview.

In each of the numerical experiments we fix $W$ to be the following
quartic double well potential
\begin{equation}
  \label{eq:double-well}
  W(\phi) = \qp{\phi^2 - 1}^2
\end{equation}
with minima at $\phi = \pm 1$.

\begin{Rem}[the quotient of the double well]
  In the computational implementation we did not use the difference
  quotient $\tfrac{W(\phi^{n+1})-W(\phi^n)}{\phi^{n+1} - \phi^n}$
  appearing in (\ref{eq:fully-discrete-mixed-form}) as it is ill-defined for
  $\phi^{n+1} = \phi^n$ and badly conditioned when $|\phi^{n+1} -
  \phi^n|$ is small. Instead we use a sufficiently high order
  approximation of this term. For (\ref{eq:double-well}) we use the
  following Taylor expansion representation
  \begin{equation}
    \frac{W(\phi^{n+1})-W(\phi^n)}{\phi^{n+1} - \phi^n}
    =
    W'(\nplush{\phi})
    +
    \tfrac{1}{24}W'''(\nplush{\phi})\qp{\phi^{n+1} - \phi^{n}}^2
  \end{equation}
  which is exact. We note that when $W$ is not polynomial a
  sufficiently high order truncation of the Taylor expansion can be
  achieved such that the possible increase in energy is of high order with
  respect to the timestep. This allows the
  construction of a method with arbitrarily small deviations of the
  energy with respect to the timestep.
\end{Rem}

\begin{Rem}[default parameters]
  In each of the following tests, unless otherwise specified, we take
  the parameters as follows: We set $\rho_1 = 1, \rho_2 = 2$, $\gamma
  = \eta = 10^{-3}, m_r = m_j = 10^{-2}$, $h \approx 0.01$, $\tau =
  0.01$ and $p=1$.
\end{Rem}

\subsection{Test 1 : 1D - benchmarking}

In this test we benchmark the numerical algorithm presented in
\S\ref{sec:fully-discrete} against a steady state solution of the
quasi-incompressible system
(\ref{eq:quasi-crazy-pde-1})--(\ref{eq:quasi-crazy-pde-3}) in one
spatial dimension on the domain $\W = [-1,1]$.

For the double well given by (\ref{eq:double-well}) a steady state
solution to the quasi-incompressible system is given by
\begin{gather}
  \phi(x,t) 
  =
  \tanh\qp{ x \sqrt{\frac{2}{\gamma}}},
  \qquad v(x,t) \equiv 0 \Foreach t.
\end{gather}
Note that on the boundary $\nabla \phi$ is not zero but of negligible
value (as $\gamma$ is small). Tables
\ref{table:p1-gamma-10-3}--\ref{table:p3-gamma-10-3} detail three
experiments aimed at testing the convergence properties for the scheme
using piecewise discontinuous elements of various orders ($p=1$ in
Table \ref{table:p1-gamma-10-3}, $p=2$ in Table
\ref{table:p2-gamma-10-3} and $p=3$ in Table
\ref{table:p3-gamma-10-3}).

\begin{table}
  \caption{\label{table:p1-gamma-10-3} In this test we benchmark a
    stationary solution of the quasi-incompressible system using the discretisation
    (\ref{eq:fully-discrete-mixed-form}) with piecewise linear elements ($p =
    1$), choosing $k = h^2$. This is done by formulating
    (\ref{eq:fully-discrete-mixed-form}) as a system of nonlinear equations, the
    solution to this is then approximated by a Newton method with
    tolerance set at $10^{-10}$. At each Newton step the solution to
    the linear system of equations is approximated using a stabilised
    conjugate gradient iterative solver with an successively
    overrelaxed preconditioner, also set at a tolerance of
    $10^{-10}$. We look at the $\leb{\infty}(0,T; \leb{2}(\W))$ errors
    of the discrete variables {$\phi_h,\, v_h$ and $\lambda_h$}, and use $e_\phi :=
    \phi - \phi_h,\ e_v := v - v_h$ {and $e_\lambda= \lambda-\lambda_h$}. In this test we choose
    $\gamma = 10^{-3}$.}
  \begin{center}
    \begin{tabular}{c|c|c|c|c|c|c}
$N$ & $\Norm{e_\phi}_{\leb{\infty}(\leb{2})}$ & EOC & $\Norm{e_v}_{\leb{\infty}(\leb{2})}$ & EOC & $\Norm{e_\lambda}_{\leb{\infty}(\leb{2})}$ & EOC \\
\hline
32 & 1.4998e-01 & 0.000 & 6.9600e-02 & 0.000 & 9.7289e-01 & 0.000 \\
64 & 9.4503e-02 & 0.666 & 5.3907e-02 & 0.369 & 6.7654e-01 & 0.524 \\
128 & 4.0138e-02 & 1.235 & 3.5739e-02 & 0.593 & 4.6306e-01 & 0.547 \\
256 & 9.8587e-03 & 2.026 & 1.6355e-02 & 1.128 & 3.3446e-01 & 0.469 \\
512 & 2.8050e-03 & 1.813 & 5.8975e-03 & 1.472 & 2.2825e-01 & 0.551 \\
1024 & 6.7240e-04 & 2.061 & 1.8467e-03 & 1.675 & 1.3269e-01 & 0.783 \\
2048 & 1.5217e-04 & 2.144 & 4.1273e-04 & 2.162 & 6.9219e-02 & 0.939 \\
4096 & 3.7793e-05 & 2.010 & 5.9895e-05 & 2.785 & 3.4988e-02 & 0.984 \\
\end{tabular}
  \end{center}
\end{table}

\begin{table}
  \caption{\label{table:p2-gamma-10-3} The test is the same as in
    Table \ref{table:p1-gamma-10-3} with the exception that
    we take $p=2$.}
  \begin{center}
    \begin{tabular}{c|c|c|c|c|c|c}
$N$ & $\Norm{e_\phi}_{\leb{\infty}(\leb{2})}$ & EOC & $\Norm{e_v}_{\leb{\infty}(\leb{2})}$ & EOC & $\Norm{e_\lambda}_{\leb{\infty}(\leb{2})}$ & EOC \\
\hline
32 & 6.8671e-02 & 0.000 & 4.7711e-02 & 0.000 & 6.8098e-01 & 0.000 \\
64 & 2.8248e-02 & 1.282 & 2.6617e-02 & 0.842 & 3.3259e-01 & 1.034 \\
128 & 6.7024e-03 & 2.075 & 7.7866e-03 & 1.773 & 2.1021e-01 & 0.662 \\
256 & 2.1369e-03 & 1.649 & 5.3622e-03 & 0.538 & 1.9486e-01 & 0.109 \\
512 & 1.7291e-04 & 3.627 & 1.8418e-03 & 1.542 & 1.2747e-01 & 0.612 \\
1024 & 1.8023e-05 & 3.262 & 4.7102e-04 & 1.967 & 6.5608e-02 & 0.958 \\
2048 & 2.1668e-06 & 3.056 & 1.1910e-04 & 1.984 & 3.2833e-02 & 0.999 \\
4096 & 2.6758e-07 & 3.018 & 2.9902e-05 & 1.994 & 1.6729e-02 & 0.973 \\
\end{tabular}
  \end{center}
\end{table}

\begin{table}
  \caption{\label{table:p3-gamma-10-3} The test is the same as in
    Table \ref{table:p1-gamma-10-3} with the exception that
    we take $p=3$.}
  \begin{center}
    \begin{tabular}{c|c|c|c|c|c|c}
$N$ & $\Norm{e_\phi}_{\leb{\infty}(\leb{2})}$ & EOC & $\Norm{e_v}_{\leb{\infty}(\leb{2})}$ & EOC & $\Norm{e_\lambda}_{\leb{\infty}(\leb{2})}$ & EOC \\
\hline
32 & 3.3914e-02 & 0.000 & 2.1390e-02 & 0.000 & 3.2962e-01 & 0.000 \\
64 & 1.0777e-02 & 1.654 & 8.5393e-03 & 1.325 & 2.2624e-01 & 0.543 \\
128 & 3.4979e-03 & 1.623 & 7.6267e-03 & 0.163 & 2.1279e-01 & 0.088 \\
256 & 2.0816e-04 & 4.071 & 1.8900e-03 & 2.013 & 9.8126e-02 & 1.117 \\
512 & 1.3447e-05 & 3.952 & 1.6423e-04 & 3.525 & 1.4974e-02 & 2.712 \\
1024 & 1.4090e-06 & 3.255 & 1.5439e-05 & 3.411 & 2.6407e-03 & 2.503 \\
2048 & 1.3055e-07 & 3.432 & 1.5523e-06 & 3.314 & 3.9831e-04 & 2.729 \\
\end{tabular}
  \end{center}
\end{table}
\begin{Rem}[optimality of the primal variables]
  Note that the results presented (and various other tests) indicate
  that
  \begin{gather}
    \Norm{e_\phi} = \Oh(k^2 + h^{p+1})
    \\
    \Norm{e_{\vec v}} = 
    \begin{cases}
      \Oh(k^2 + h^{p+1}) \text{ if } p \text { is odd }
      \\
      \Oh(k^2 + h^{p}) \text{ if } p \text { is even }
    \end{cases}
    \\
    \Norm{e_\lambda} = 
    \begin{cases}
      \Oh(k^2 + h^{p}) \text{ if } p \text{ is odd }
      \\
      \Oh(k^2 + h^{p-1}) \text{ if } p \text{ is even }
    \end{cases}
  \end{gather}
  As such, we see the convergence rates are optimal for $\phi$ and
  $\vec v$ if $p$ is odd. This suboptimality in $\vec v$ for even
  order finite element spaces has been observed previously
  \cite{GiesselmannMakridakisPryer:2012}. Regarding the suboptimality
  of $\lambda$ we note that the energy dissipation equality provides
  no stability for $\lambda.$
\end{Rem}

\subsection{Test 2 : 2D - random initial data}
\label{sec:2d-random}

In this test we examine the behaviour of the solution when the initial
conditions for $\phi$ are random perturbations of the unstable extremum
of the double well. More precisely, let $\{ x_i \}_{i = 1}^M$ denote
the mesh points of the triangulation $\T{}$ of $\W = [-1,1]^2$. We then let $Y_{i} \sim
\text{Uniform\qp{-1, 1}}$ denote a set of uniformly distributed random
values, defined at each of the mesh points. We set $Y(\vec x)$ to be the
Lagrange interpolant of these random values and define
\begin{gather}
  \label{eq:1d-random-ics}
  \phi_h^0 = \frac{1}{100} Y(\vec x) \AND  \vec v_h^0 \equiv 0
\end{gather}
to be the initial conditions for this test.  Figure
\ref{fig:2d-random-solution-plot} shows solution plots at various
times together with the energy/mass/energy deviation plot. The energy
deviation in this case is a visual representation of the energy
dissipation equality stated in Theorem
\ref{the:fully-discrete-energy}. In this sense, we are defining the
energy deviation for $n\in [0,N-1]$ to be the quantity
\begin{equation}
  \label{eq:energy-dev}
  \begin{split}
    \int_\W W(\phi_h^{n+1}) + \frac{1}{2} \rho(\phi_h^{n+1})\norm{\vec
      v_h^{n+1}}^2 + \frac{\gamma}{2} \norm{\vec q_h^{n+1}}^2 -
    \int_\W W(\phi_h^{n}) + \frac{1}{2} \rho(\phi_h^{n})\norm{\vec
      v_h^{n}}^2 + \frac{\gamma}{2} \norm{\vec q_h^{n}}^2 
    \\
    + k_n
    \Bigg(\int_\W m_r \qp{\nplush{a_h}}^2 - m_j \cA_1\qp{\nplush{a_h},
      \nplush{a_h}} - \eta \cA_2\qp{\nplush{\vec v_h}, \nplush{\vec
        v_h}} \Bigg).
  \end{split}
\end{equation}

Note that the mass is conserved, the energy is monotonically
decreasing and the energy deviation is zero.

\renewcommand{\figscale}{1.0}
\renewcommand{\figwidth}{\textwidth}

\begin{figure}[h!]
  \caption[]
          {
            \label{fig:2d-random-solution-plot} 
            \ref{sec:2d-random} Test 2 -- The solution, $\phi_h$ to
            the quasi-incompressible system with random initial conditions
             at various values of $t$. }
          \begin{center}
            \subfigure[][$t=0$]{
              \includegraphics[scale=\figscale, width=0.4\figwidth]
                              {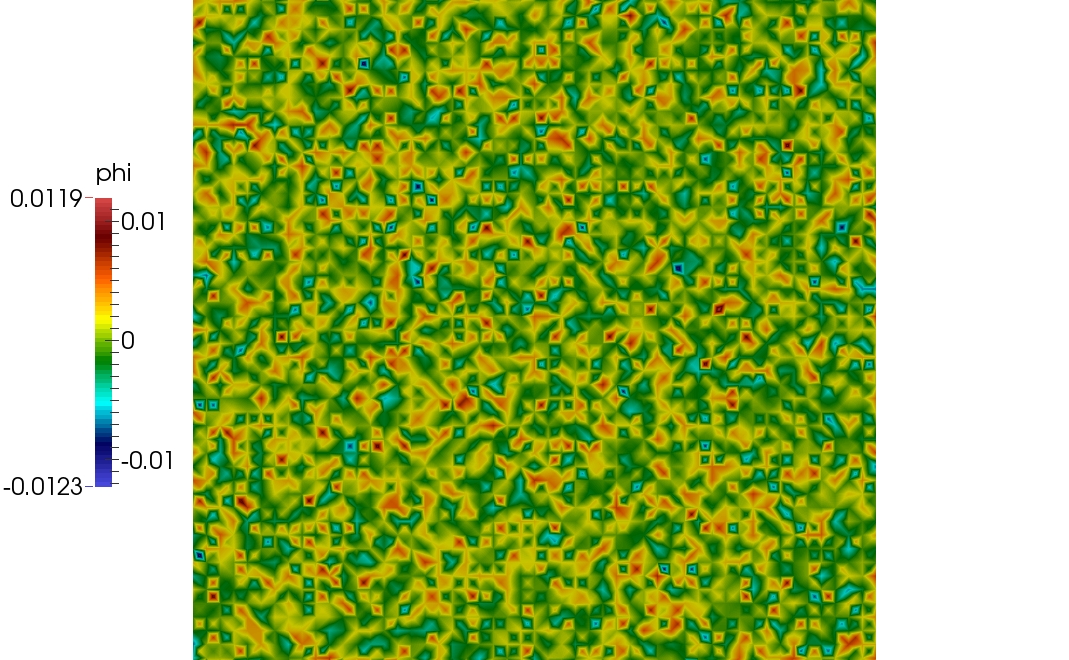}
            }
            \hfill
            \subfigure[][$t=0.05$]{
              \includegraphics[scale=\figscale, width=0.4\figwidth]
                              {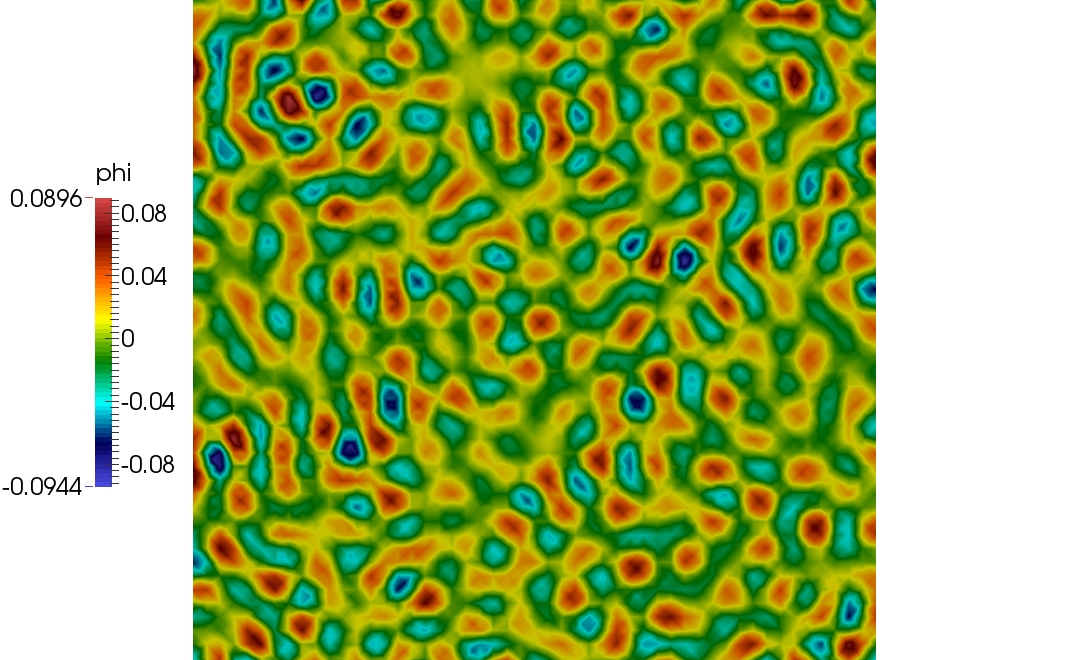}
                             
            }
            \hfill
            \subfigure[][$t=0.16$]{
              \includegraphics[scale=\figscale, width=0.4\figwidth]
              {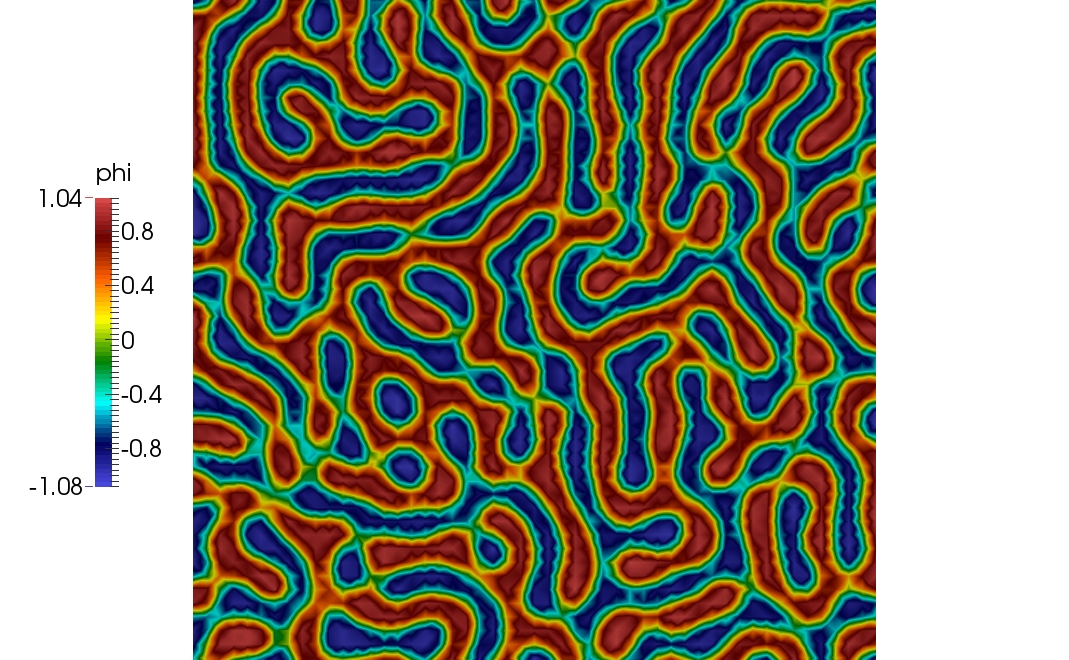}
                             
            }
             \hfill
             \subfigure[][$t=0.3$]{
               \includegraphics[scale=\figscale, width=0.4\figwidth]
               {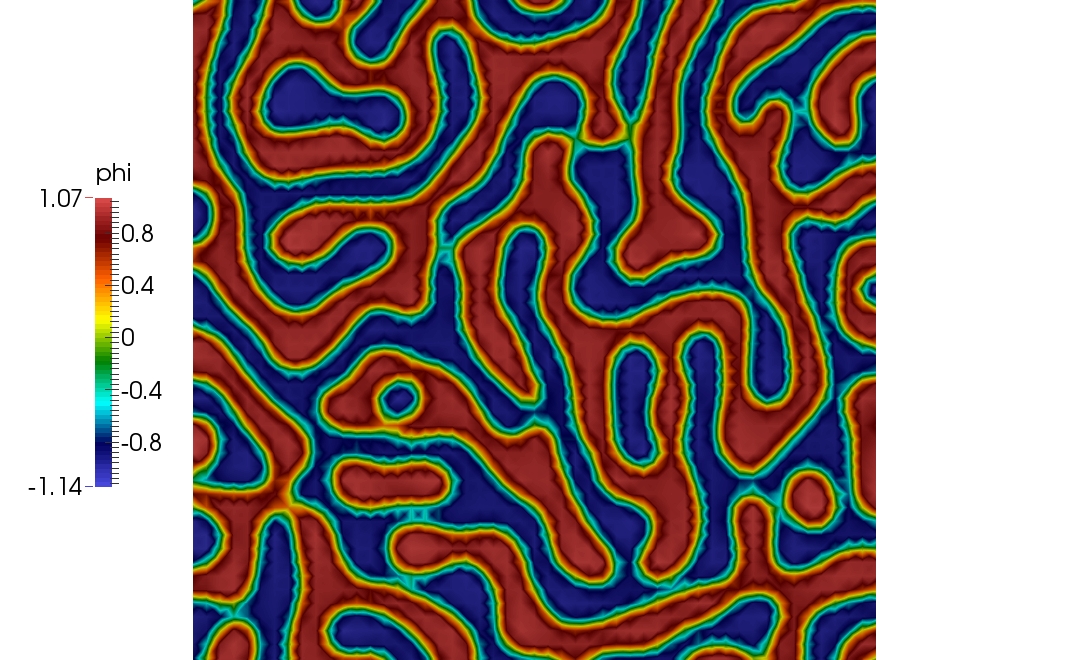}
                             
             }
             \hfill
             \subfigure[][$t=0.5$]{
               \includegraphics[scale=\figscale, width=0.4\figwidth]
               {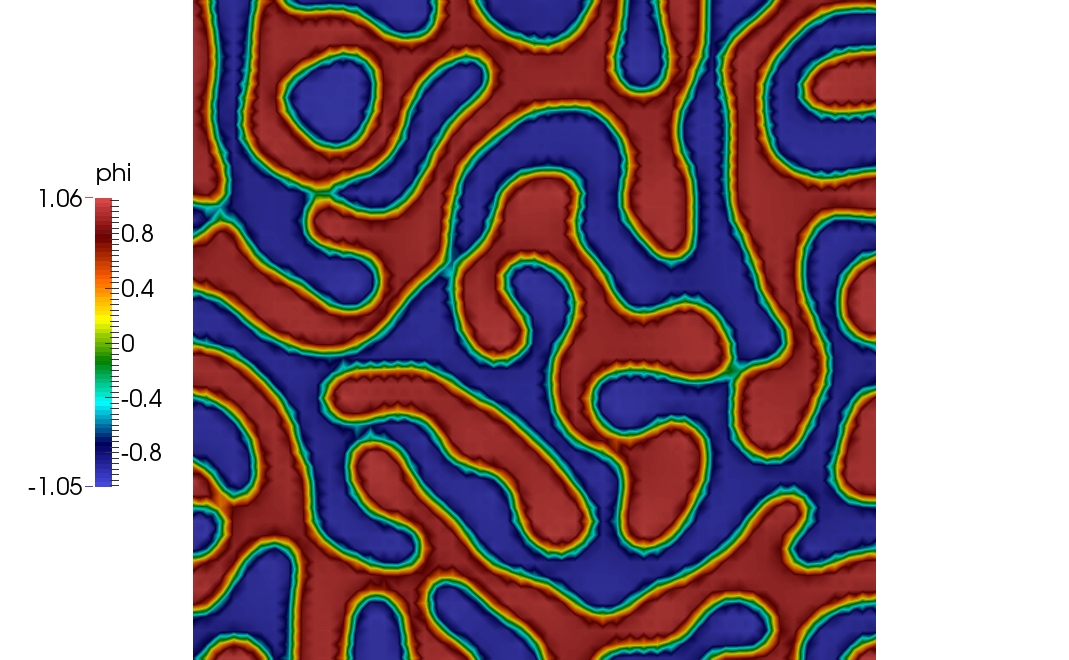}
                             
             }
             \hfill
             \subfigure[][$t=1.$]{
               \includegraphics[scale=\figscale, width=0.4\figwidth]
               {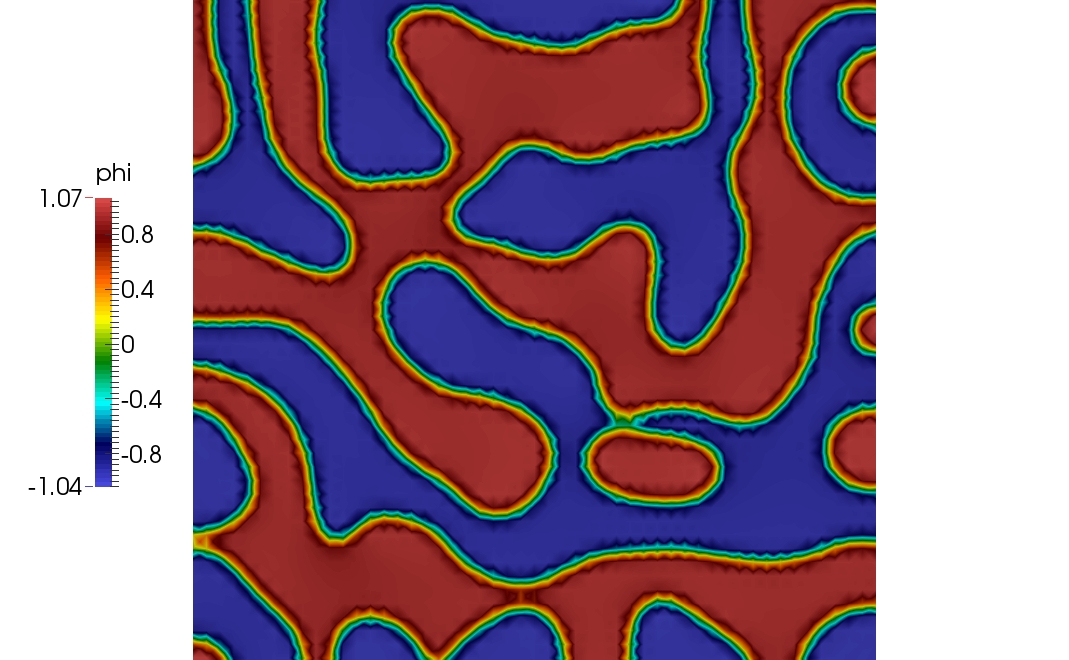}             
             }
              \hfill
            \subfigure[][Conservativity/consistency plot]{
              \includegraphics[scale=\figscale, width=0.7\figwidth]
               {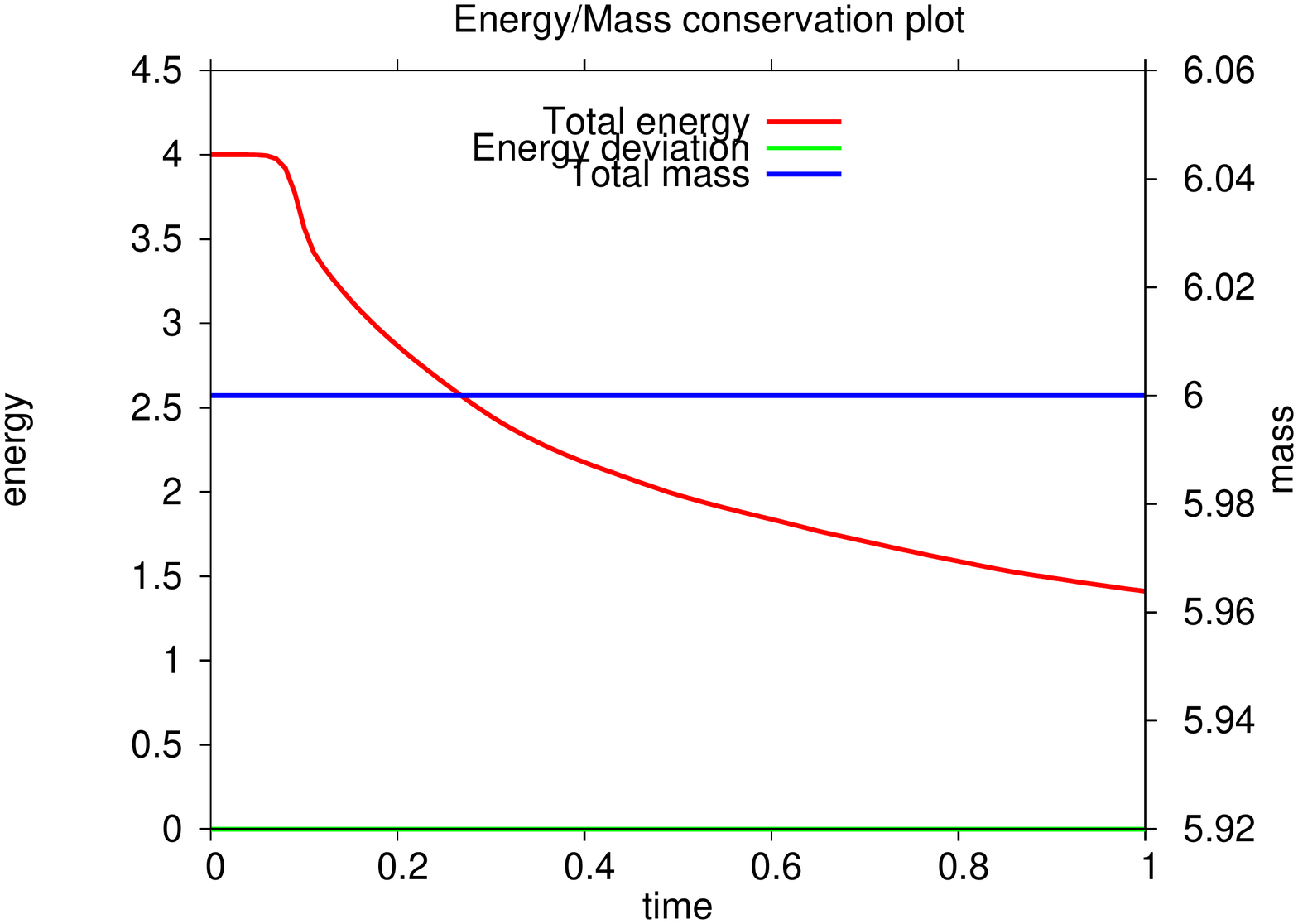}
            }
          \end{center}
\end{figure}

\subsection{Test 3 : 2D - parameter values}
\label{sec:2d-sims}

In this experiment we investigate the effects differing magnitudes of
parameter values have on the dynamics of the system. We vary the
diffusive terms $m_r$ and $m_j$. 

The initial conditions we consider are given by considering $\W =
[0,1]^2$ and defining subsets
\begin{gather}
  \W_1 = \{\vec x : \qp{\abs{x_1 - 1/4}^2 + \abs{x_2 - 1/4}^2} \leq 0.05^2 \}
  \\
  \W_2 = \{\vec x : \qp{\abs{x_1 - 1/4}^2 + \abs{x_2 - 3/4}^2} \leq 0.01^2 \}
  \\
  \W_3 = \{\vec x : \qp{\abs{x_1 - 3/4}^2 + \abs{x_2 - 1/4}^2} \leq 0.01^2 \}
  \\
  \W_4 = \{\vec x : \qp{\abs{x_1 - 3/4}^2 + \abs{x_2 - 3/4}^2} \leq 0.01^2 \},
\end{gather}
and choosing
\begin{gather}
  \label{eq:param-values-ics}
  \phi^0
  = 
  \begin{cases}
    -1 \text{ if } \vec x \in\W_1\cup\W_2\cup\W_3\cup\W_4
    \\
    1 \text{ otherwise }
  \end{cases}
  \qquad
  \vec v = \vec 0.
\end{gather}

Figure \ref{fig:2d-solution-plot} gives some comparitive solution
plots at various times in the simulation. Note that by decreasing the
magnitude of the dissipative terms, the system takes longer to reach a
steady state. The simulation with the smallest values reaches a steady
state at $t\approx 32$. Note that when each simulation reaches a
steady state $\Norm{\vec v_h}_{\leb{\infty}(\W)} \leq 10^{-5}$ which
means that there are no relevant parasitic currents.

\renewcommand{\figscale}{1.0}
\renewcommand{\figwidth}{\textwidth}

\begin{figure}[h]
  \caption[]
          {
            \label{fig:2d-solution-plot} 
            \ref{sec:2d-sims} Test 3 -- The solution, $\phi_h$ to the
            quasi-incompressible system with initial conditions
            (\ref{eq:param-values-ics}) at various values of
            $t$. Notice that there are no parasitic currents appearing
            in the interfacial layer. The velocity tends to zero over
            the entire domain as time increases. }
          \begin{center}
            \subfigure[][$t=0.1$, left $m_j=m_r=1$, middle $m_j=m_r=0.1$, right $m_j=m_r=0.01$]{
              \includegraphics[scale=\figscale, width=0.33\figwidth]
                              {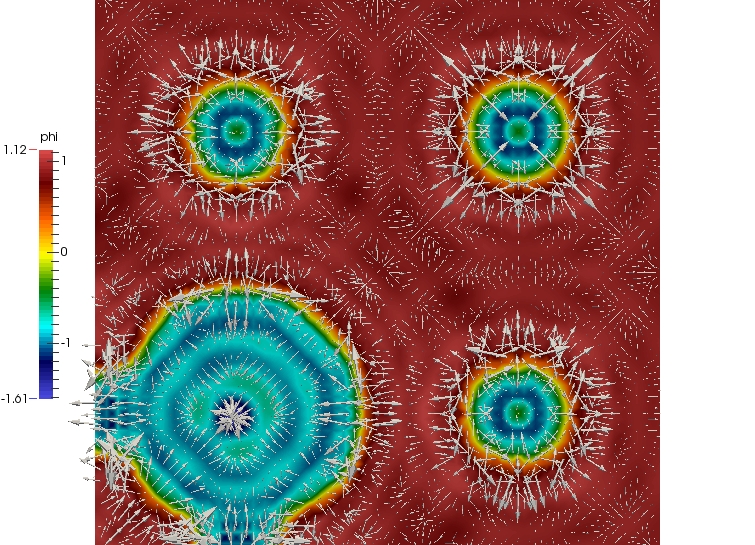}
              \includegraphics[scale=\figscale, width=0.33\figwidth]
                              {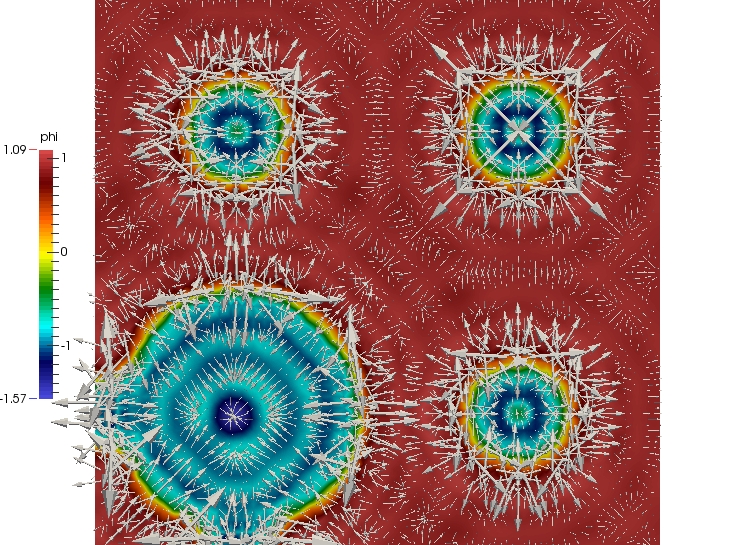}
              \includegraphics[scale=\figscale, width=0.33\figwidth]
                              {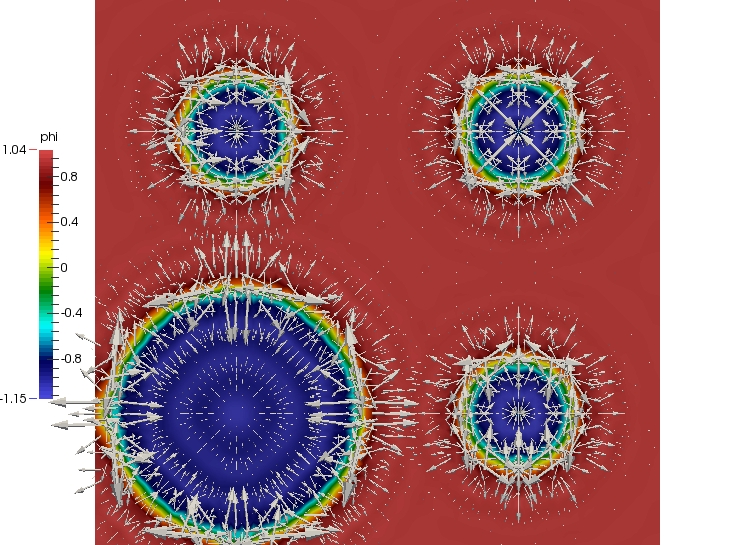}
                             
            }
            \hfill
            \subfigure[][$t=0.25$, left $m_j=m_r=1$, middle $m_j=m_r=0.1$, right $m_j=m_r=0.01$]{
              \includegraphics[scale=\figscale, width=0.33\figwidth]
                              {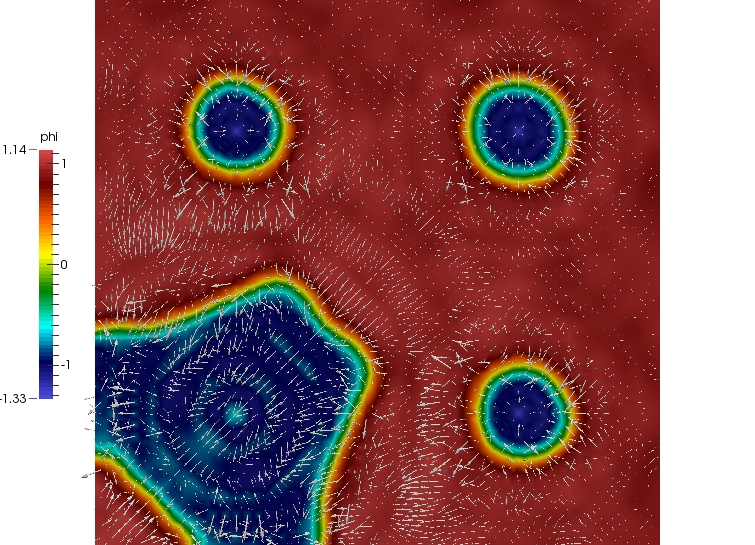}
              \includegraphics[scale=\figscale, width=0.33\figwidth]
                              {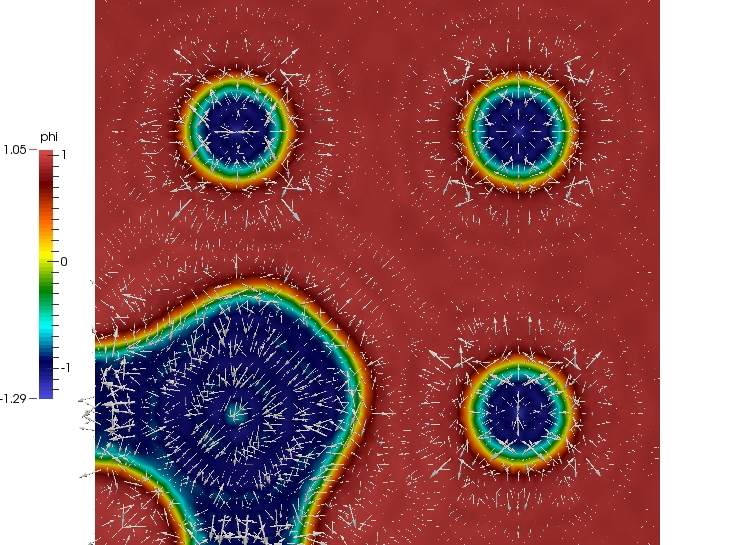}
              \includegraphics[scale=\figscale, width=0.33\figwidth]
                              {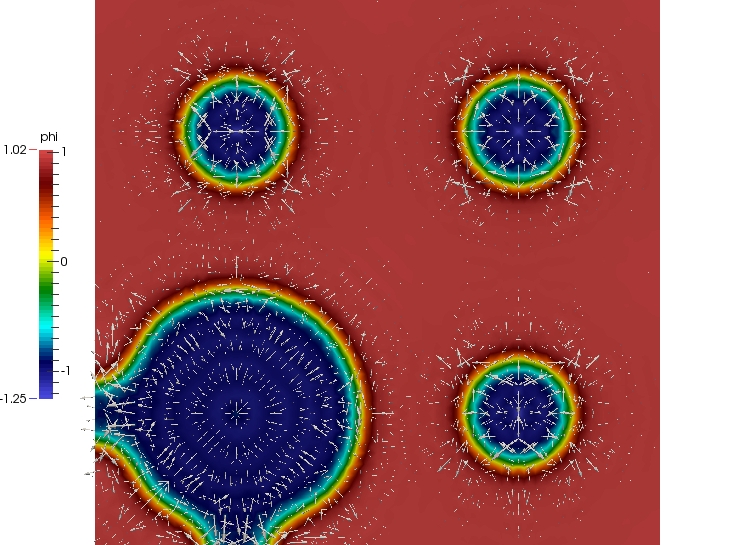}
            }
            \hfill
            \subfigure[][$t=0.5$, left $m_j=m_r=1$, middle $m_j=m_r=0.1$, right $m_j=m_r=0.01$]{
              \includegraphics[scale=\figscale, width=0.33\figwidth]
                              {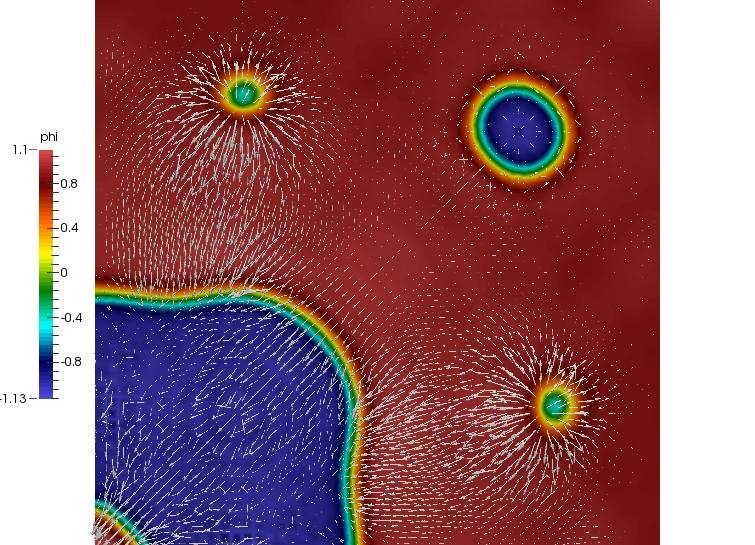}
              \includegraphics[scale=\figscale, width=0.33\figwidth]
                              {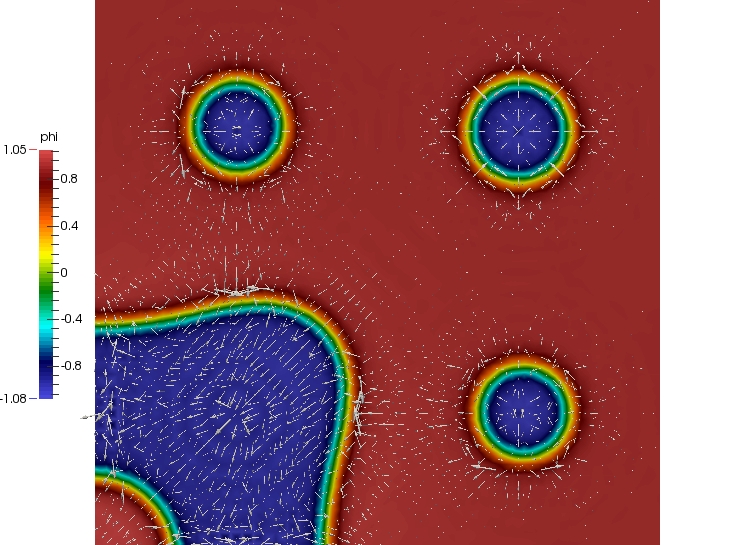}
              \includegraphics[scale=\figscale, width=0.33\figwidth]
                              {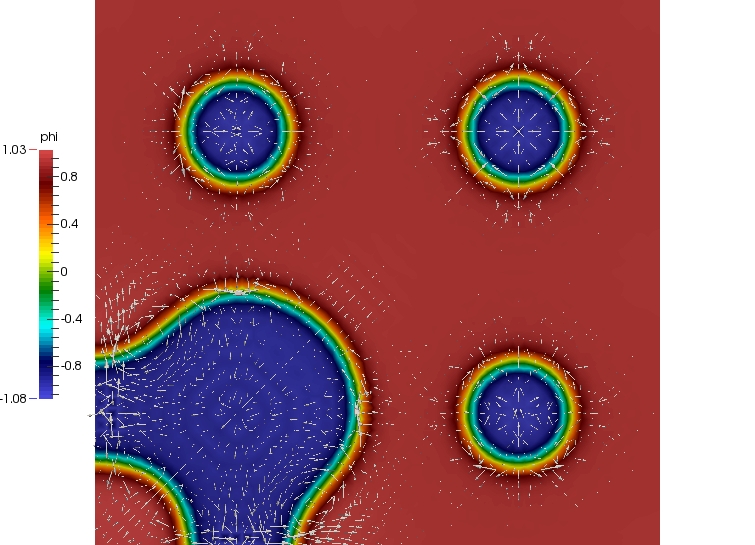}                            
            }
            \hfill
            \subfigure[][$t=1.4$, left $m_j=m_r=1$, middle $m_j=m_r=0.1$, right $m_j=m_r=0.01$]{
              \includegraphics[scale=\figscale, width=0.33\figwidth]
                              {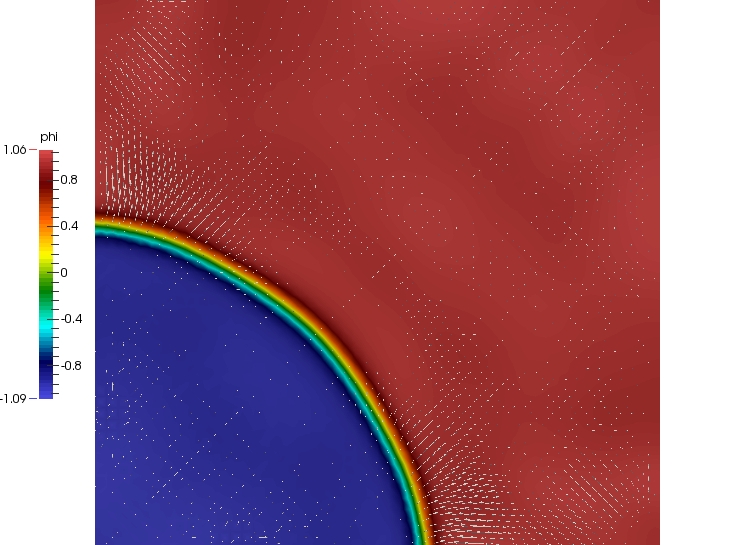}
              \includegraphics[scale=\figscale, width=0.33\figwidth]
                              {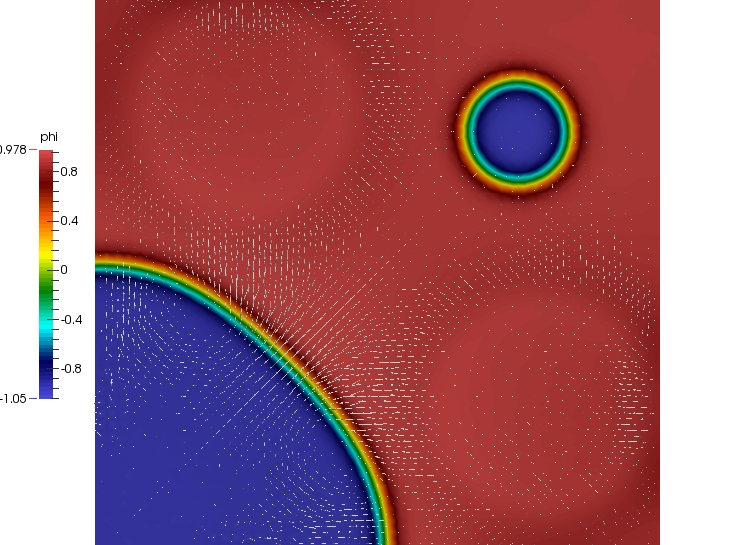}
              \includegraphics[scale=\figscale, width=0.33\figwidth]
                              {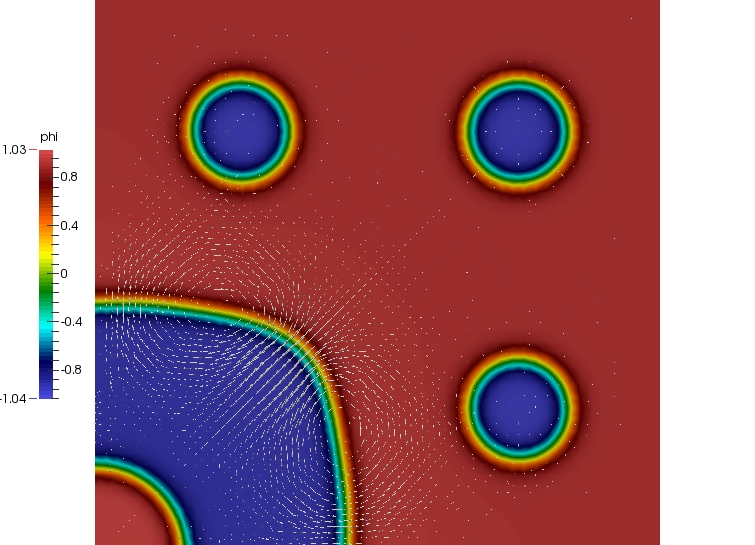}                                                         
            }
            \hfill
            \subfigure[][$t=5$, left $m_j=m_r=1$, middle $m_j=m_r=0.1$, right $m_j=m_r=0.01$]{
              \includegraphics[scale=\figscale, width=0.33\figwidth]
                              {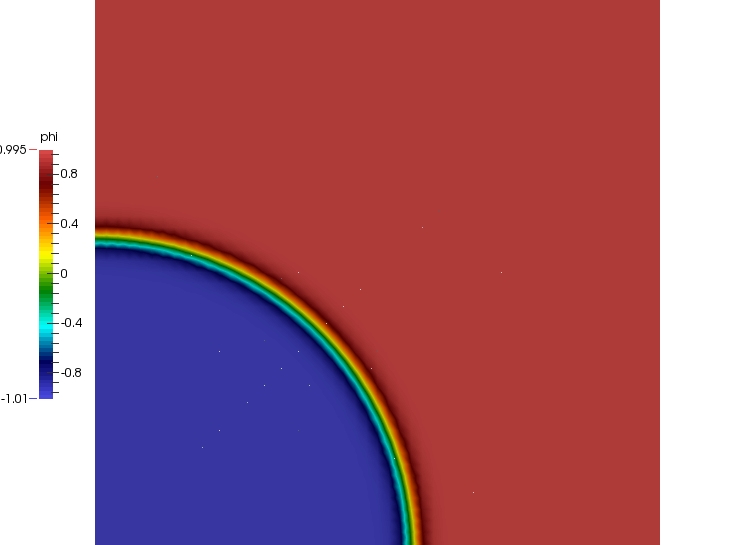}
              \includegraphics[scale=\figscale, width=0.33\figwidth]
                              {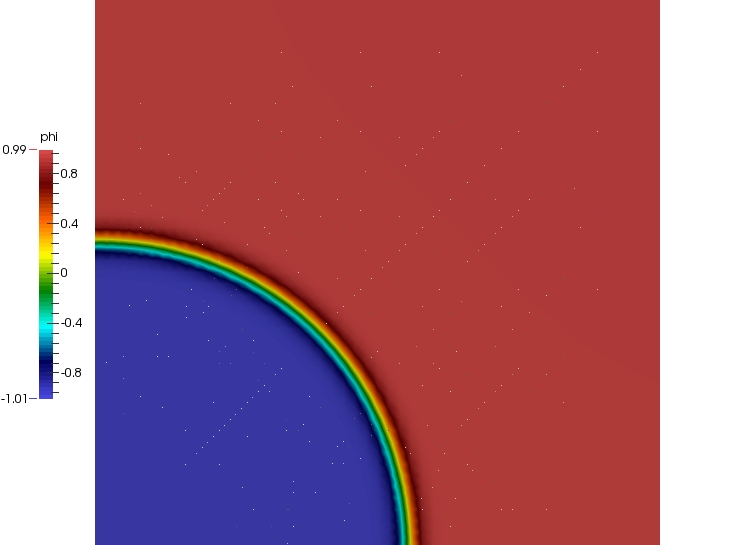}
              \includegraphics[scale=\figscale, width=0.33\figwidth]
                              {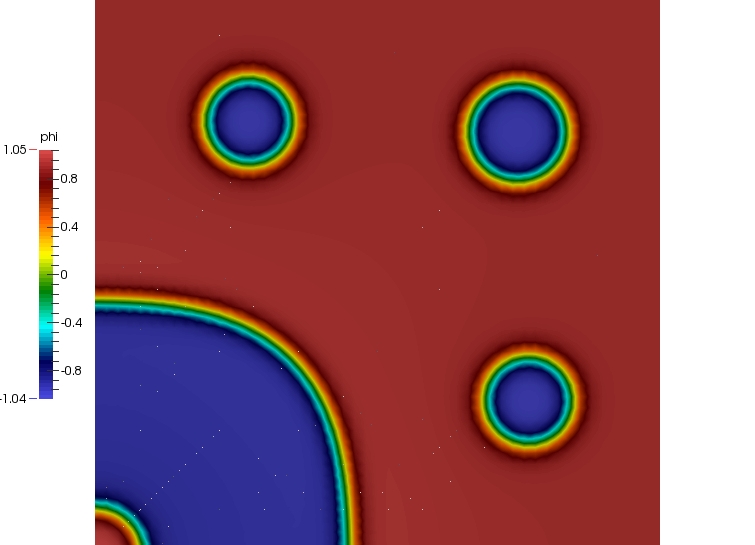}                                                         
                            
            }
            
          \end{center}
\end{figure}

\subsection{Test 4 : 2D rotating coordinate system}
\label{sec:rotate}

Due to the invariance properties of the model
\eqref{eq:quasi-crazy-pde-m1} including the full Navier-Stokes tensor
should we desire computations in a rotating coordinate system the
required changes are very simple. We need only account for inertial or
fictitious forces.  This is in contrast to the model described in
\cite{AGG} which does not behave well with respect to coordinate
changes involving rotating coordinate systems. The fictitious forces
we need to introduce are the Coriolis and the centrifugal force. In
case we consider a planar model problem where the system rotates with
angular velocity $\omega$ around an axis which is perpendicular to the
computational domain then the modified sytem of equations reads
\begin{equation}
  \label{eq:quasi-crazy-pde-coriolis}
  \begin{split}
  \pdt \phi + \div\qp{\phi\vec v} 
  &=
  c_+ \qp{m_j \Delta - m_r}\qp{c_+ \mu(\phi) + c_- \lambda} 
  \\
  \rho(\phi)
  \qp{
    \pdt {\vec v}
    +
    \qp{\Transpose{\vec v} \nabla}\vec v
  }
  +
  \nabla \qp{p(\phi)+\lambda}
  &= 
   \div( \vec \sigma_{NS})
  +
  \gamma \phi \nabla\Delta \phi 
  - \rho(\phi) {\vec \Omega} \times ({\vec \Omega } \times {\vec x}) 
\\
&\qquad - 2\rho(\phi) \vec \Omega \times \vec v
  \\
  \div{\vec v} 
  &=
  c_-
  \qp{m_j \Delta - m_r}\qp{c_+ \mu(\phi) + c_- \lambda} 
  \end{split}
\end{equation}
where $\vec \Omega = \Transpose{(0,0,\omega)}$ and we embed $\vec v $
to $\rR^3$ as $(\vec v; 0)$ for the sake of the vector product.

We now use the original system including the Navier Stokes tensor
(\ref{eq:quasi-crazy-pde-m1}) and energy consistent approximations for
this problem follow our arguments given a standard (signed)
discretisation of the Navier--Stokes tensor. Indeed, the
discretisation is identical to (\ref{eq:fully-discrete-mixed-form})
with the exception of equation
(\ref{eq:fully-discrete-mixed-form})$_2$ which now reads
\begin{equation}
  \label{eq:ns-cor-alternate}
  \begin{split}
    0 &= \int_\W \rho(\nplush{\phi_h}) \frac{\geovec v_h^{n+1} - \vec
      v_h^n}{k_n} \cdot \geovec \Xi + \rho(\nplush{\phi_h}) \qp{
      \qp{\nplush{\geovec v_h} \cdot \nabla }\nplush{\geovec v_h}}
    \cdot \geovec \Xi
    \\
    &\qquad\qquad- \frac{1}{2} \rho(\nplush{\phi_h})
    \nabla\qp{\norm{\nplush{\geovec v_h}}^2} \cdot \geovec \Xi -
    \eta\cA_2\qp{{\nplush{\geovec v_h}},{\geovec \Xi}} + \nabla
    \nplush{b_h} \cdot \geovec \Xi
    \\
    &\qquad\qquad + \frac{\nplush{\phi_h}}{c_+}\nabla(\nplush{a_h} -
    c_- \nplush{b_h}) \cdot \geovec \Xi + \qp{\rho(\nplush{\phi_h})\vec\W \times \qp{\vec \W \times \vec x} }{\cdot \geovec \Xi}
    \\
    &\qquad\qquad +
    \qp{2 \rho(\nplush{\phi_h})\vec \W \times \nplush{\vec v_h}}{\cdot \geovec \Xi}
    \\
    & \qquad + \int_\E \qp{{-}\avg{ \geovec \Xi} \otimes \avg{
        \rho(\nplush{\phi_h}) \nplush{\geovec v_h}}}:
    \tjump{\nplush{\geovec v_h}} \\&\qquad \qquad {+} \frac{1}{2}
    \jump{\norm{\nplush{\geovec v_h}}^2}\cdot \avg{\rho(\nplush{\phi_h})
      \geovec \Xi}
    \\
    &\qquad\qquad - \jump{\nplush{b_h}} \cdot \avg{\geovec \Xi} -
    \frac{1}{c_+}\jump{\nplush{ a_h} - c_-\nplush{ b_h}} \cdot
    \avg{\nplush{\phi_h} \geovec \Xi},
  \end{split}
\end{equation}
where
\begin{equation}
  \label{eq:ns-tensor-disc}
  \begin{split}
  \cA_2\qp{\vec v_h, \vec \Xi} 
  &=
  -
  \int_\W 
  \qp{ \eta_1 - \frac{2}{d}\eta_2}
  \frob{\qp{\div\qp{\vec v_h} \geomat I_d}}{\D\vec \Xi}
  +
  \eta_2 
  \frob{\qp{\D \vec v_h + \Transpose{\D \vec v}}}{\D\vec \Xi}
  \\
  &\qquad+
  \qp{ \eta_1 - \frac{2}{d}\eta_2}
  \int_{\E\cup\partial\W}
  \frob{\qp{\jump{\vec v_h} \geomat I_d}}{\avg{\D\vec \Xi}}
  +
  \frob{\qp{\avg{\div\qp{\vec v_h}}\geomat I_d}}{\tjump{\vec \Xi}}
  \\
  &\qquad+
  \eta_2
  \int_{\E\cup\partial\W}
  \frob{\qp{\tjump{\vec v_h} + \Transpose{\tjump{\vec v_h}}}}{\avg{\D\vec\Xi}}
  +
  \frob{\tjump{\vec \Xi}}{\avg{\qp{\D\vec v_h+\Transpose{\qp{\D\vec v_h}}}}}
  \\&\qquad -
  \int_\E \frac{\sigma}{h} 
  \frob{\tjump{\vec v_h}}{\tjump{\vec \Xi}},
  \end{split}
\end{equation}
represents an interior penalty type discretisation of the
Navier--Stokes tensor which is signed when the penalty parameter
$\sigma$ is chosen large enough.

We also have access to a Lyapanov functional representing the
\emph{energy} of the system. In this case
\begin{equation}
  \begin{split}
    \ddt \bigg( \int_\W W(\phi) &+ \frac{\rho(\phi)}{2} \norm{\vec v}^2 +
    \frac{\gamma}{2} \norm{\nabla \phi}^2 - \w^2 \frac{\rho(\phi)}{2}
      \norm{\vec x}^2 \bigg)
      \\  &  = - \int_\W m_j  \norm{\nabla \qp{c_+
        \mu(\phi) + c_- \lambda}}^2
    + m_r \qp{ c_+ \mu(\phi) + c_- \lambda }^2 + \frob{\D\vec v}{\vec \sigma_{NS}}.
  \end{split}
\end{equation}

Using the arguments presented above it can be shown that the fully
discrete scheme (\ref{eq:fully-discrete-mixed-form}) with
(\ref{eq:fully-discrete-mixed-form})$_2$ replaced by
(\ref{eq:ns-cor-alternate}) satisfies both mass conservation as well
as the following energy dissipation equality
\begin{equation}
  \begin{split}
    \int_\W 
    W(\phi_h^{n+1}) 
    &+
    \frac{1}{2} \rho(\phi_h^{n+1})\norm{\vec v_h^{n+1}}^2
    +
    \frac{\gamma}{2} \norm{\vec q_h^{n+1}}^2
    - \w^2 \frac{\rho(\phi_h^{n+1})}{2}
    \norm{\vec x}^2 
    \\
    &=
      \int_\W 
      W(\phi_h^{n}) 
      +
      \frac{1}{2} \rho(\phi_h^{n})\norm{\vec v_h^{n}}^2
      +
      \frac{\gamma}{2} \norm{\vec q_h^{n}}^2
      - \w^2 \frac{\rho(\phi_h^{n})}{2} \norm{\vec x}^2 
      \\
      &\qquad -
      k_n 
      \Bigg(\int_\W m_r \qp{\nplush{a_h}}^2
      -
      m_j \cA_1\qp{\nplush{a_h}, \nplush{a_h}} 
      \\
      &\qquad\qquad -
      \cA_2\qp{\nplush{\vec v_h}, \nplush{\vec v_h}}
      \Bigg),
    \end{split}
  \end{equation}
  with $\cA_2$ given by (\ref{eq:ns-tensor-disc}).

\begin{figure}[h!]
  \caption[]
          {
            \label{fig:2d-rotate} 
            \ref{sec:rotate} Test 4 -- The solution, $\phi_h$ to
            the quasi-incompressible system with initial conditions
            (\ref{eq:rot-ics}) at various values of $t$. }
          \begin{center}
            \subfigure[][$t=0.01$]{
              \includegraphics[scale=\figscale, width=0.3\figwidth]
                              {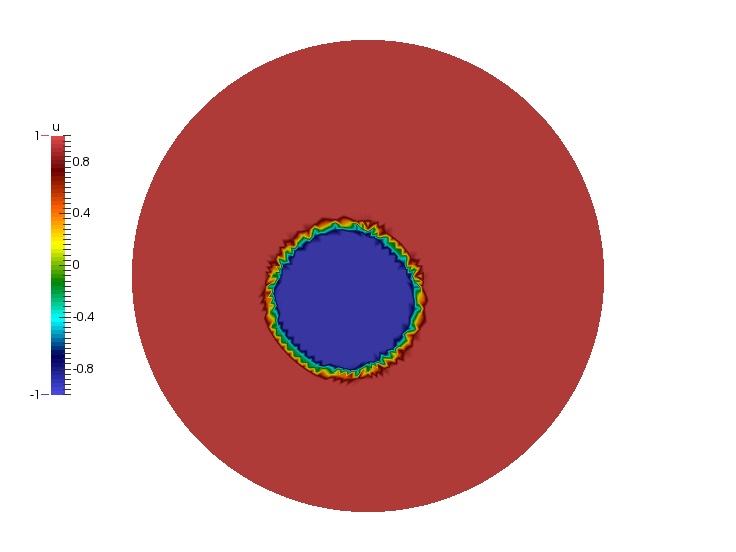}
            }
            \hfill
            \subfigure[][$t=1.75$]{
              \includegraphics[scale=\figscale, width=0.3\figwidth]
                              {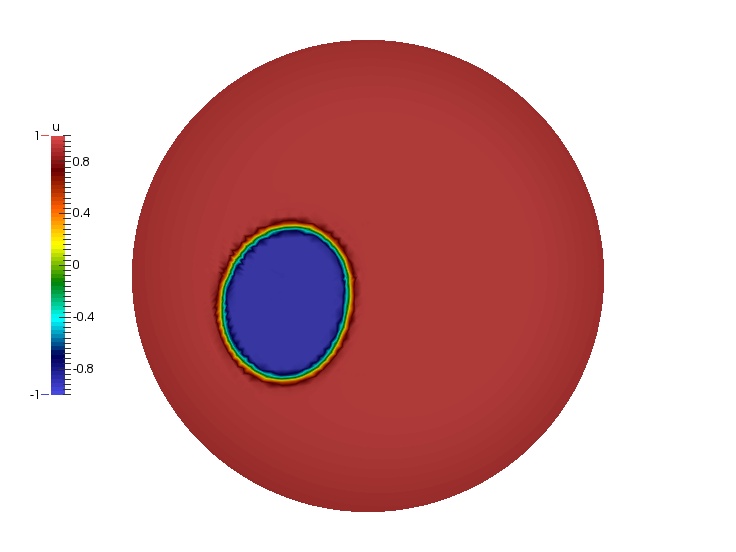}
                             
            }
            \subfigure[][$t=2.61$]{
              \includegraphics[scale=\figscale, width=0.3\figwidth]
                              {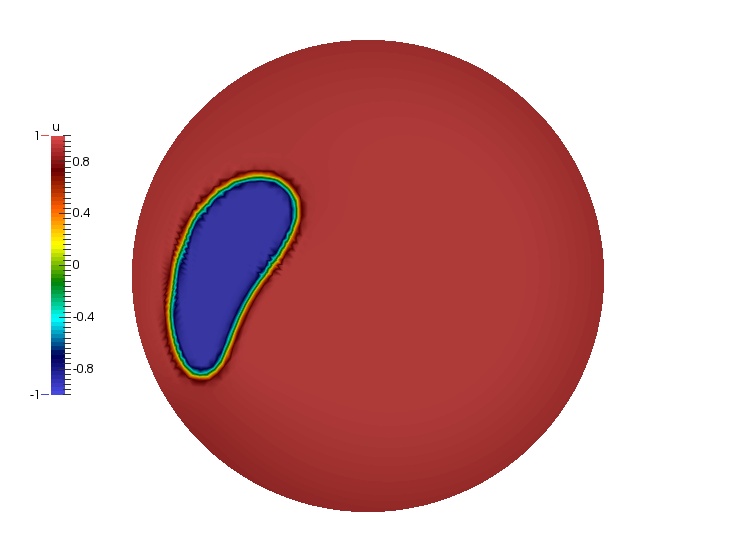}
            }
            \hfill
            \subfigure[][$t=2.91$]{
              \includegraphics[scale=\figscale, width=0.3\figwidth]
                              {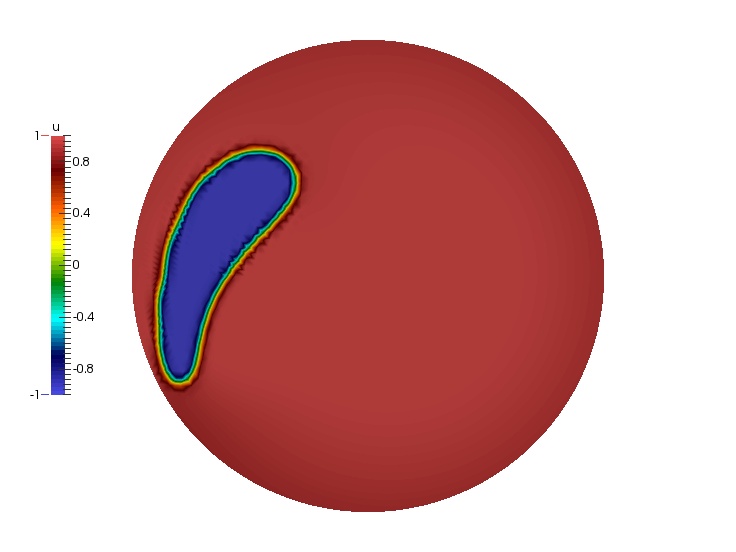}
                             
            }
            \hfill
            \subfigure[][$t=4$]{
              \includegraphics[scale=\figscale, width=0.3\figwidth]
                              {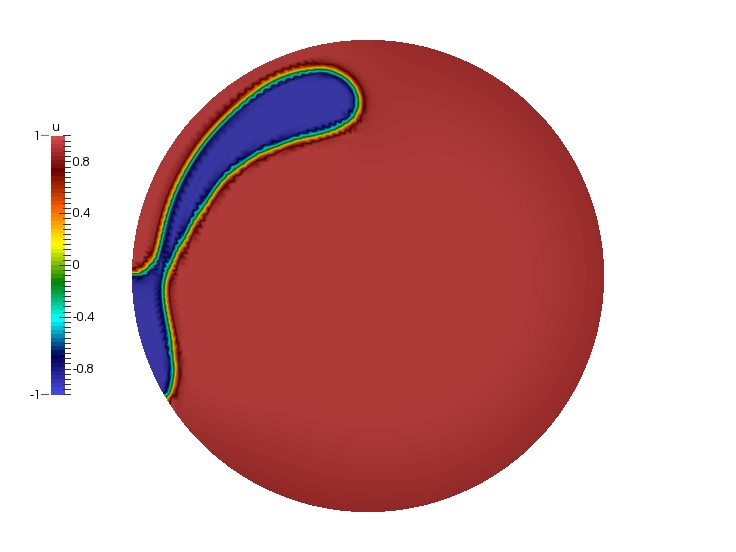}
            }
            \hfill
            \subfigure[][$t=4.5$]{
              \includegraphics[scale=\figscale, width=0.3\figwidth]
                              {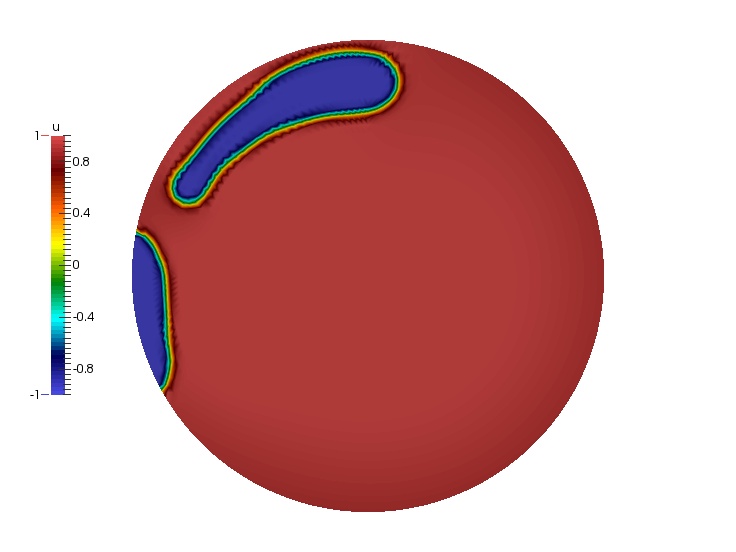}
            }
            \hfill
            \subfigure[][$t=4.98$]{
              \includegraphics[scale=\figscale, width=0.3\figwidth]
                              {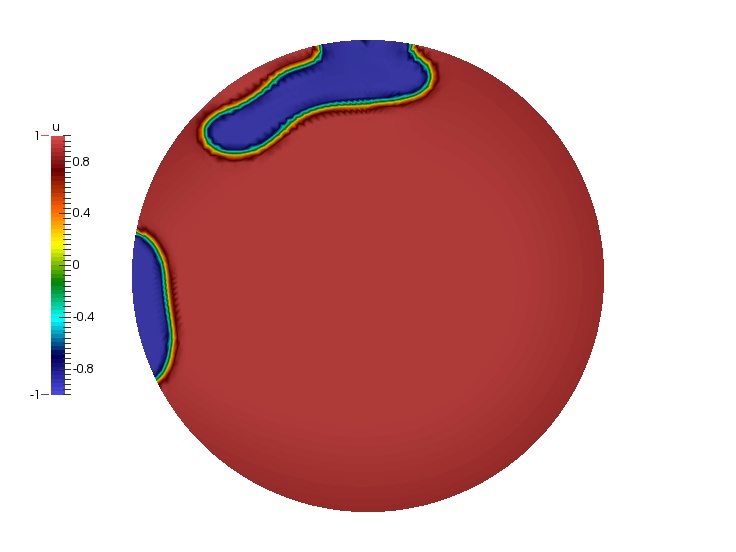}
            }
            \hfill
            \subfigure[][$t=6.52$]{
              \includegraphics[scale=\figscale, width=0.3\figwidth]
                              {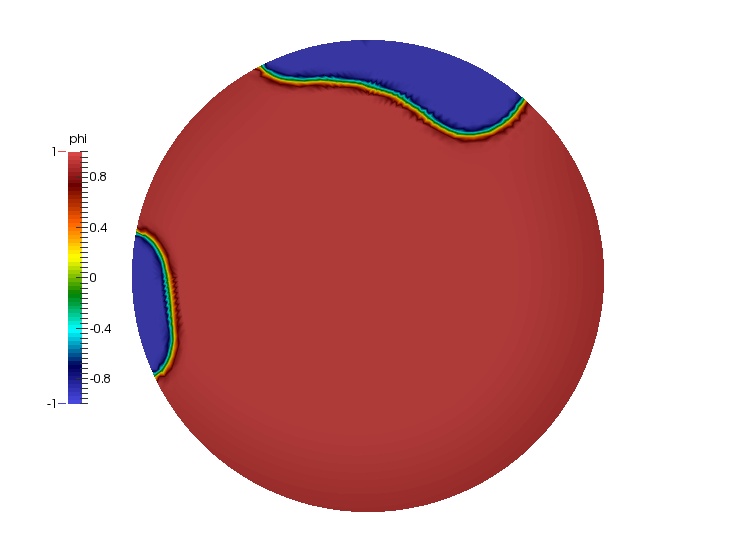}
            }
            \hfill
            \subfigure[][$t=7.64$]{
              \includegraphics[scale=\figscale, width=0.3\figwidth]
                              {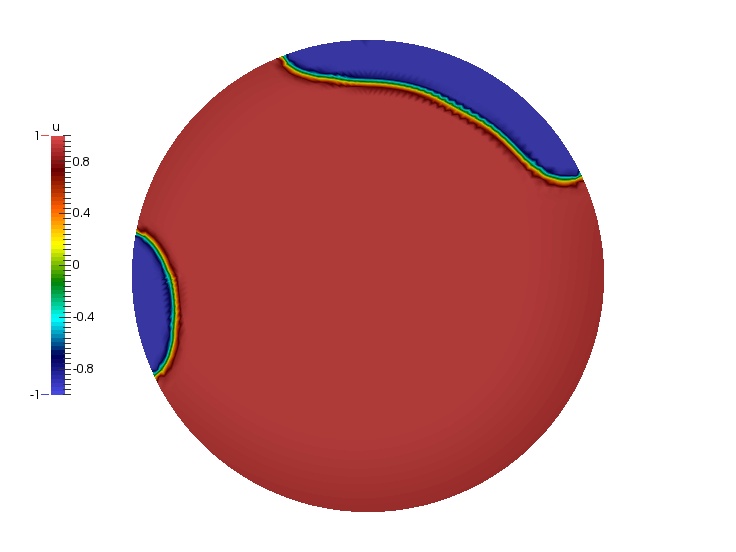}
            }
            \hfill
            \subfigure[][Conservativity/consistency plot]{
              \includegraphics[scale=\figscale, width=0.7\figwidth]
              {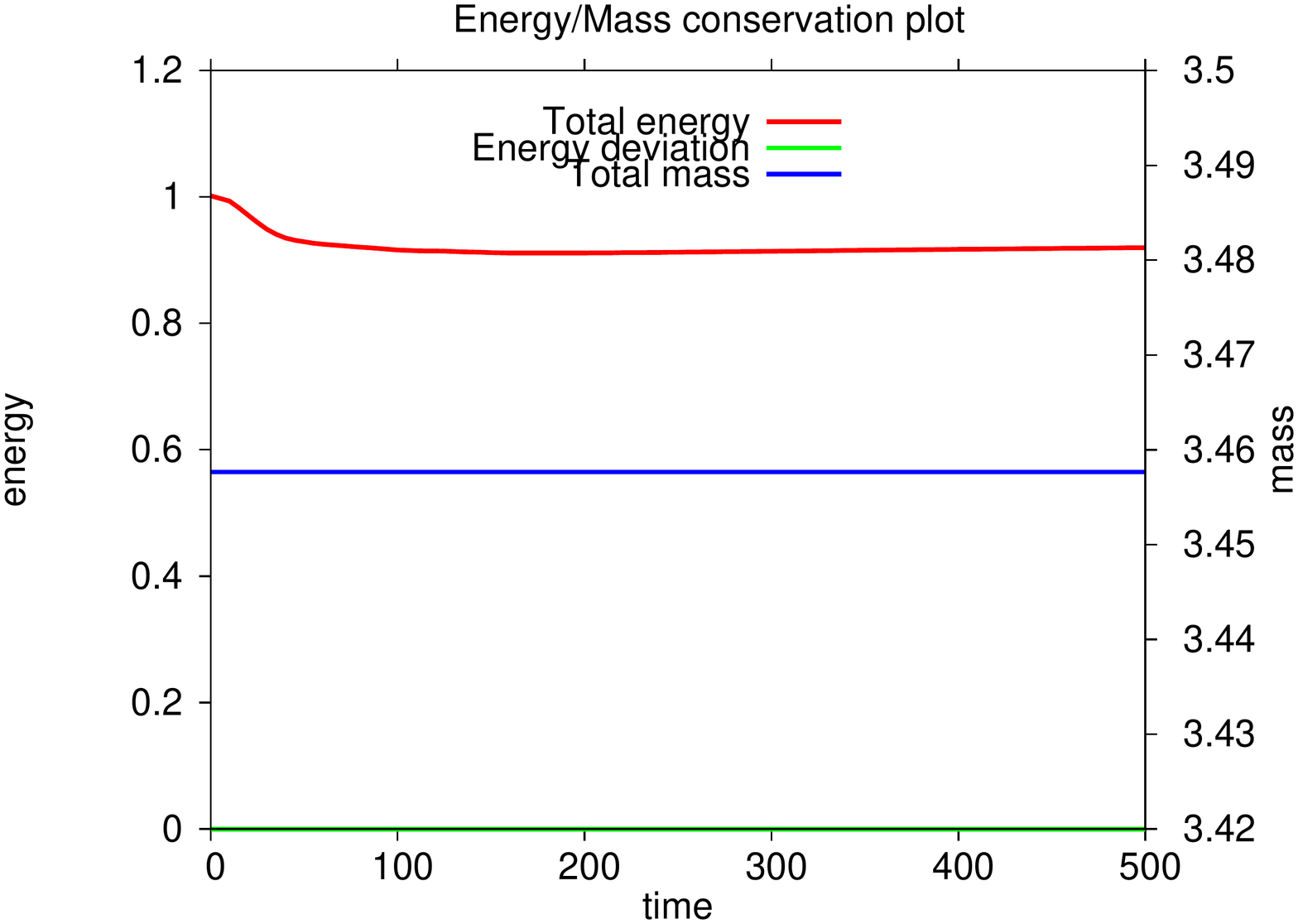}
            }
          \end{center}
\end{figure}

  In Figure \ref{fig:2d-rotate} we illustrate a numerical simulation
  using these principles. We take $\W$ to be a polyhedral
  approximation to the unit circle. We set $\eta_1 = 0.001$ and
  $\eta_2 = 0.005$. We use an initial condition which is a offset
  bubble from the coordinate axis, \ie
\begin{equation}
  \label{eq:rot-ics}
  \phi^0 :=
  \begin{cases}
    -1 \text{ if } \qp{\abs{x_1 + 0.1}^2 + \abs{x_2 + 0.1}^2} \leq 0.1^2
    \\
    1 \text{ otherwise}
  \end{cases},
  \qquad \vec v^0 = \vec 0.
\end{equation}
We show some solution plots at various times as well as the mass/energy
plot.

\subsection{Test 5 : 2D - Rayleigh Taylor instability}
\label{sec:rt}

{
In this test we examine the robustness of the scheme when a denser
fluid lies on top of a lighter one. In this case it is expected that
waves will form over the interface which can give rise to
the formation of plumes. 
}

\begin{figure}[h!]
  \caption[]
          {
            \label{fig:2d-rt} 
            \ref{sec:rt} Test 5 -- The solution, $\phi_h$ to
            the quasi-incompressible system with initial conditions
            (\ref{eq:rt-ics}) at various values of $t$.}
          \begin{center}
            \subfigure[][$t=0.01$]{
              \includegraphics[scale=\figscale, width=0.3\figwidth]
                              {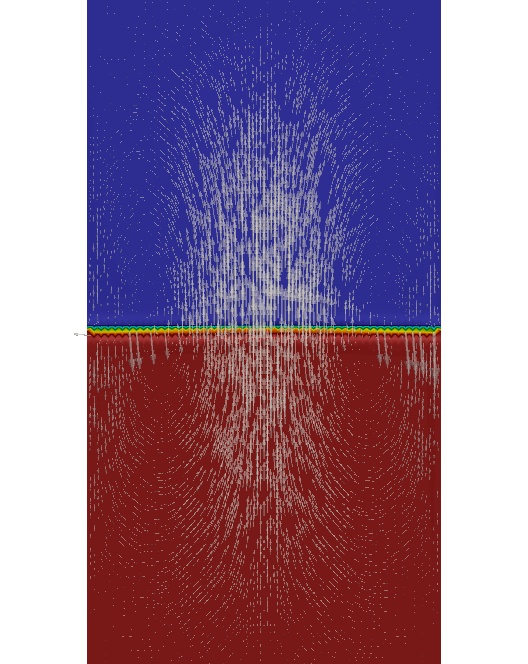}
            }
            \hfill
            \subfigure[][$t=3$]{
              \includegraphics[scale=\figscale, width=0.3\figwidth]
                              {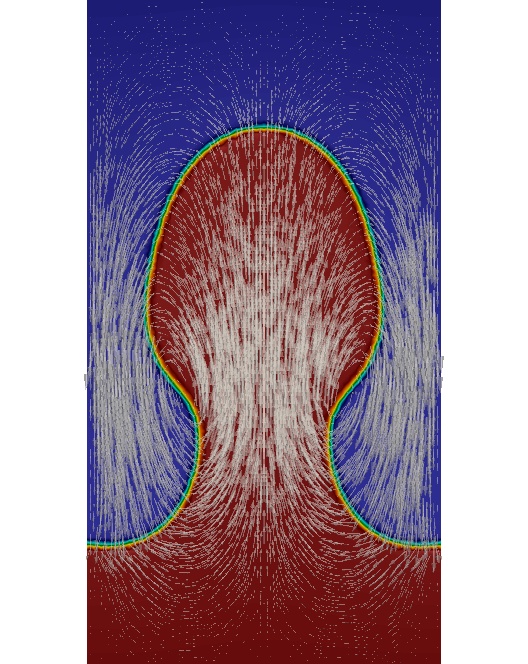}
                             
            }
            \subfigure[][$t=5$]{
              \includegraphics[scale=\figscale, width=0.3\figwidth]
                              {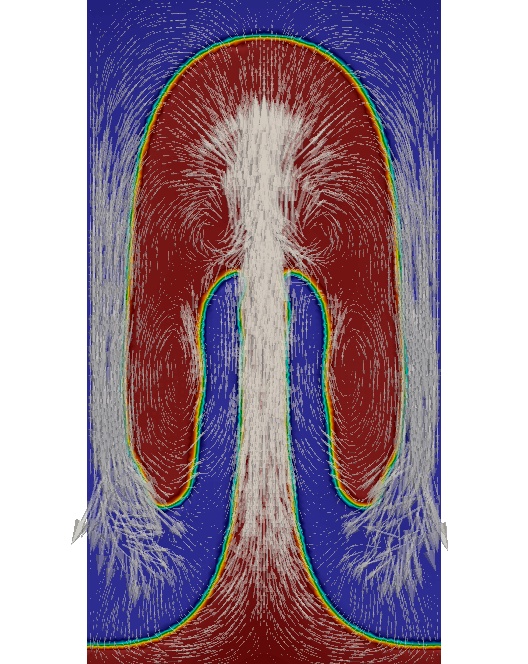}
            }
            \hfill
            \subfigure[][$t=6.65$]{
              \includegraphics[scale=\figscale, width=0.3\figwidth]
                              {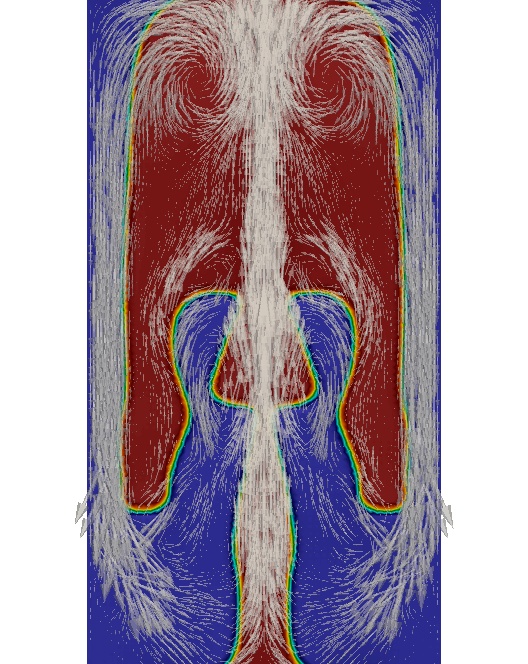}
                             
            }
            \hfill
            \subfigure[][$t=8.1$]{
              \includegraphics[scale=\figscale, width=0.3\figwidth]
                              {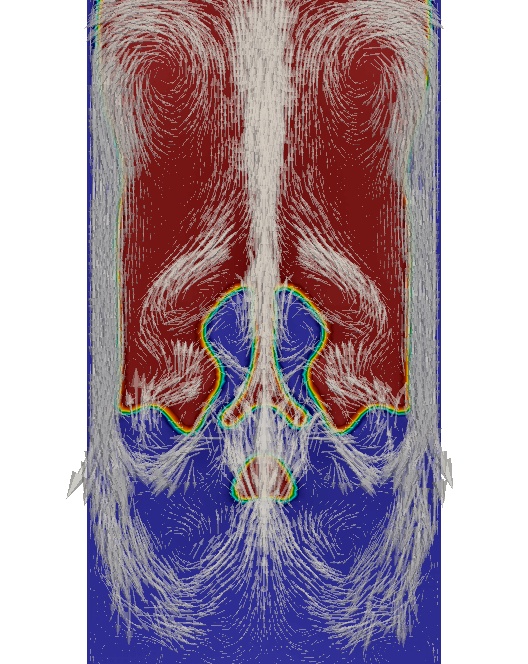}
            }
            \hfill
            \subfigure[][$t=9.11$]{
              \includegraphics[scale=\figscale, width=0.3\figwidth]
                              {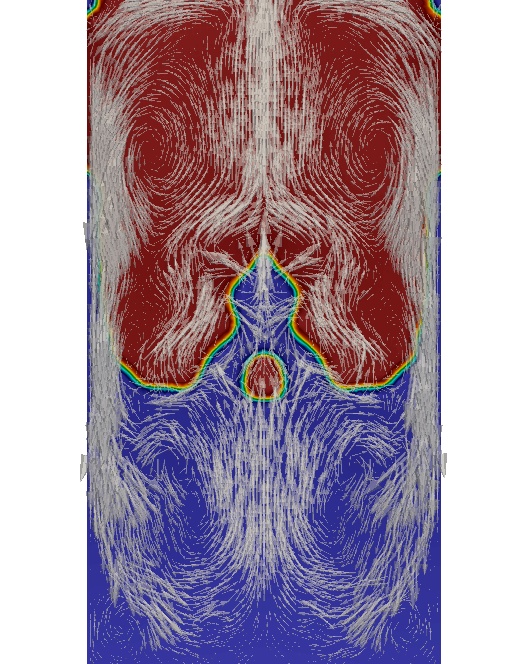}
            }
            \hfill
            \subfigure[][$t=9.54$]{
              \includegraphics[scale=\figscale, width=0.3\figwidth]
                              {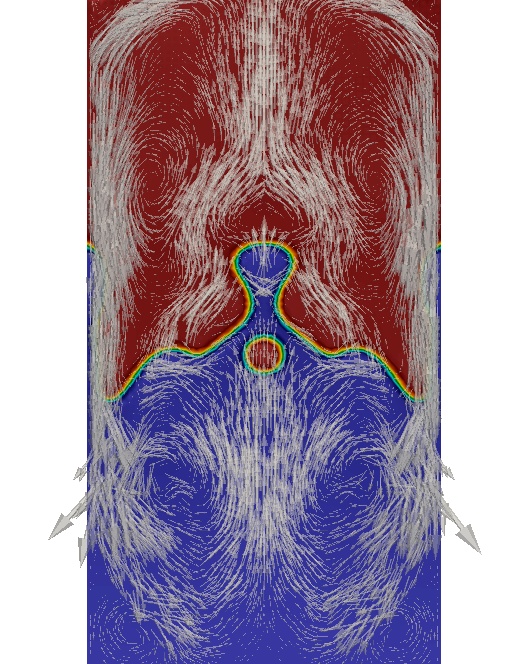}
            }
            \hfill
            \subfigure[][$t=13$]{
              \includegraphics[scale=\figscale, width=0.3\figwidth]
                              {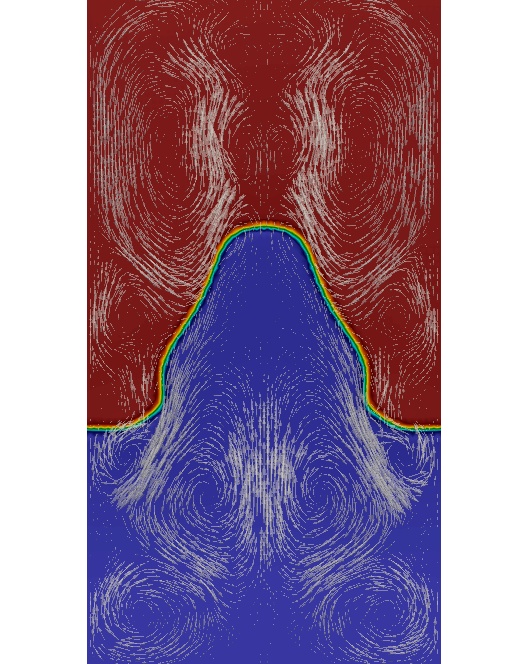}
            }
            \hfill
            \subfigure[][$t=39.95$]{
              \includegraphics[scale=\figscale, width=0.3\figwidth]
                              {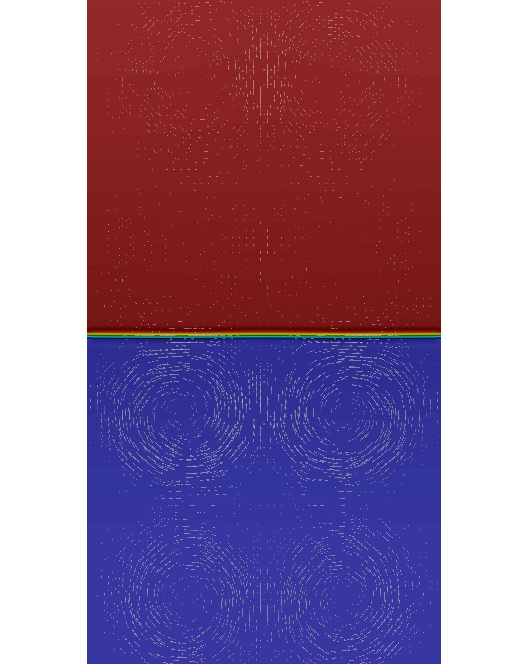}
            }
          \end{center}
\end{figure}

{
We take $\W = [-1,1]\times [-2,2]$ and choose
\begin{equation}
  \label{eq:rt-ics}
  \phi^0 :=
  \begin{cases}
    1 \text{ if } x_2 \leq 0
    \\
    -1 \text{ otherwise}
  \end{cases},
  \qquad \vec v^0 = 
  \bigg(
    \begin{array}{c}
      0
      \\
      {\qp{1+\cos{{\pi x_1}}} \qp{1+\cos{{{\pi x_2}/{2}}}}}/{4}
    \end{array}
    \bigg).
\end{equation}
We also modify (\ref{eq:fully-discrete-mixed-form})$_2$ to take
gravitational effects into account. In this case
(\ref{eq:fully-discrete-mixed-form})$_2$ takes the form
\begin{equation}
  \begin{split}
    0 &=
  \int_\W 
  \rho(\nplush{\phi_h}) \frac{\geovec v_h^{n+1} - \vec v_h^n}{k_n} \cdot \geovec \Xi
  +
    \rho(\nplush{\phi_h}) \qp{ \qp{\nplush{\geovec v_h} \cdot \nabla }\nplush{\geovec v_h}} \cdot \geovec \Xi
    \\
    &\qquad\qquad-
    \frac{1}{2} \rho(\nplush{\phi_h}) \nabla\qp{\norm{\nplush{\geovec v_h}}^2} \cdot \geovec \Xi
    -
    \eta\cA_2\qp{{\nplush{\geovec v_h}},{\geovec \Xi}}
    +
    \nabla \nplush{b_h} \cdot \geovec \Xi
    \\
    &\qquad\qquad +
    \frac{\nplush{\phi_h}}{c_+}\nabla(\nplush{a_h} - c_- \nplush{b_h}) \cdot \geovec \Xi
    +
    \rho(\nplush\phi)\geovec g  \cdot \geovec \Xi
    \\
    & \qquad +
    \int_\E
    \qp{{-}\avg{ \geovec \Xi} \otimes \avg{ \rho(\nplush{\phi_h}) \nplush{\geovec v_h}}}: \tjump{\nplush{\geovec v_h}} 
    \\&\qquad \qquad {+}
    \frac{1}{2} \jump{\norm{\nplush{\geovec v_h}}^2}\cdot \avg{\rho(\nplush{\phi_h}) \geovec \Xi}
    \\ 
    &\qquad\qquad
    -
    \jump{\nplush{b_h}} \cdot \avg{\geovec \Xi}
    -
    \frac{1}{c_+}\jump{\nplush{ a_h} - c_-\nplush{ b_h}} \cdot \avg{\nplush{\phi_h} \geovec \Xi},
    \end{split}
\end{equation}
where $\geovec g = \Transpose{\qp{0, 0.01}}$ is a gravitational
constant. In Figure \ref{fig:2d-rt} we show results from a numerical
experiment with the initial conditions given in (\ref{eq:rt-ics}).  }

\begin{Rem}[guaranteeing positivity of $\rho(\phi)$ and solvability of
  the numerical scheme]
  \label{rem:positive-density}
  {The energy dissipation equality of the numerical scheme
    given in Theorem \ref{the:fully-discrete-energy} gives us no
    information on the solvability of the discrete scheme. In
    addition, the positivity of the density $\rho(\phi)$ is not
    guaranteed. Numerically, for low denisty ratios, like those in
    tests 1--5 where $\rho_2/\rho_1 = 2$, positivity and solvability is
    observed. However, for higher density ratios, this is no longer
    the case. To overcome this difficulty, there are at least three
    possibilities:
  }

  {
    The first is to use a different energy density, which penalises
    values of $\phi$ outside the interval $[-1,1]$. To that end, we
    introduce
    \begin{equation}
      \label{eq:modified-energy}
      W(\phi) = (1+\phi)^2(1-\phi)^2 + A\qp{\qp{\phi - 1 + \norm{\phi - 1}}^2 + \qp{-\phi - 1 + \norm{-\phi-1}}^2},
    \end{equation}
    where $A$ is a large parameter chosen relative to the density
    ratio $\rho_2/\rho_1$ to ensure the density is positive. From a
    modelling point of view, the energy density $W$ is purely
    artificial and thus can be chosen reasonably freely. }

  { The second approach is to use a cutoff of the density
    function as detailed in \cite{Gru}. The main idea is to use the
    densities of the pure phases when $\phi\not\in [-1,1]$.
  }

  { The third approach is to modify the mobilities such that
    they are functions of $\phi$ that are degenerate when $\phi\not\in
    [-1,1]$ in a similar light to \cite{GruRumpf}.
  }

  { The first approach fits into the analytical framework
    developed in this contribution, the second and third do not.  As
    such, we will not persue the case of denisty cutoff functions or
    nonconstant mobilities further but we believe that our results are
    extendable to these cases.}
\end{Rem}

\subsection{Test 6 : 1D - High density ratios}
\label{sec:1d-random}

\begin{figure}[h!]
  \caption[]
          {
            \label{fig:1d-dens-rat-2} 
            \ref{sec:1d-random} Test 6 -- The solution, $\phi_h$ to
            the quasi-incompressible system, using the modified double
            well in (\ref{eq:modified-energy}) with $A =
            (\rho_1 / \rho_2)^2$, with initial conditions
            (\ref{eq:1d-random-ics}) at various values of $t$. In this
            case $\rho_2 / \rho_1 = 2$ and $\max{\phi} = 1.2175$ hence
            $\rho(\phi) > 0$ for all time.}
          \begin{center}
            \subfigure[][$t=0$]{
              \includegraphics[scale=\figscale, width=0.4\figwidth]
                              {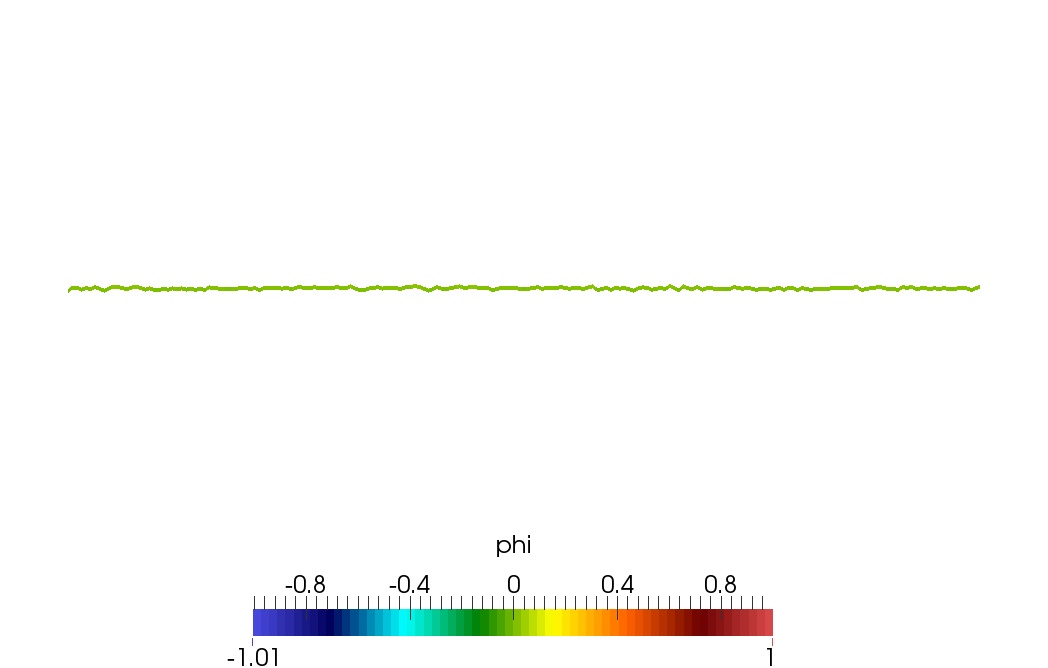}
            }
            \hfill
            \subfigure[][$t=0.09$]{
              \includegraphics[scale=\figscale, width=0.4\figwidth]
                              {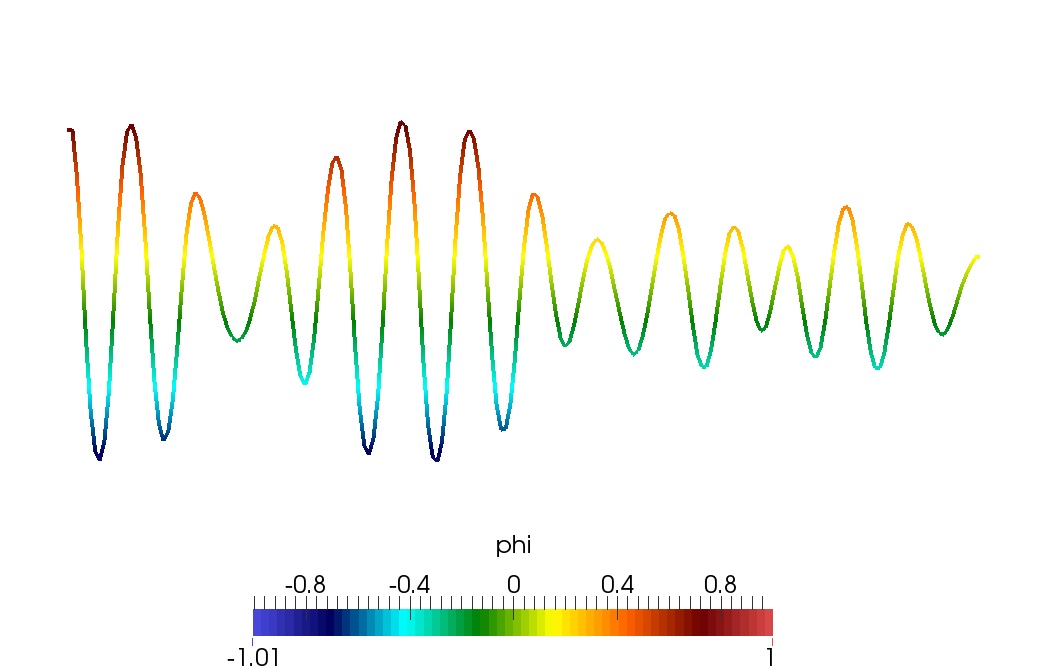}
                             
            }
            \hfill
            \subfigure[][$t=0.43$]{
              \includegraphics[scale=\figscale, width=0.4\figwidth]
              {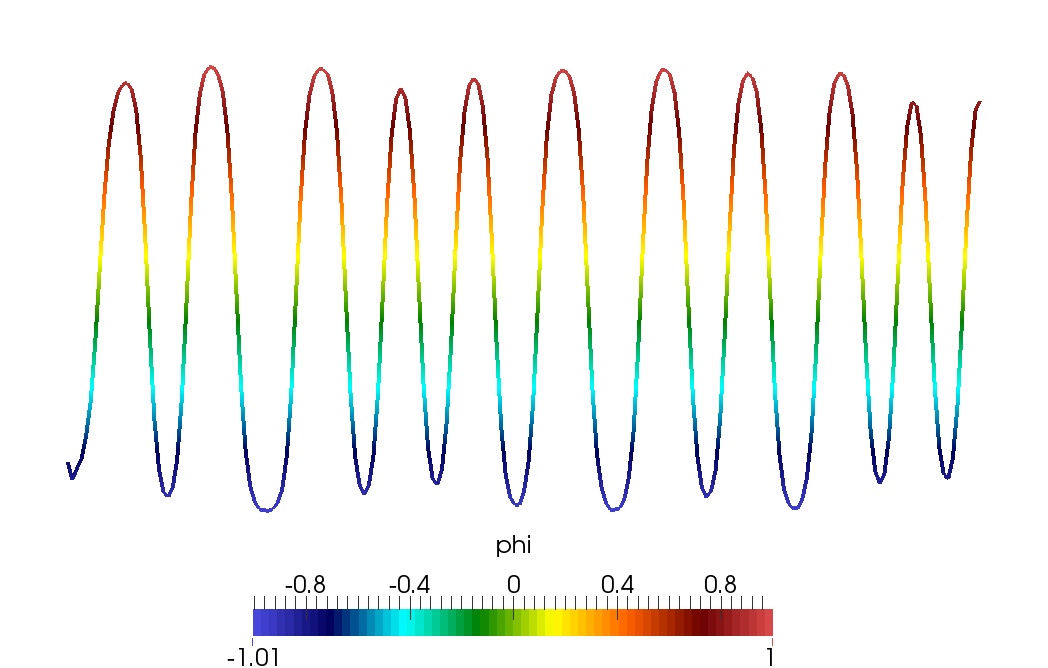}
                             
            }
            \hfill
            \subfigure[][$t=5$]{
              \includegraphics[scale=\figscale, width=0.4\figwidth]
              {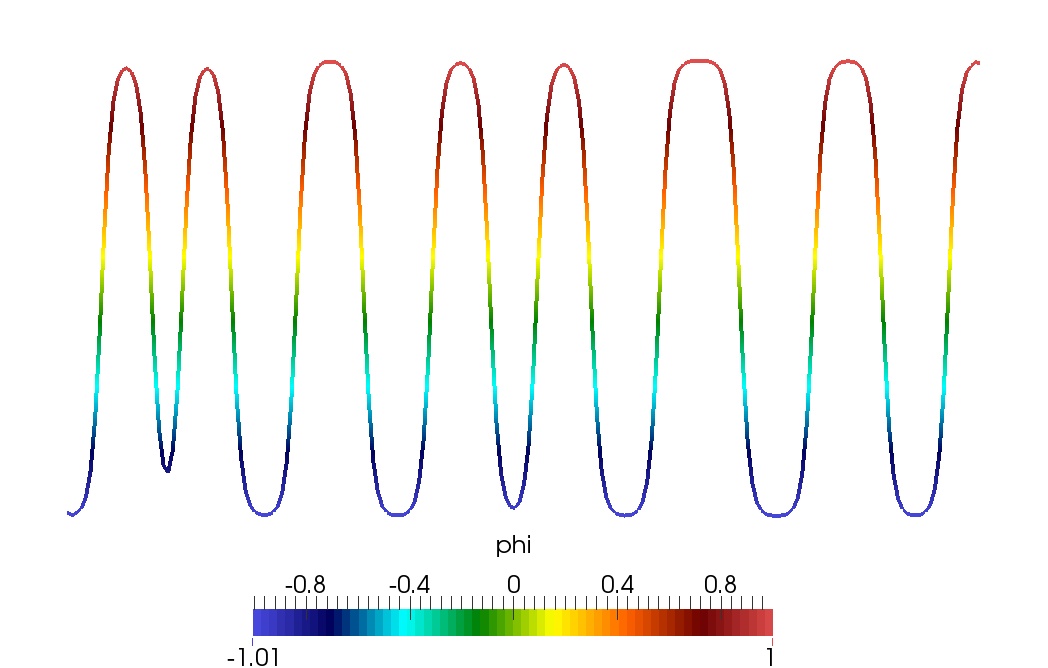}
                             
            }
            \hfill
            \subfigure[][$t=10$]{
              \includegraphics[scale=\figscale, width=0.4\figwidth]
              {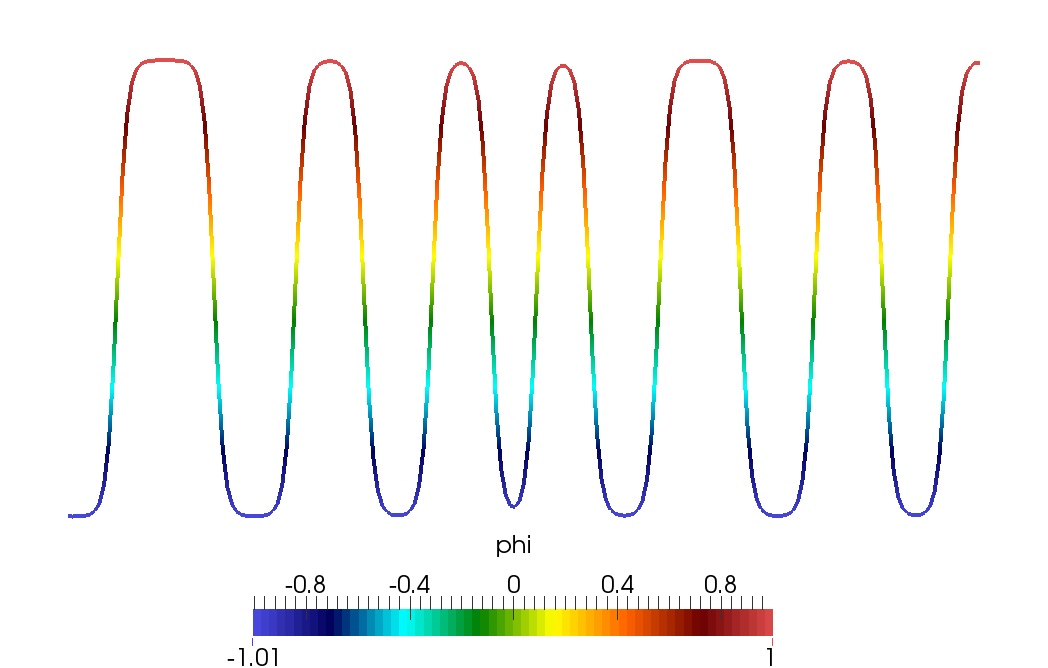}
                             
            }
            \hfill
            \subfigure[][$t=100$]{
              \includegraphics[scale=\figscale, width=0.4\figwidth]
              {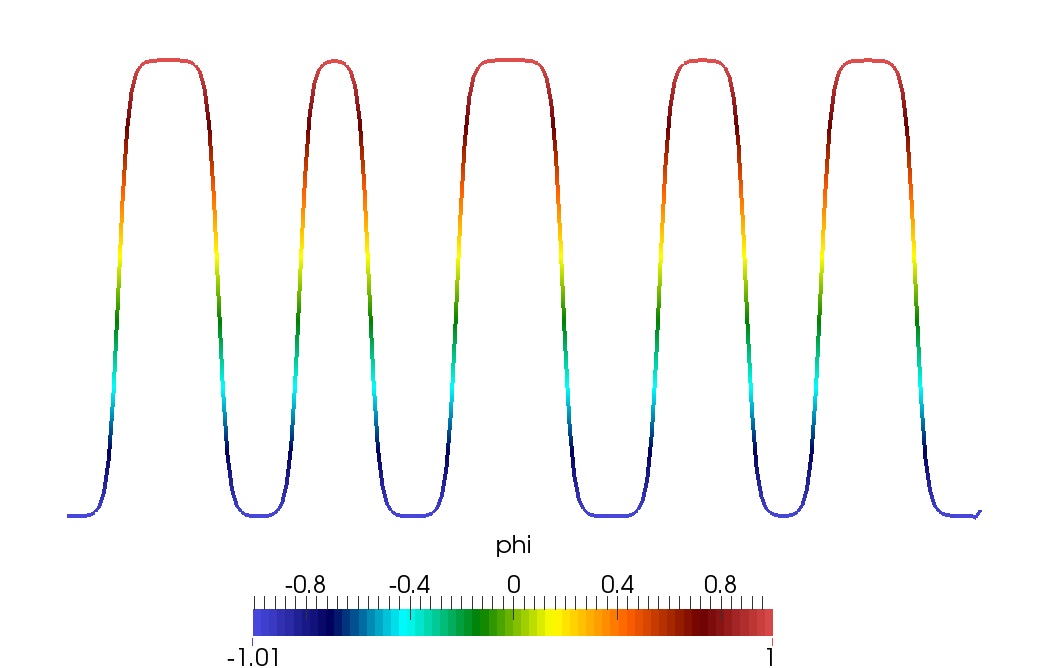}
                             
            }
            \hfill
            \subfigure[][Conservativity/consistency plot]{
              \includegraphics[scale=\figscale, width=0.7\figwidth]
              {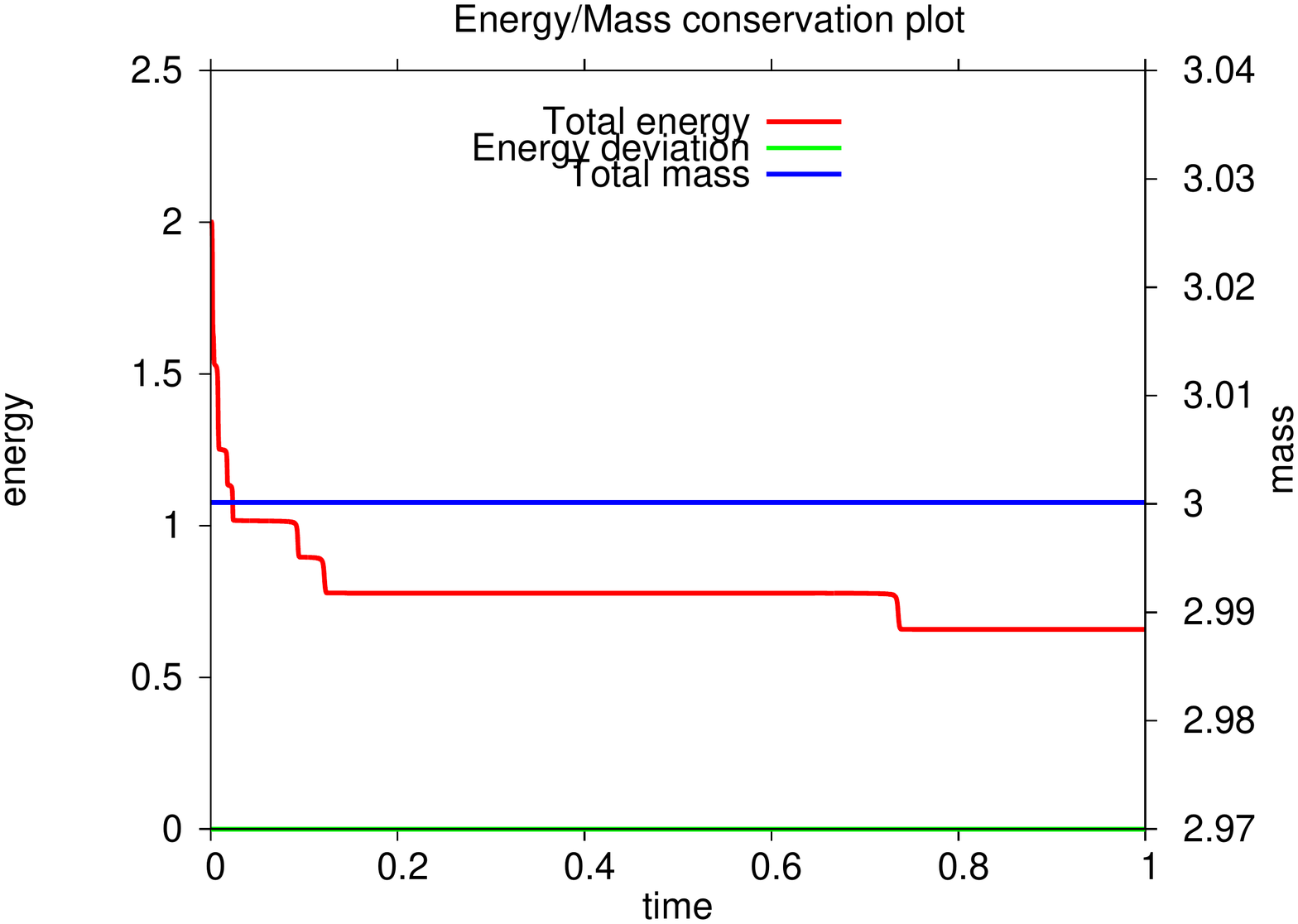}}
          \end{center}
\end{figure}

\begin{figure}[h!]
  \caption[]
          {
            \label{fig:1d-dens-rat-10} 
            \ref{sec:1d-random} Test 6 -- As Figure
            \ref{fig:1d-dens-rat-2} but in this case $\rho_2 / \rho_1
            = 10$ and $\max{\phi} = 1.0233$ hence $\rho(\phi) > 0$ for
            all time.}
          \begin{center}
            \subfigure[][$t=0$]{
              \includegraphics[scale=\figscale, width=0.4\figwidth]
                              {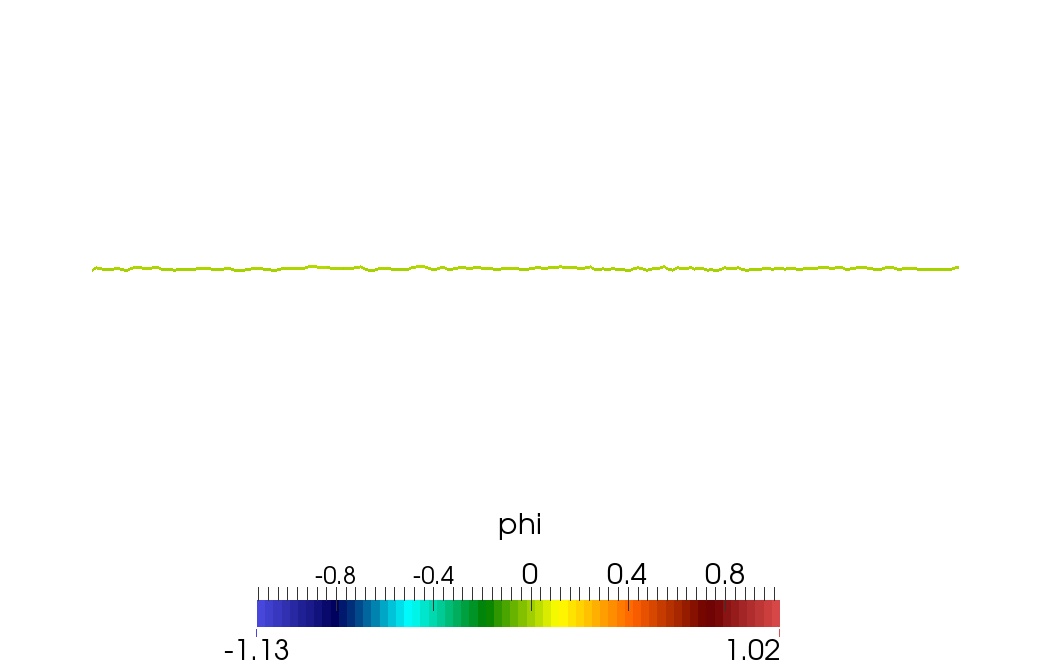}
            }
            \hfill
%
            \hfill
            \subfigure[][$t=0.43$]{
              \includegraphics[scale=\figscale, width=0.4\figwidth]
              {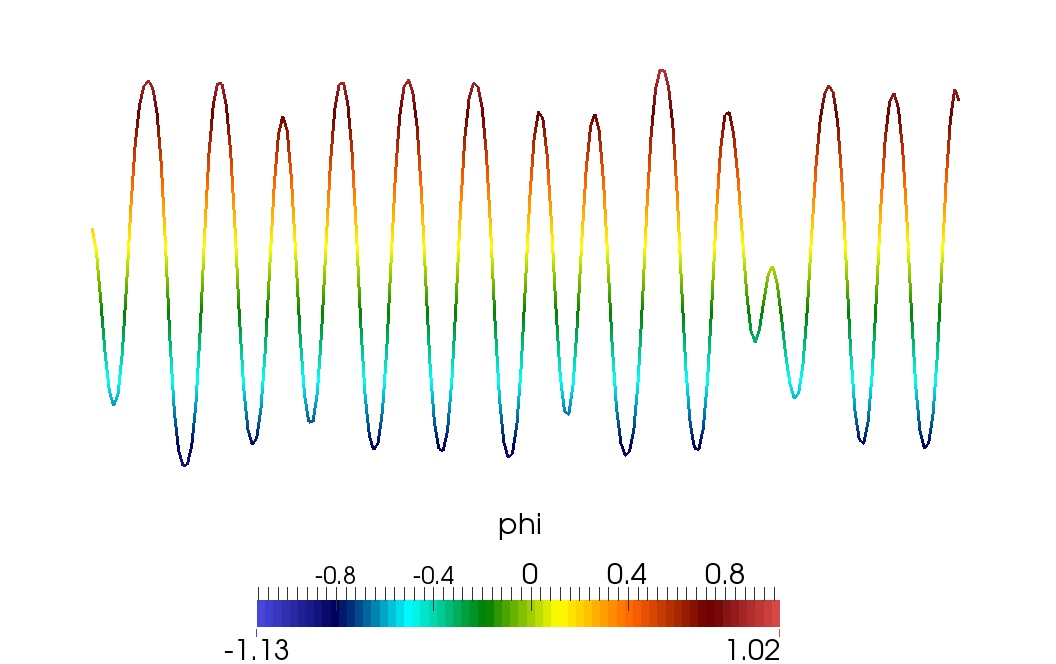}
                             
            }
            \hfill
            \subfigure[][$t=10$]{
              \includegraphics[scale=\figscale, width=0.4\figwidth]
              {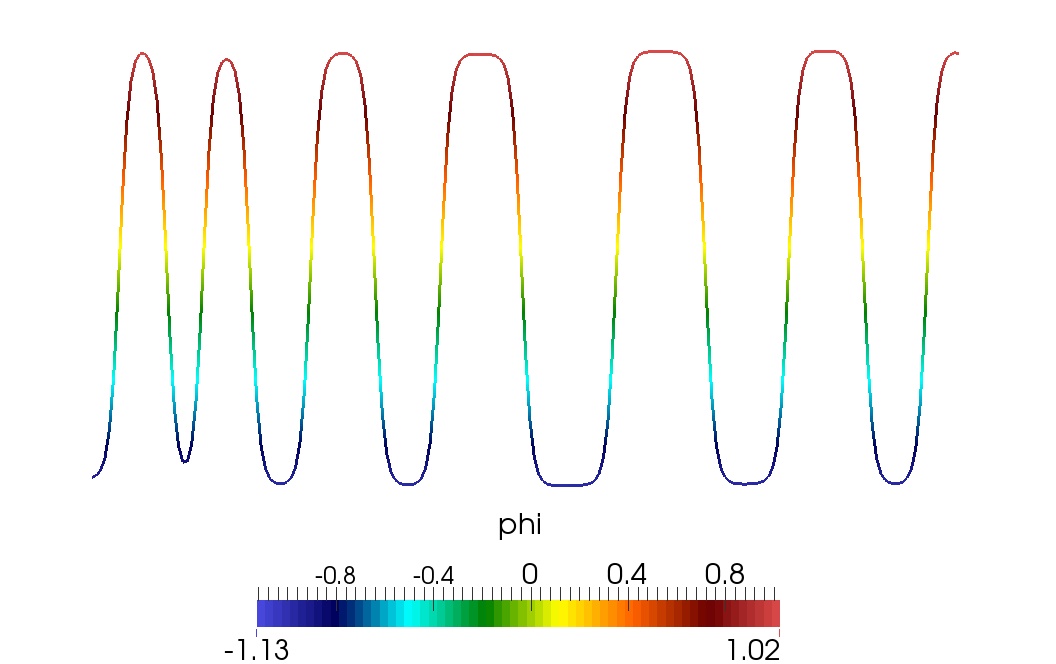}
                             
            }
            \hfill
            \subfigure[][$t=100$]{
              \includegraphics[scale=\figscale, width=0.4\figwidth]
              {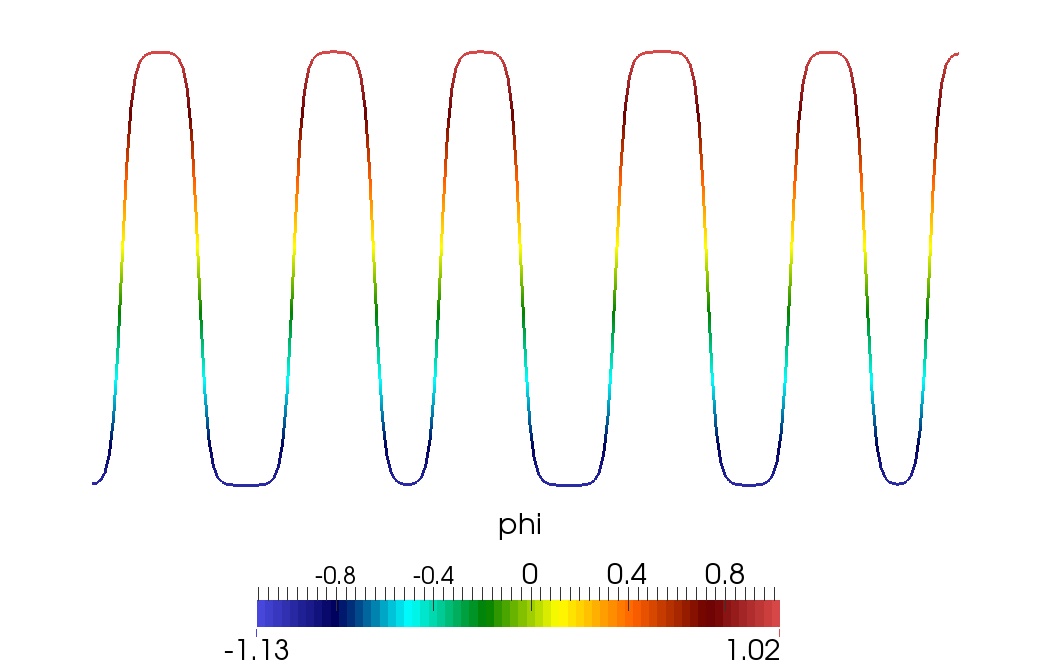}                            
            }
          \end{center}
\end{figure}

\begin{figure}[h!]
  \caption[]
          {
            \label{fig:1d-dens-rat-100} 
            \ref{sec:1d-random} Test 6 -- As Figure
            \ref{fig:1d-dens-rat-2} but in this
            case $\rho_2 / \rho_1 = 100$ and $\max{\phi} = 1.0052$ hence
            $\rho(\phi) > 0$ for all time.}
          \begin{center}
            \subfigure[][$t=0$]{
              \includegraphics[scale=\figscale, width=0.4\figwidth]
                              {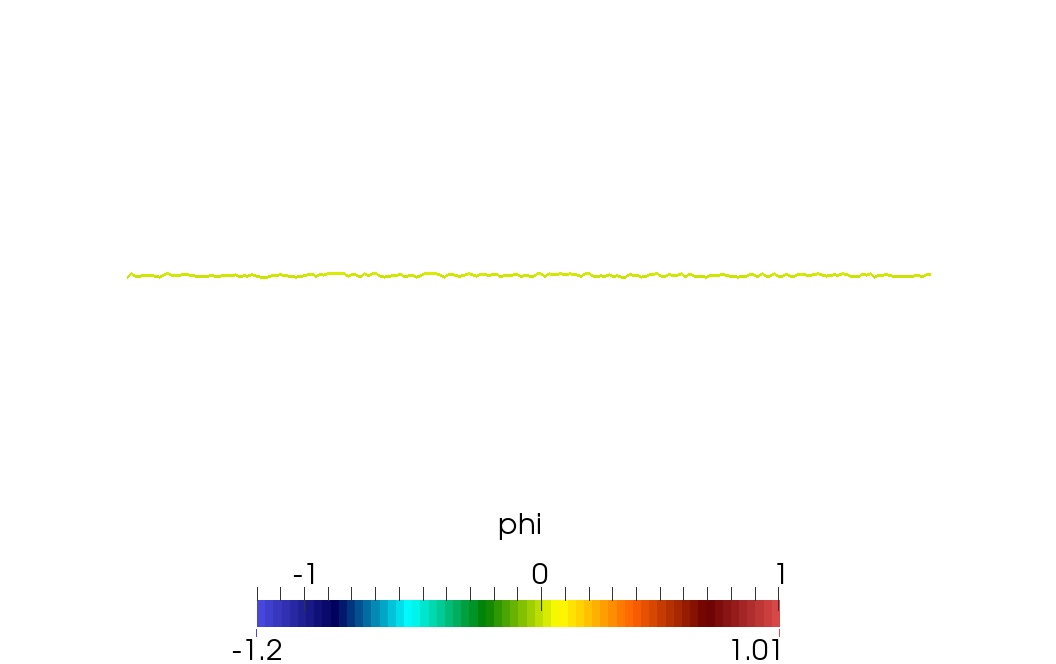}
            }
%
            \hfill
            \subfigure[][$t=0.43$]{
              \includegraphics[scale=\figscale, width=0.4\figwidth]
              {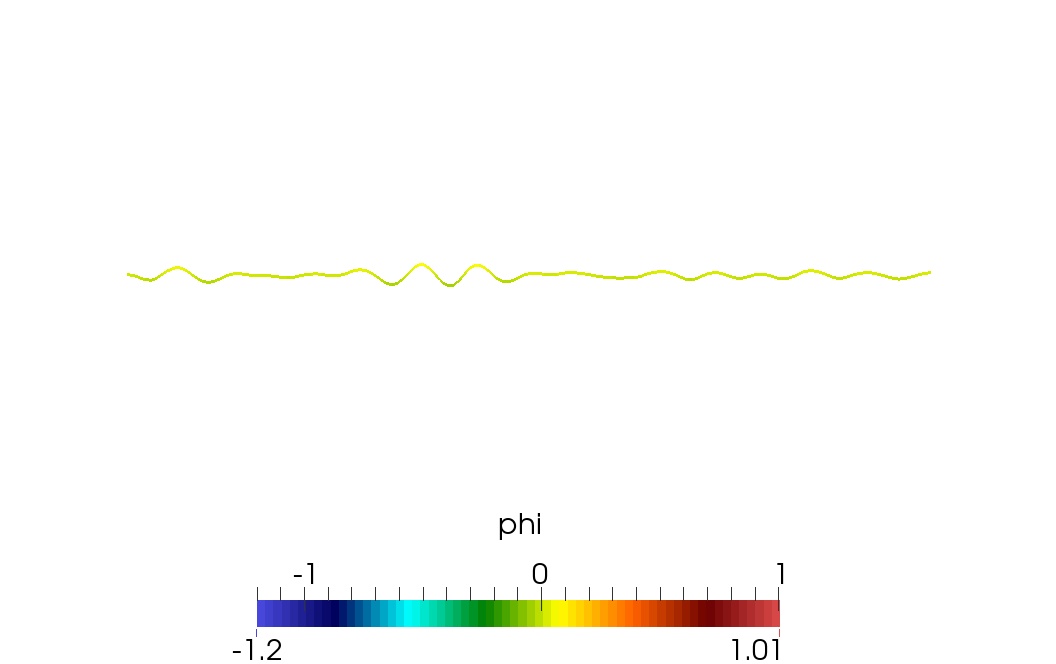}
                             
            }
%
            \hfill
            \subfigure[][$t=10$]{
              \includegraphics[scale=\figscale, width=0.4\figwidth]
              {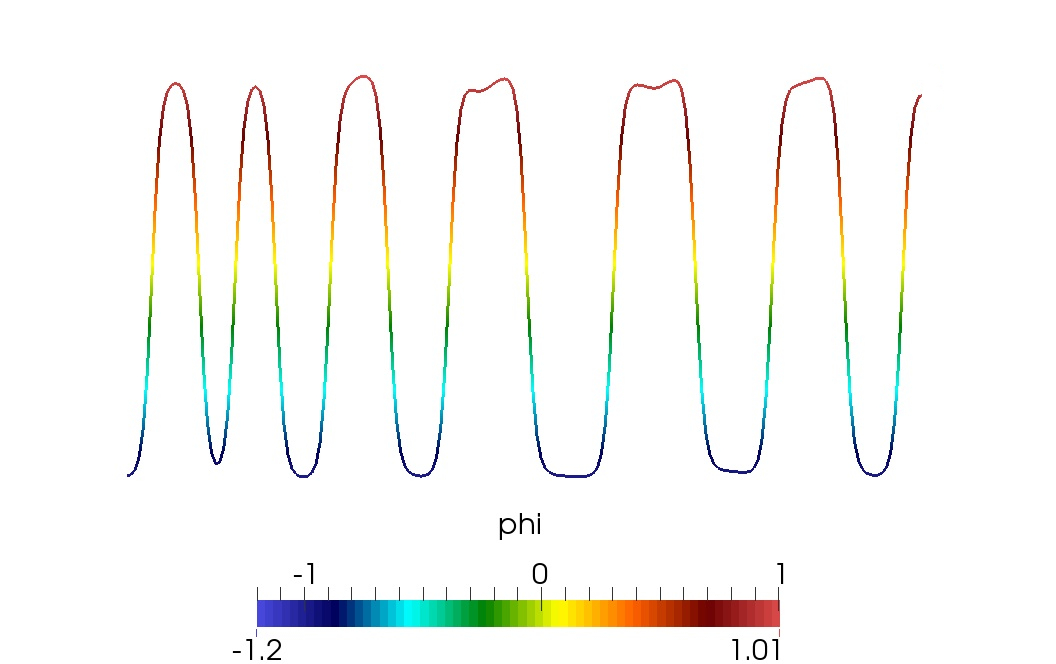}
                             
            }
            \hfill
            \subfigure[][$t=100$]{
              \includegraphics[scale=\figscale, width=0.4\figwidth]
              {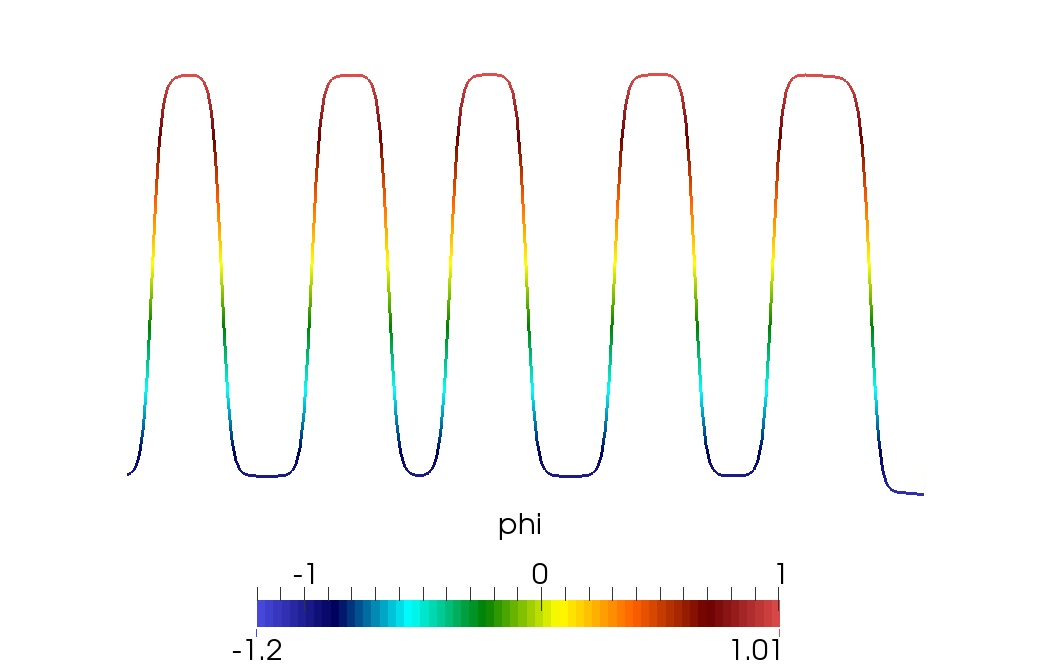}
                             
            }
          \end{center}
\end{figure}

\begin{figure}[h!]
  \caption[]
          {
            \label{fig:1d-dens-rat-1000} 
            \ref{sec:1d-random} Test 6 -- As Figure
            \ref{fig:1d-dens-rat-2} but in this
            case $\rho_2 / \rho_1 = 1000$ and $\max{\phi} = 1.0006$ hence
            $\rho(\phi) > 0$ for all time.}
          \begin{center}
            \subfigure[][$t=0$]{
              \includegraphics[scale=\figscale, width=0.4\figwidth]
                              {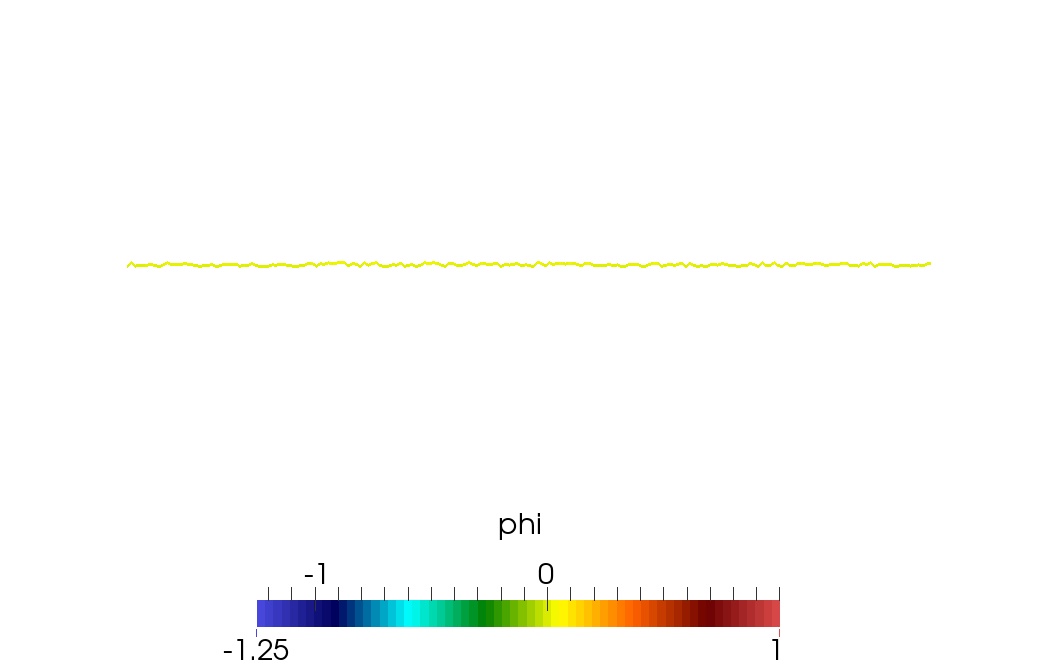}
            }
%
            \hfill
            \subfigure[][$t=0.43$]{
              \includegraphics[scale=\figscale, width=0.4\figwidth]
              {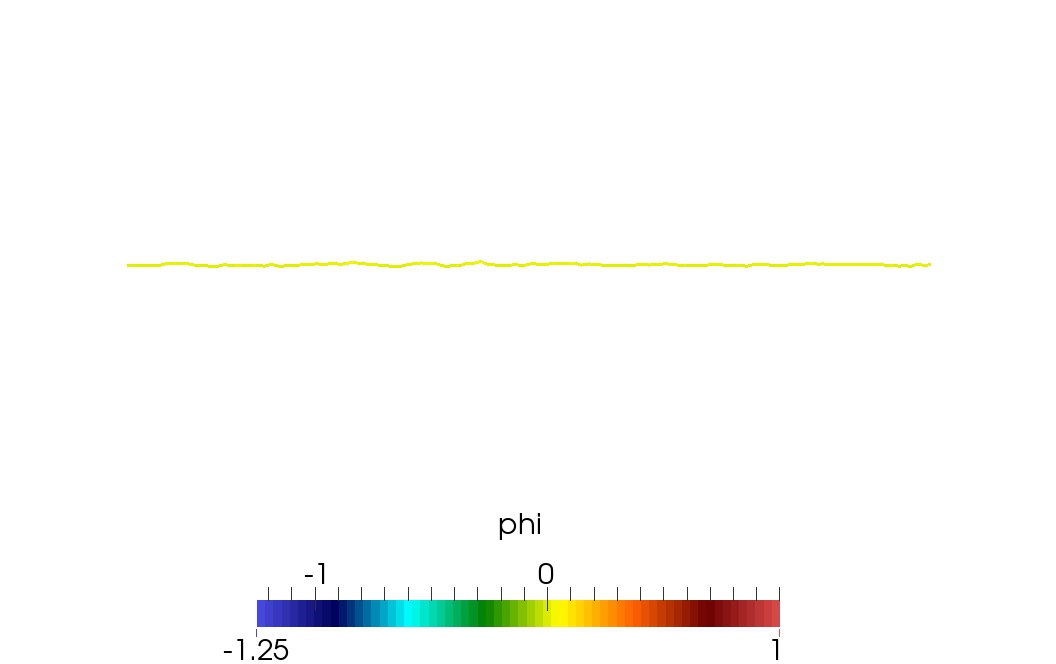}
                             
            }
%
            \hfill
            \subfigure[][$t=10$]{
              \includegraphics[scale=\figscale, width=0.4\figwidth]
              {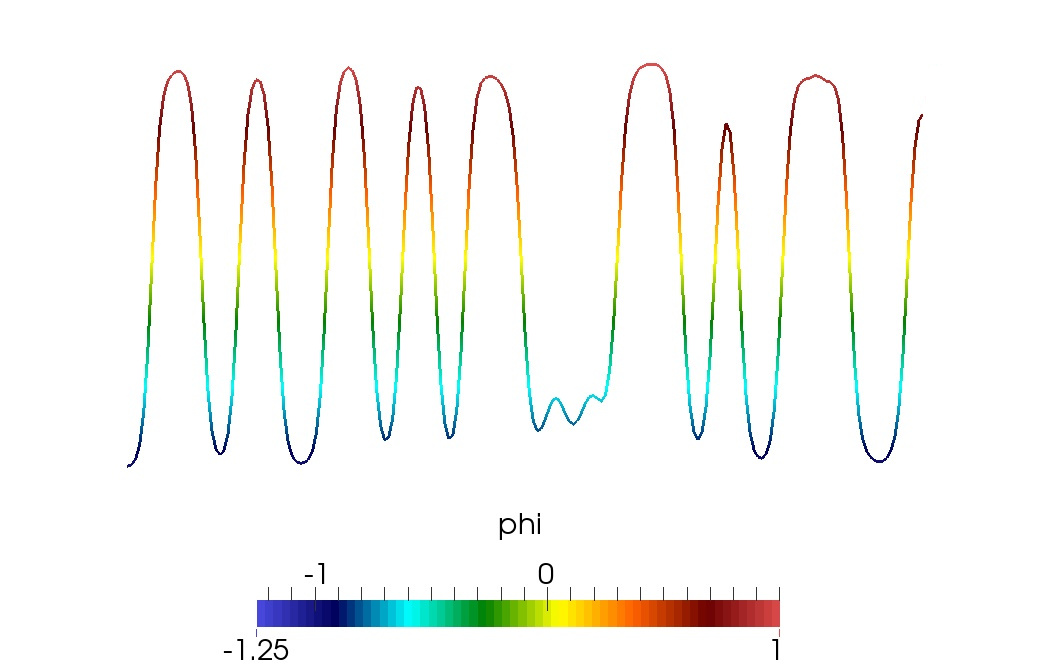}
                             
            }
            \hfill
            \subfigure[][$t=100$]{
              \includegraphics[scale=\figscale, width=0.4\figwidth]
              {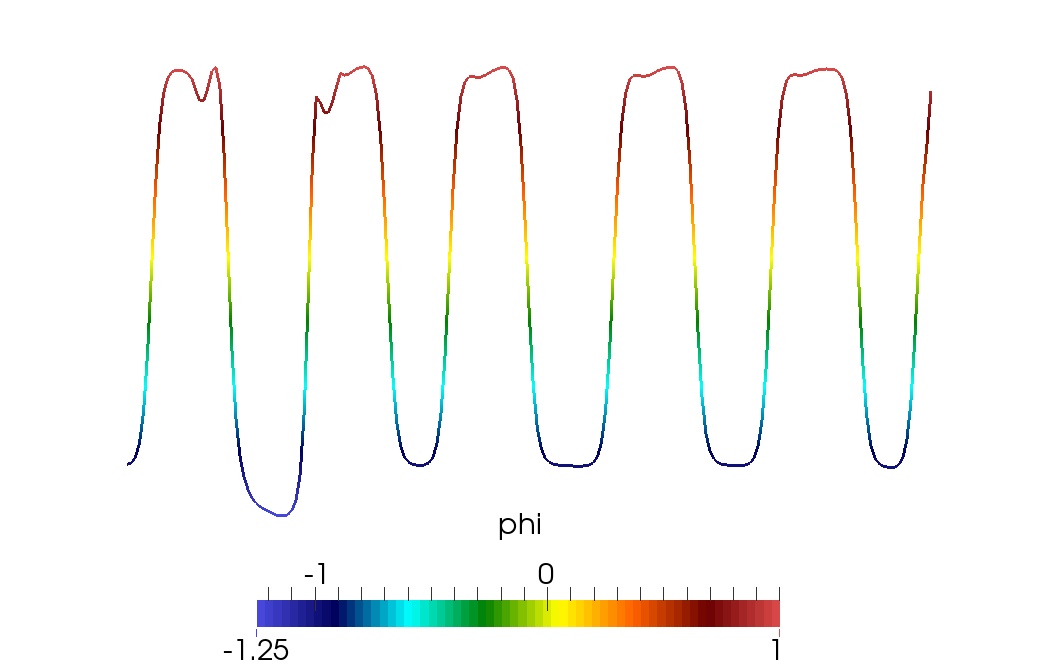}
                             
            }
          \end{center}
\end{figure}

{ In this test we examine the numerical schemes behaviour
  for various density ratios based on the modified energy density
  (\ref{eq:modified-energy}). In Figures
  \ref{fig:1d-dens-rat-2}--\ref{fig:1d-dens-rat-1000} we study a 1D
  equivalent problem to that given in Test 2 for various density
  ratios ranging from $\rho_2 / \rho_1 = 2$ to $\rho_2 / \rho_1 =
  1000$. We note that with $A = \qp{\rho_2 / \rho_1}^2$ the density
  $\rho(\phi) > 0$ and as the density ratio increases the simulation
  takes longer to achieve a steady state.}

\bibliographystyle{alpha}
\bibliography{nskbib,tristansbib,tristanswritings}

\end{document}